\documentclass[10pt,a4paper]{article}
\usepackage[utopia]{mathdesign}

\usepackage{amsfonts}
\usepackage{color}
\usepackage{float}
\usepackage{soul}
\usepackage{graphicx}
\definecolor{labelkey}{rgb}{0,0.08,0.45}
\definecolor{refkey}{rgb}{0,0.6,0.0}
\definecolor{Brown}{rgb}{0.45,0.0,0.05}
\definecolor{lime}{rgb}{0.00,0.8,0.0}
\definecolor{lblue}{rgb}{0.5,0.5,0.99}
\usepackage{amssymb}
\usepackage{theorem}
\usepackage{fancyhdr}
\usepackage{enumerate}
\usepackage[table]{xcolor}
\usepackage{algorithm}
\usepackage{algpseudocode}
\usepackage{amsmath}
\usepackage[normalem]{ulem}
\usepackage[margin=1.1in]{geometry}
\usepackage{hyperref}
\usepackage[colorinlistoftodos]{todonotes}

\definecolor{labelkey}{rgb}{0.6,0.6,0.6}
\definecolor{refkey}{rgb}{0,0.6,0.0}
\def\disp{\displaystyle}
\def\ve{\varepsilon}
\def\dd{\delta}

\def\lm{\lambda}
\def\O{\Omega}
\def\tx{\widetilde{x}}

\def\oa{\bar a}

\def\ox{\bar{x}}
\def\oy{\bar{y}}

\def\ov{\bar{v}}
\def\ou{\bar{u}}

\def\oT{\bar{T}}

\def\tu{\tilde{u}}
\def\oy{\bar{y}}

\def\tT{\Tilde T}

\def\gg{\gamma}
\def\dn{\downarrow}

\def\Lto{\Longrightarrow}

\def\tto{\rightrightarrows}

\def\H{H\"{o}lder}
\def\Limsup{\mathop{{\rm Lim}\,{\rm sup}}}

\def\tto{\rightrightarrows}

\def\Tilde{\widetilde}
\def\Bar{\overline}
\def\bar{\overline}
\def\ra{\right\rangle}
\def\la{\left\langle}

\def\ve{\varepsilon}
\def\B{I\!\!B}
\def\h{\hfill\Box}
\def\R{\mathbb{R}}
\def\N{\mathbb{N}}

\def\co{\mbox{\rm co}\,}
\def\gph{\mbox{\rm gph}\,}
\def\epi{\mbox{\rm epi}\,}

\def\dom{\mbox{\rm dom}\,}

\def\dist{\mbox{\rm dist}\,}

\def\bd{\mbox{\rm bd}\,}

\def\cU{\mathcal{U}}
\def\dn{\downarrow}
\def\O{\Omega}

\def\ph{\varphi}
\def\emp{\emptyset}

\def\oR{\Bar{\R}}
\def\lm{\lambda}
\def\Lm{\Lambda}
\def\gg{\gamma}

\def\dd{\delta}

\def\al{\alpha}
\def\be{\beta}

\def\vt{\vartheta}

\def \N{I\!\!N}
\def\th{\theta}
\def\vTh{\vartheta}

\def\tx{\widetilde{x}}
\def\tu{\widetilde{u}}
\def\ta{\widetilde{a}}

\def\({\left(}
\def\){\right)}
\def\[{\left[}
\def\]{\right]}
\def\n{\left \|}
\def\en{\right \|}
\def\nn {\left \{ }
\def\hnn {\right \}}
\def\var{\mbox{\rm var}}
\def\cA{\mathcal A}
\def\l{\left}
\def\r{\right}
\def\T{\bar T}
\def\H{\bar H}
\def\eq{\begin{equation}}
\def\eeq{\end{equation}}
\def\l{\left}
\def\r{\right}
\newtheorem{theorem}{Theorem}[section]

\newtheorem{proposition}[theorem]{Proposition}
\newtheorem{definition}[theorem]{Definition}
\theoremstyle{plain}{\theorembodyfont{\rmfamily}
}
\theoremstyle{plain}{\theorembodyfont{\rmfamily}
}
\theoremstyle{plain}{\theorembodyfont{\rmfamily}
\theoremstyle{plain}{\theorembodyfont{\rmfamily}
}

\theoremstyle{plain}{\theorembodyfont{\rmfamily}
}

\renewcommand{\theequation}{\thesection.\arabic{equation}}
\newcommand\numberthis{\addtocounter{equation}{1}\tag{\theequation}}
\allowdisplaybreaks

\begin{document}

{\large
\title{Optimal Control of Nonconvex Sweeping Processes with Variable Time via Finite-Difference Approximations}

\author{TAN H. CAO \footnote{Department of Applied Mathematics and Statistics, SUNY (State University of New York) Korea, Incheon, Korea (tan.cao@stonybrook.edu). Research of this author was supported by the National Research Foundation of Korea grant funded by the Korea Government (MIST) NRF-2020R1F1A1A01071015.},\; BORIS S. MORDUKHOVICH\footnote{Department of Mathematics, Wayne State University, Detroit, Michigan, USA (aa1086@wayne.edu). Research of this author was partly supported by the US National Science Foundation under grant DMS-2204519 and by the Australian Research Council under Discovery Projects DP-190100555 and DP-250101112. },\\  DAO NGUYEN\footnote{Department of Mathematics, University of Michigan, Ann Arbor, MI 48109, USA, and Department of Mathematics and Statistics, San Diego State University, San Diego, CA 92182, USA (dnguyen28@sdsu.edu). Research of this author was supported by the AMS-Simon Foundation.},\; TRANG NGUYEN\footnote{Department of Mathematics, Wayne State University, Detroit, Michigan, USA (daitrang.nguyen@wayne.edu). Research of this author was partly supported by the US National Science Foundation under grant DMS-2204519.}, \; NGUYEN N. THIEU\footnote{Institute of Mathematics, Vietnam Academy of Science and Technology, Hanoi, Vietnam, and SUNY Korea, Incheon, Korea (nnthieu@math.ac.vn).  Research of this author was supported by the National Research Foundation of Korea grant funded by the Korea Government (MIST) NRF-2020R1F1A1A01071015.}}
\maketitle
\vspace*{-0.35in}
\begin{abstract}
The paper is devoted to the study of a new class of optimal control problems for nonsmooth dynamical systems governed by nonconvex discontinuous differential inclusions of the sweeping type with involving variable time
into optimization. We develop a novel version of the method of discrete approximations of its own qualitative and numerical importance with establishing its well-posedness and strong convergence to optimal solutions of the controlled sweeping process. Using advanced tools of variational analysis and generalized differentiation leads us to deriving new necessary conditions for optimal solutions to discrete approximation problems, which serve as suboptimality conditions for the original continuous-time controlled sweeping process. The obtained results are applied to a class of motion models of practical interest, where the established necessary conditions are used to investigate the agents' interactions and to develop an algorithm for calculating optimal solutions.
\end{abstract}\vspace*{-0.1in}
\em Keywords:} Nonsmooth dynamical systems, discontinuous differential inclusions, sweeping processes, optimal control, discrete approximations, variational analysis, optimality conditions, motion models \vspace*{-0.03in}

\vspace{-0.2in}

\section{Introduction and Problem Formulation}\label{S1}
\setcounter{equation}{0}\vspace*{-0.1in}

The {\em sweeping process} (``processus de rafle'' in French) was introduced by Moreau in the 1970s (see \cite{moreau}) and then was extensively investigated in the literature. The nonsmooth dynamics of Moreau's sweeping process is described by the discontinuous differential inclusion
\begin{equation}\label{sp}
\left\{\begin{matrix}
\dot{x}(t)\in-N\big(x(t);C(t)\big)\;\textrm{ a.e. }\;t\in[0,T],\\x(0)=x_0\in C(0)\subset\R^n,
\end{matrix}\right.
\end{equation}
where $C(t)$ is a continuously moving convex set, and where the symbol $N(\cdot;C)$ in \eqref{sp} signifies the normal cone of convex analysis defined by
\begin{equation}\label{nc-conv}
N(x;C)=N_C(x):=\left\{\begin{array}{ll}
\big\{v\in\R^n\;\big|\;\la v,u-x\ra\le 0\big\}&\mbox{ if }\;x\in C,\\
\emptyset&\mbox{ otherwise}.
\end{array}\right.
\end{equation}
Over the years, mathematical theory of sweeping processes and related issues for complementarity systems and evolution variational inequalities have been widely developed with numerous applications to models of mechanics, electrical circuits, robotics, economics, traffic equilibria, engineering, etc.; see, e.g., \cite{ant,mb,Bro,brog,CT,hb1,mv0,MM,ve} among many other publications. In a parallel way, dynamical systems of the sweeping type have been independently studied in the Russian literature by Krasnosel'ski\v{i} and his followers in connection with systems of hysteresis; see the book \cite{KP} and the references therein.\vspace{-0.05in}

One of the most fundamental results in the theory of sweeping processes establishes the existence and {\em uniqueness} of solutions to dynamical systems of type \eqref{sp} with convex and mildly nonconvex sets $C(t)$  moving in the Lipschitzian or absolutely continuous way; see \cite{CT,MM}. The latter excludes considering any optimization problem for the sweeping dynamics \eqref{sp}.  This phenomenon for the {\em dissipative discontinuous} differential inclusions of type \eqref{sp} strongly distinguishes them from differential inclusions with {\em Lipschitz continuous} right-hand sides for which the corresponding Cauchy problem admits multiple solutions and  optimal control theory has been largely developed; see, e.g., the books \cite{cl,m-book1,v}.\vspace*{-0.05in}

Formulating optimal control problems for sweeping processes requires entering control actions into \eqref{sp} via some additional terms. In \cite{ET,tolstonogov}, it was done by using additive controlled perturbations of the sweeping dynamics with establishing the existence of optimal controls and developing relaxation procedures of the Bogolyubov-Young type from the calculus of variations. The control formulations in \cite{chhm0,chhm} involved parameterizing polyhedral moving sets by control functions in the form
$C(t)=C(u(t))$ with deriving necessary optimality conditions for the corresponding Mayer and Bolza problems. Yet another formulation of optimal control problems for sweeping processes with establishing necessary optimality conditions was suggested in \cite{bk}, where controlled terms appeared in the adjacent ordinary differential equation; see also \cite{adam} for the further development of this model.\vspace*{-0.05in}

In the recent years, {\em necessary optimality conditions} of different types for various optimal control problems under diverse assumptions have become the dominating theme of sweeping control theory and its numerous applications. We refer the reader to \cite{ac,haddad,brog,CCMN20,CCMN21,CKMNP21a,CMNP21b,CMNN21,cm1,cm2,cm3,cm4,cg21,CMN19,CMN20,cmnn24,CMNNO2024,pfs,pfs23,henrion23,her-pall,khalil,mn20,mnn,vera,zeidan} and the bibliographies therein for a variety of results in this direction. It has been realized from the very beginning \cite{bk,chhm0,chhm} that deriving necessary conditions for optimal solutions to controlled sweeping processes is a {\em highly challenging task}, which is immensely more involved in comparison with standard control theory and optimal control of Lipschitzian differential inclusions. In particular, optimization of the sweeping dynamics is unavoidably accompanied by {\em pointwise state constraints} and exhibits behavior similar to {\em irregular mixed control-state constraints}, which have not been satisfactory resolved even in classical smooth ODE control settings.\vspace*{-0.05in}

Despite significant advances made in the theory of necessary optimality conditions for controlled sweeping processes, many fundamental issues arising from both theory and applications remain unresolved. In particular, some practical models of nanotechnology and traffic equilibria require including the {\em duration} of the sweeping process into optimization, i.e., considering problems with ``free time." To the best of our knowledge, the only paper systematically addressing this issue is the quite recent publication \cite{cmnn24} concerning a Mayer-type  sweeping control problem with the controlled dynamics and uncontrolled {\em polyhedral} moving sets. It should be mentioned that the imposed polyhedrality assumption constitutes a serious restriction from both viewpoints of sweeping control theory and its practical applications.\vspace*{-0.05in}

In this paper, we investigate a new class of constrained Bolza-type sweeping control problems with {\em free time} and  {\em highly nonpolyhedral} (even nonconvex) controlled moving sets as well as with controlled perturbations acting in the sweeping dynamics. Observe that sweeping processes with nonpolyhedral {\em uniformly prox-regular} moving sets were considered in \cite{ant,CT,timoshin,ve}, while without any control and/or optimization. To the best of our knowledge, the only paper addressing the study of optimal control with uniformly prox-regular sets is the one of \cite{cm3}, where a fixed-time sweeping control problem of this type has been investigated in the absence of endpoint constraints on trajectories. However, the class of feasible controls in the perturbed dynamics was chosen there as {\em $W^{1,2}$-smooth}, which is a very restrictive and unnatural requirement never imposed in classical control theory. Such a restrictive requirement came from the drawback of the technique used in \cite{cm3}.\vspace*{-0.05in}

Now we develop a {\em new device} allowing us, in particular, to overcome this drawback and investigate a significantly more general family of sweeping control problems with the class of constrained {\em measurable} controls in additive perturbations of the sweeping dynamics. The free-time sweeping control problem $(P)$ studied in this paper is formulated as follows: minimize the Bolza-type cost functional 
\begin{equation}
\label{e:BP}
J[x,u,a,T]:=\ph\big(x(T),T\big)+\int^T_0\ell\big(t,x(t),u(t),a(t),\dot x(t),\dot u(t)\big) dt
\end{equation} 
over control functions $u(\cdot)$ and $a(\cdot)$ such that the constrained control pair 
$\big(u(\cdot),a(\cdot)\big)\in \cU\times\cA$ generates the corresponding trajectories $x(\cdot)$ of the sweeping differential inclusion 
\begin{equation}
\label{e:SP}
\begin{cases}
-\dot x(t)\in N_{C_{u(t)}}\big(x(t)\big)+f\big(x(t),a(t)\big)\;\;\mbox{a.e.}\;\;t\in[0,T], \\
x(0)=x_0\in C_u(0),\\
(x(T), T)\in \Xi_x\times \Xi_T\subset \R^n\times[0,\infty),  
\end{cases}
\end{equation}
where $\Xi_x$ and $\Xi_T$ are given subsets of $\R^{n}$ and $\R$, respectively, and where the {\it nonconvex} (hence nonpolyhedral) controlled moving sets are defined by
\begin{equation}
\label{e:MS}
C_u:=C+u=\bigcap^m_{i=1}C_i+u,\;\; C_i:=\big\{x\in\R^n\;\big|\;g_i(x)\ge0\big\}\;\;\mbox{for}\;\;i=1,\ldots,m
\end{equation}
 via some convex and twice continuously differentiable functions $g_i:\R^n\to\R$, $i=1,\ldots,m$. The collections of control functions $\cU$ and $\cA$ are described, respectively, as
\eq 
\label{U}
\cU:=\bigcup_{T\ge0}\nn u:[0,T]\to \R^n\;\; \mbox{is absolutely continuous and}\;\; \nu_i(u(t)) \leq 0 \;\;\mbox{a.e.}\;t\in [0,T],\;i=1,\ldots,s \hnn,
\eeq
where $\nu_i(\cdot)$ are continuously differentiable functions with the Jacobian of full rank, and as
\eq
\label{A}
\cA:=\bigcup_{T\ge0}\nn a:[0,T]\to \R^d\;\; \mbox{is measurable and}\;\; a(t)\in A \;\;\mbox{a.e.}\;\; t\in [0,T] \hnn
\eeq
for given nonempty sets $A\subset\R^d$ and natural numbers $m, n, d$, $s$. Here the extended-real-valued terminal cost function $\ph:\R^n\times [0, \infty] \to\Bar\R:=(-\infty,\infty]$ and the running cost/integrand $\ell:[0,T]\times\R^{4n+d}\to\Bar\R$ in \eqref{e:BP} are lower semicontinuous (l.s.c) on their domains. \vspace*{-0.05in}

Since the sets $C_u$ in \eqref{e:SP} are nonconvex, we use therein the appropriate normal cone of nonconvex variational analysis, which is defined and discussed in Section~\ref{S2} and reduces to \eqref{nc-conv} in the case of convex sets. Taking into account that $N_{C_u}(x):=\emp$ if $x\notin C_u$, it follows from \eqref{e:SP} that, besides the pointwise constraints on the controls $u(\cdot)\in{\cal U}$, $a(\cdot)\in{\cal A}$ and the endpoint constraints  in \eqref{e:SP}, the optimal control problem $(P)$ intrinsically contains the pointwise irregular {\em mixed state-control constraints} 
\begin{equation}
\label{e:MC}
g_i\big(x(t)-u(t)\big)\ge 0\;\;\mbox{for all}\;\;t\in[0,T]\;\;\mbox{and}\;\;i=1,\ldots,m,
\end{equation}
which have been realized as the most challenging and largely underinvestigated even for classical optimal control theory for dynamical systems governed by  smooth ODEs.\vspace*{-0.05in}

This paper provides a systematic study of the sweeping control problem $(P)$ under standing general assumptions formulated and discussed in Section~\ref{S2}, which also presents some preliminaries from variational analysis and generalized differentiation broadly used in what follows.\vspace*{-0.05in}

Section~\ref{sec:disc-app} plays a crucial role in our approach. Given {\em any feasible} solutions $(x(\cdot),u(\cdot),a(\cdot),\Bar T)$ to the constrained sweeping dynamics \eqref{e:SP}, we construct there a sequence of {\em finite-difference/discrete approximations} whose solutions, properly extended to the continuous-time interval $[0,\bar T]$, converge to $(x(\cdot),u(\cdot),a(\cdot))$ {\em strongly} in the $W^{1,2}\times W^{1,2}\times L^2$-norm. This constructive result has both theoretical and numerical values while allowing us to basically replace the infinite-dimensional dynamical system \eqref{e:SP} by a finite-dimensional system with the discrete dynamics.\vspace*{-0.05in}

In Section~\ref{sec:exist}, we establish the {\em existence} theorem of optimal solutions to the sweeping control problem $(P)$ under certain convexity assumptions, define a {\em relaxation} of $(P)$ of the Bogolyubov-Young type, and introduce the notions of local minimizers of our subsequent study.\vspace*{-0.05in}

Section~\ref{sec:disc-local} is devoted to constructing a sequence of optimal control problems $(P_k)$ for {\em finite-difference} dynamical systems that admit optimal solutions when $k$ is sufficiently large, while being such that optimal solutions to $(P_k)$, extended to the continuous-time intervals, {\em strongly $W^{1,2}\times W^{1,2}\times L^2$-approximate} the designated {\em local minimizer} for the original sweeping control problem $(P)$. Each problem $(P_k)$ can be equivalently reduced to a (nondynamic) problem of {\em mathematical programming} with increasingly many {\em geometric constraints} defined by {\em graphical sets} that are generated by the normal cone sweeping dynamics. To handle constraints of this type requires explicit evaluations of {\em second-order subdifferential constructions} of variational analysis, which is done in Section~\ref{tools}. \vspace*{-0.05in}

Section~\ref{nec} presents the derivation of {\em necessary optimality conditions} for problems $(P_k)$, which provide {\em suboptimality} conditions of any degree of accuracy for local minimizers of the original problem $(P)$ and thus can be used to solve sweeping control problems. Our device is based on the advanced tools of variational analysis and nonsmooth optimization with applying second-order subdifferential calculations.\vspace*{-0.05in}

The efficiency of the obtained optimality conditions is demonstrated in Section~\ref{sec:motion} for solving a class of practical {\em motion models}, which are described by a nonconvex prox-regular controlled sweeping process with constraints. The established theoretical results allows us investigate {\em possible interactions} between agents of the motion models and to develop a {\em numerical algorithm} to calculate optimal solutions. The final Section~\ref{sec:conc} summarizes the main results of the paper and discuss directions of our future research.\vspace*{-0.05in}

In what follows, we consider trajectories $x: [0,T]\to \R^n$ to the sweeping inclusion \eqref{e:SP} in the space $W^{1,2}([0,T];$ $\R^n)$ generated by control pairs $(u(\cdot),(a\cdot))$ in $W^{1,2}([0,T];\R^n)\times L^2([0,T];\R^d)$ with the endtime variable $T$. To avoid confusions, we identify the trajectory $x:[0,T]\to\R^n$ of the sweeping process with its extension $x_e(\cdot)$ to the interval $[0,\infty]$ by 
$$
x_e(t): = x(T)\;\mbox{ for all }\;t > T.
$$
Given such a trajectory  $x\in W^{1,2}([0,T];\R^n)$. define its the norm 
$$
\n x\en_{W^{1,2}}:= \| x(0)\| + \n \dot x_e\en_{L^2}.  
$$\vspace*{-0.5in}

\section{Preliminaries and Standing Assumptions}\label{S2}
\setcounter{equation}{0}\vspace*{-0.1in}

The notation used in this paper is standard in variational analysis, generalized differentiation, and optimal control; see, e.g., the books \cite{m-book,rw,v}. Unless otherwise stated, the norm $\|x\|$ of $x\in\R^n$ is Euclidean with $\la \cdot, \cdot\ra$ standing for the inner product in $\R^n$. We also denote $\N:=\{1,2,\ldots\}$. Recall that $A^*$ signifies the matrix transposition (adjoint operators).\vspace*{-0.05in}
  
Let $\O\subset\R^n$ be a locally closed set around $\ox\in\O$. The {\em proximal normal cone} to $\O$ at $\ox$ is defined by 
\begin{equation}\label{e1}
N^P_\O(\ox):=\begin{cases}
\big\{v\in\R^n\mid\exists\,\al>0\;\mbox{ with }\;\ox\in\Pi(\ox+\al v;\O)\big\}\quad & \mbox{if}\; \ox\in\O,\\
\emptyset &\mbox{if}\; \ox\not\in\O,
\end{cases}
\end{equation}
where $\Pi(x;\O)$ stands for the Euclidean projector of $x\in\R^n$ onto $\O$, i.e., 
$$
\Pi(x;\O):=\big\{w\in\O\;\big|\;\|x-w\|=\dist(x;\O)\big\},
$$
and where $\dist(x;\O):=\inf_{y\in\O}\|x-y\|$ is the (Euclidean) distance from $x\in\R^n$ to $\O$. Note that $\Pi(x;\O)\ne\emp$ for all $x\in\R^n$, and that the cone \eqref{e1} is always convex but may not be closed in $\R^n$, and it may be trivial $N^P_\O(\ox)=\{0\}$ at boundary points $\ox\in\mbox{bd}\,\O$ of $\O$. Much better properties are exhibited by the {\em limiting normal cone} (known also as the collections of basic, general, or Mordukhovich normals) to $\O$ at $\ox$ defined by
\begin{equation}\label{e1a}
N_\O(\ox):= \big\{v\in \R^n\big|\;\exists x_k \to \ox,\,w_k\in \Pi(x_k;\O),\;\alpha_k\ge 0\;\mbox{ with } \;\al_k(x_k- w_k)\to v\;\mbox{ as }\;k\to\infty\big\}
\end{equation}
if $\ox\in\O$ and $N_\O(\ox):=\emptyset$ if $\ox\notin \O$.\vspace*{-0.05in}

When $\O$ is a convex set, both normal cones \eqref{e1} and \eqref{e1a} reduce to the normal cone of convex analysis \eqref{nc-conv}, while the limiting one \eqref{e1a} may be nonconvex even for simple convex sets as, e.g., for $\O:=\{(x,y)\in\R^2\;|\;y=|x|\}$ at $(0,0)$. Nevertheless, the normal cone \eqref{e1a} is robust and enjoy a {\em full calculus} together with the corresponding generalized differentiation constructions (subdifferential and coderivatives) for nonsmooth functions and set-valued mappings/multifunctions. All of this is due to {\it variational/extremal principles} of variational analysis; see \cite{m-book,m18,rw} for more details and references.\vspace*{-0.05in}

It has been well recognized that the two normal cones in \eqref{e1} and \eqref{e1a} agree for the broader class of nonconvex sets introduced in variational analysis by Poliquin and Rockafellar \cite{pr} under the name of {\em prox-regular} sets. In what follows, we use the {\em uniform} version of prox-regularity introduced below.\vspace*{-0.1in}

\begin{definition}\label{D3} Given $\eta>0$, the set $\O\subset\R^n$ is said to be $\eta$-{\sc prox-regular} if for all $\ox\in\bd\O$ and $v\in N^P_\O(\ox)$ with $\|v\|=1$ we have 
$B(x+\eta v, \eta)\cap\O =\{\ox\}$, where $B(z,r)$ stands for the ball centered at $z$ with radius $r$.
\end{definition}\vspace*{-0.05in}
Note that the $\eta$-prox-regularity of $\O$ can be verified in the following inequality: 
$$
\la v, y-\ox\ra \leq \dfrac{\|v\|}{2\eta}\|y-\ox\|^2 \mbox{ for all } y\in \O, \; \ox\in\bd\O, \mbox{ and } v\in N^P_\O(\ox) . 
$$
In other words,  $\O$ is $\eta$-prox-regular if any external ball with radius smaller than $\eta$ can be rolled around it.\vspace*{-0.05in}

Next we recall the generalized differential constructions for nonsmooth functions and multifunctions associated with the limiting normal cone \eqref{e1a}. Given an extended-real-valued function $\phi: \R^n\to \bar\R$ that is l.s.c.\  around
$\ox\in\dom\phi:=\{x\in\R^n\;|\;\phi(x)<\infty\}$, the {\em subdifferential} of $\phi$ at $\ox$ is defined via the limiting normal cone to the epigraph $\epi\phi:=\{(x,\al)\in\R^{n+1}\;|\;\al\ge\phi(x)\}$ by
\begin{equation}\label{e:sub}
\partial\phi(\ox):=\big\{v\in\R^n\;\big|\;(v,-1)\in N_{\epi\phi}\big(\ox,\phi(\ox)\big)\big\}.
\end{equation}	
Let $F\colon\R^n\tto\R^m$ be a set-valued mapping. Assume that its graph 
$$
\gph F:=\big\{(x,y)\in\R^n\times\R^m\big|\;y\in F(x)\big\}
$$
is closed around $(\ox,\oy)\in\gph F$ and define the {\em coderivative} of $F$ at $(\ox,\oy)$ via the limiting normal cone \eqref{e1a} to the graph of the mapping $F$ by 
\begin{equation}\label{e:cor}
D^*F(\ox,\oy)(u):=\big\{v\in\R^n\;|\;(v,-u)\in N_{\gph F}(\ox,\oy)\;\big\}\;\mbox{for all }\;u\in\R^m.
\end{equation}
When $F\colon\R^n\to\R^m$ is single-valued (in this case, $\oy=F(\ox)$ is omitted in the coderivative notation) and continuously differentiable around $\ox$, we have the representation
$$
D^*F(\ox)(u)=\big\{\nabla F(\ox)^*u\big\},\quad u\in\R^m,
$$
via the adjoint/transposed Jacobian matrix $\nabla F(\ox)^*$. In \cite{m-book,m18,rw}, the reader can find analytic representations of the subdifferential and coderivative with their various properties, calculations, and applications.\vspace*{-0.05in}

Due to the normal cone description of the sweeping dynamics in \eqref{sp}, a significant role in the study of sweeping optimal control is played by the following second-order construction introduced in \cite{m92}. Given $\phi\colon\R^n\to\oR$, the {\em second-order subdifferential} (or generalized Hessian) of $\phi$ at $\ox\in\dom\phi$ relative to $\ov\in\partial\phi(\ox)$ is defined as the coderivative of the first-order subdifferential mappings by
\begin{equation}\label{2nd}
\partial^2\phi(\ox,\ov)(u):=\big(D^*\partial\phi\big)(\ox,\ov)(u),\quad u\in\R^n.
\end{equation}
If $\phi$ is twice continuously differentiable around $\ox$, the generalized second-order subdifferential of $\phi$ in \eqref{2nd} reduces to the (symmetric) Hessian matrix $\partial^2\phi(\ox)(u)=\{\nabla^2\phi(\ox)u\}$ for all 
$u\in\R^n$. We refer the reader to the recent book \cite{m24}
for a comprehensive theory and numerous 
applications of the second-order construction \eqref{2nd} with its explicit calculations for large classes of 
functions appeared in variational analysis and optimization; see also Section~\ref{tools} for the calculation of \eqref{2nd} related to the sweeping dynamics.\vspace*{-0.05in}

Throughout the paper, we impose the following {\it standing assumptions} on the given data of the free-time sweeping control problem $(P)$ ensuring, in particular, that for all {$t\in [0,T]$} the moving constraint set is {\em uniformly prox-regular}, and hence the proximal and limiting normal cones in \eqref{e1} and \eqref{e1a} agree. This enables us to use the normal cone notation ``$N$" without any upper-script in the sweeping differential inclusion \eqref{e:SP} in the rest of the paper and to employ the properties available in the variational analysis for either one of these cones. Our standing assumptions are imposed as follows.\vspace*{-0.05in}

{\bf(H1)} The perturbation mapping $f\colon\R^n\times\R^d\to\R^n$ in \eqref{e:SP} is globally Lipschitzian, i.e., there exists a positive constant $L_f$ such that
\begin{equation}\label{e:Lips}
\|f(x_1,a_1)-f(x_2,a_2)\|\le L_f\(\|x_1-x_2\|+\|a_1-a_2\|\)\;\mbox{ for all }\;(x_1,x_2, a_1,a_2)\in\R^{2n+2 d}.
\end{equation}
Furthermore, there exists a positive constant $M$ ensuring the growth condition
\begin{equation}\label{e:growth-con}
\|f(x,a)\|\le M(1+\|x\|)\;\mbox{ for any }\;x\in\bigcup_{t\in[0,T], T\ge0}C_u(t),\;a\in\R^d.
\end{equation}\vspace*{-0.15in}

{\bf(H2)} The functions $g_i(\cdot),\; i=1,\ldots,m$, satisfy the following conditions:

\hspace{0.2in}{\bf (H2.1)} $g_i(\cdot)$ are twice continuously differentiable on $\R^n$.

\hspace{0.2in}{\bf (H2.2)} There exist positive constants $c$ and $M_i$, $i=1,2,3$, together with open sets $V_i\supset C_i$ {, $i=1,\ldots,m$}, such that $d_H\(C_i,\R^n \backslash V_i \)>c$ and 
\begin{equation}\label{e:bound-grad}
M_1\le\|\nabla g_i(x)\|\le M_2,\quad \nabla^2 g_i(x)\|\le M_3\;\mbox{ for all }\;x\in V_i,{\; i=1,\dots,m},
\end{equation}
where $d_H$ stands for the Hausdorff distance between sets. 

\hspace{0.2in}{\bf (H2.3)} There exist positive numbers $\be$ and $\rho$ such that
\begin{equation}\label{e:w-inverse-triangle-in}
\sum_{i\in I_\rho(x)}\lm_i\|\nabla g_i(x)\|\le\be\Big\|\sum_{i\in I_\rho(x)}\lm_i\nabla g_i(x)\;\Big\|\;\mbox{ for all }\;x\in C\;\mbox{ and }\;\lm_i\ge 0,
\end{equation}
where the index set for the perturbed constraints is defined by
\begin{equation}\label{e:rho-index}
I_\rho(x):=\big\{i\in\{1,\ldots,m\}\big|\;g_i(x)\le\rho\big\}.
\end{equation}\vspace*{-0.2in}

{\bf(H3)} Denote $U:=\disp\bigcap_{i=1}^sU_i:=\bigcap_{i=1}^s\nu^{-1}_i((-\infty,0])$ and suppose that: \vspace*{-0.05in}

\hspace{0.2in}{\bf (H3.1)}  There exist a positive number $L_\nu$  such that 
\eq
\label{nu-Lip}
|\nu_i(u_2)-\nu_i(u_1)| \leq L_\nu\| u_2-u_1\|
\eeq 
for all $u_1$ and $u_2$ in $\R^n$ and $i=1,\ldots, s$. Moreover, $\max_{1\leq i \leq s}\n \nabla\nu_i(u)\en \leq L_\nu$ for all $u\in \R^n$. \vspace*{-0.1in}

\hspace{0.2in}{\bf(H3.2)} If $u\in U $ and  $v_i\in N_{U_i}(u)$  such that $\disp\sum^s_{i=1}v_i=0$, then $v_i=0$ for all $i=1, \ldots, s$.\vspace*{-0.05in}

{\bf(H4)} The set $A$ is a compact subset of $\R^d$, and for each $x\in \R^n$ the set $f(x,A)$ is convex.\vspace*{-0.05in}

{\bf(H5)} The terminal cost $\ph\colon [0,\infty]\times\R^n\to\Bar{\R}$ is l.s.c.\ with respect to the second variable, while the running cost $\ell$ in \eqref{e:BP} is such that $\ell_t:=\ell(t,\cdot)\colon\R^{4n+2d}\to\Bar{\R}$ is l.s.c.\ for a.e.\ $t\in[0,T]$, bounded from below on bounded sets, and the function $t\mapsto\ell(t,x(t),u(t),a(t),\dot x(t),\dot u(t),\dot a(t))$ is summable on $[0,T]$ for each feasible trajectory-control triple $(x(t),u(t),a(t))$.\vspace*{-0.05in}

The following results are taken from \cite[Proposition~2.2]{cm3} and \cite[Proposition~2.3]{cm3}, where the detailed proofs can be found; see also \cite{ve} for more discussions and further developments. The first proposition confirms the aforementioned uniform prox-regularity of the moving set in \eqref{sp}.\vspace*{-0.1in}
\begin{proposition}
\label{Prop}
Under the assumptions in {\rm(H2.2)} and {\rm(H2.3)}, the moving set $C_u(t)$, $t\in[0,T]$, in \eqref{e:SP} and \eqref{e:MS} is $\eta$-prox-regular with $\eta=\frac{\alpha}{M_3\be}$.
\end{proposition}\vspace*{-0.1in}

The next proposition establishes the existence and uniqueness of a feasible sweeping trajectory $x(\cdot)$ corresponding to feasible controls $u(\cdot)$ and $a(\cdot)$ together with solution estimates.\vspace*{-0.1in} 

\begin{proposition}\label{Prop1} Suppose that the assumptions in {\rm(H1)}, and {\rm (H2)} are fulfilled for any fixed $T>0$. Then for arbitrary controls $u(\cdot)\in W^{1,2}([0,T];\R^n)$ and $a(\cdot)\in L^2([0,T];\R^d)$, the sweeping process \eqref{e:SP} admits the unique solution $x(\cdot)\in W^{1,2}([0,T];\R^n)$ generated by these controls. Moreover, we have the estimates
$$
\|x(t)\|\le l:=\|x_0\|+e^{2MT}\left(2MT(1+\|x_0\|)+\int^T_0\|\dot{u}(s)\|ds\right)\;\mbox{ for all }\;t\in[0,T],
$$
\begin{equation}\label{e:boundedness}
\|\dot{x}(t)\|\le 2(1+l)M+\|\dot{u}(t)\|\;\mbox{ a.e. }\;t\in[0,T].
\end{equation}
\end{proposition}
\vspace*{-0.25in}

\section{Discrete Approximations of Constrained Sweeping Dynamics}\label{sec:disc-app}
\setcounter{equation}{0}\vspace*{-0.1in}

This section is devoted to constructing 
finite-difference/discrete approximation of the constrained sweeping process\eqref{e:SP} without considering its optimization.\vspace*{-0.05in}

Let us first rewrite the controlled sweeping differential inclusion \eqref{e:SP} in an equivalent differential inclusion form for extended trajectories. Define the set-valued mapping $F:\R^n\times\R^n\times\R^d\tto\R^n$ by 
\begin{equation}
\label{F}
F(x,u,a):=N_C(x-u)-f(x,a).
\end{equation}
It follows from Motzkin's theorem of the alternative that $F$ admits the explicit representation 
\begin{equation}\label{e:F-rep}
F(x,u,a)=-\left\{\sum_{i\in I(x-u)}\lm_i\nabla g_i(x-u)\; \Big| \; \lm_i\ge 0\right\}-f(x,a)
\end{equation}
via the index set of active constraints
\begin{equation}\label{e:a-index}
I(y):=\big\{i\in\{1,\ldots,m\}\;\big|\; g_i(y)=0\big\}
\end{equation}
at the point $y:=x-u\in C$. Then the sweeping differential inclusion \eqref{e:SP} can be rewritten as
\begin{equation}
\label{e:SP1}
\big(-\dot x(t),u(t),a(t)\big) \in F(x(t),u(t),a(t)) \times U \times A\quad\mbox{a.e.}\; t\in [0,T],
\end{equation}
with the initial-boundary conditions 
$$
\begin{cases}
(x_0,u(0),a(0))\in C_u(0)\times U \times A, \\
(x(T), T) \in \Xi_x\times\Xi_T.  
\end{cases}
$$ 
For each $k\in\N$, consider a real number $T^k$ approximating $T$ and the {\em uniform discrete partition/mesh} on the interval $[0,T^k]$ defined by
\begin{equation} \label{e:DP}
\Delta_k(T^k)\colon=\left\{0=t^k_0<t^k_1<\ldots<t^k_{k-1}<t^k_{k}=T^k\right\}\;\mbox{ with }\;h_k:=t^k_{j+1}-t^k_j =  \frac{T^k}{{k}}, \quad j=0,\ldots,k-1.
\end{equation}\vspace*{-0.1in}
  
Fix now a feasible solution $(x(\cdot), u(\cdot), a(\cdot), \T)$ to \eqref{e:SP1} such that the functions $a(\cdot),\dot u(\cdot), \dot x(\cdot)$ are of {\em bounded variation} (BV) on $[0,\Bar T]$. Given the partition \eqref{e:DP} with $T^k=\T$, the following major theorem  establishes the existence of a {\em strong discrete approximation} of the solution in the space $W^{1,2}([0,\T];\R^{2n})\times L^2([0,\T];\R^d)$. Its proof provides the {\em iteration scheme} for constructing this approximation. \vspace*{-0.1in}

\begin{theorem}\label{Th1}
Under the assumptions in {\rm (H1)}--{\rm (H4)} be satisfied, let $(x(\cdot),u(\cdot)$, $a(\cdot),\T)$ be a feasible solution to the sweeping inclusion in \eqref{e:SP} for which we have the aforementioned BV properties together with estimates
\begin{equation}
\label{e:var}
\begin{cases}
\max\nn \var\(\dot x(\cdot);[0,\T]\),\var\(\dot u(\cdot);[0,\T]\), \var\( a(\cdot);[0,\T]\)\hnn \le \mu, \\[2ex]
\max\l\{\l\| \frac{ u(t^k_1)-u(0)}{h_k}\r\|, \l\| \frac{ u(\T)-u(t^k_{k-1})}{h_k}\r\| \r\}\leq \mu
\end{cases}
\end{equation}
with some constant $\mu>0$. Then there exist sequence of piecewise linear functions {$\big\{\big(x^k(t^k_j),u^k(t^k_j)\big) |$ $j=0,\ldots k\big\}$} and piecewise constant functions { $\{a_k(t^k_j)|\; j=0,\ldots k\}$} to discrete inclusions 
\begin{equation}
\label{e:ua-dc}
\begin{cases}
u^k(t^k_j)\in U_k:=\bigcap^s_{i=1}U^i_k:=\bigcap^s_{i=1}\nu^{-1}_i\((-\infty, L_\nu\delta_k)\),\\
{a^k(t^k_j)\in A}, \\ 
(x^k(\T), \T)\in (\Xi_x+\mu^x_kB(0,1))\times\Xi_T,
\end{cases}
\end{equation}
\begin{equation}
\label{e:x-dc}
x^k(t)=x^k(t_j) + (t-t_j)v^k_j,\;\;t^k_j\le t\le t_{j+1}\;\;\mbox{with}\;\; -v^k_j\in F(x^k(t_j),u^k(t_j),a^k(t_j))
\end{equation}
for $j=0,\ldots,k-1$ with the endpoint conditions 
\begin{equation}
\label{ec}
g_i\(x^k(0)-u^k(0)\)\ge0\;\;\mbox{for all}\;\;i=1,\ldots,m,
\end{equation}
such that $\nn(x^k(\cdot),u^k(\cdot))\hnn\to (x(\cdot),u(\cdot))$ uniformly on $[0, \T]$ and that $\nn(\dot x^k(\cdot),\dot u^k(\cdot), a^k(\cdot))\hnn\to (\dot x(\cdot),$ $\dot u(\cdot), a(\cdot))$ in the $L^2$-norm, where $\nn\(\delta_k,\mu^x_k\)\hnn$ is a sequence of positive numbers converging to zero as $k\to\infty$ ,and where $F$ is defined in \eqref{F}. Moreover, for every $k\in \N$ we have the estimate 
\begin{equation}
\begin{array}{ll}
\label{uk-est}
\mbox{\rm var}\l(\dot u^k;[0,\T] \r) \leq \Tilde\mu \quad \mbox{and} \quad \max\l\{\l\| \frac{ u^k(t^k_1)-u^k(0)}{h_k}\r\|, \l\| \frac{ u^k(\T)-u^k(t^k_{k-1})}{h_k}\r\|   \r\}\leq \Tilde\mu,
\end{array} 
\end{equation}
where $\Tilde\mu$ is a constant defined as
\begin{equation}
\label{tmu}
{\Tilde\mu\colon =\max\{3\mu+1+4L_f\T\mu e^{L_f\T}(2L_f+1) + 2L_f\mu, \; 2e^{L_f\T }(2L_f+1)(\mu+1) + \mu\}.}
\end{equation}
\end{theorem}\vspace*{-0.1in}
{\bf Proof}. 
For each $k\in\N$, consider the uniform partition $\Delta_k(\T)$ in \eqref{e:DP}. For convenience, we split the proof into the following major steps and skip the superscript ``$k$" for the mesh points $t^k_j$.\\[0.5ex]
\noindent
{\bf Step~1:} {\it Constructing $a^k(\cdot)$ to approximate $a(\cdot)$}. Since $a(t)\in A$ a.e., there exists $\tau_j\in [t_j,t_{j+1})$ such that $a(\tau_j)\in A$ for all $j=0,\dots k-1$. Let us construct $a^k:[0,T]\to \R^d$ by
\begin{equation*}
a^k(t):=
\begin{cases}
a(\tau_j) \quad&\mbox{if} \; t\in [t_j,t_{j+1}),\\
a(\tau_{k-1}) &\mbox{if} \;t= \T.
\end{cases}
\end{equation*}
For each $j=0,\dots,k-1$, find  $s^a_{\tau_j}\in [t_j,t_{j+1}]$ satisfying 
$$
\sup_{t\in [t_j,t_{j+1}]}\| a(t)-a(\tau_j)\| \leq \|a(s^a_{\tau_j})-a(\tau_j)\|+2^{-k}.
$$
This brings us to the relationships
\begin{equation*}
\begin{array}{ll}
&\disp\int^{\T}_0\n a(t)-a^k(t) \en^2 dt = \sum_{j=0}^{k-1} \int_{t_j}^{t_{j+1}} \n a(t)-a(\tau_j)\en^2 dt\disp\leq  \sum_{j=0}^{k-1} \int_{t_j}^{t_{j+1}} \l[\|a(s^a_{\tau_j})-a(\tau_j)\|+2^{-k}\r]^2 dt\\
&\disp\leq\sum_{j=0}^{k-1} 2 h_k \l(\|a(s^a_{\tau_j})-a(\tau_j)\|^2 +4^{-k}\r)= 2 h_k\l(\sum_{j=0}^{k-1}\|a(s^a_{\tau_j})-a(\tau_j)\|^2 +k4^{-k}\r)\\
&\leq \dfrac{2\T \, \var^2\( a(\cdot);[0,\T]\)}{k}+2^{-2k+1}\T\leq \dfrac{2\T}{k}\mu^2 + 2^{-2k+1}\T,
\end{array}
\end{equation*}
which ensure in turn that
\eq
\label{mu-k}
\begin{array}{ll}
\mu_k^a:={\disp \int^{\T}_0}\n a(t)-a^k(t)\en^2 dt \leq \dfrac{2\T}{k}\mu^2 + 2^{-2k+1}\T,\quad k\in\N.
\end{array}
\eeq
Therefore, $\mu_k^a \to 0$ and the sequence $\{a^k(\cdot)\}$ converges strongly to $a(\cdot)$ in $L^2([0,T];\R^d)$ as $k\to\infty$.\\[0.5ex]
\noindent
{\bf Step~2:} {\it Constructing $\(x^k(\cdot),u^k(\cdot)\)$ to approximate $\(x(\cdot),u(\cdot)\)$.} Let us start with denoting
\begin{equation*}
\begin{cases}
x^k_0: =x_0,\\
u^k_0: = x^k_0 - x(t_0) +u(t_0),\\
v^k_0:= \Pi\(-\dot x(t_0);F\(t_0,x_0,u_0,a_0\)\)		
\end{cases}
\end{equation*}
and call the collection $\{v^k_J\;|\;j=0,\ldots,k-1\}$ a discrete velocity. Then we proceed with the following iterations: for $j=0,\ldots,k$, consider the vectors
\begin{equation*}
\begin{cases}
u^k_j:=x^k_j-x(t_j)+u(t_j),\\
v^k_j:=\Pi\(-\dot x(t_j);F\(t_j,x^k_j,u^k_j,a^k_j\)\),\\
x^k_{j+1} := x^k_j+h_kv^k_j
\end{cases}
\end{equation*}
and define the approximate trajectory and the approximate control by
\begin{equation*}
x^k(t) =\begin{cases}
x^k_j+(t-t_j)v^k_j &\mbox{for}\;\;t\in [t_j,t_{j+1}),\quad j=1,\dots,k-1,\\
x^k_k+h_kv^k_k&\mbox{for }\; t=\T,
\end{cases}
\end{equation*}
$$
u^k(t)=\begin{cases}
u^k_j+\frac{(t-t_j)}{h_k}\(u^k_{j+1}-u^k_j\)&\mbox{ for }\; t\in [t_j,t_{j+1}),\quad j=1,\dots,k-1,\\
u^k_{k}&\mbox{for }\;t=\T.
\end{cases}
$$
It follows from the above constuctions that we have the relationship
$$
F\big(x^k_j,u^k_j,a^k_j\big)=F\big(x(t_j),u(t_j),a(t_j)\big)+f\big(x^k_j,a^k_j\big)-f\big(x(t_j),a(t_j)\big),
$$
which yields the estimate of the difference between the discrete velocity $v^k_j$ and the continuous one $\dot x(t_j)$:
\begin{equation}\label{dv1}
\begin{aligned}
\n v^k_j-\dot x(t_j) \en &= \dist\(-\dot x(t_j) ;F\(x^k_j,u^k_j,a^k_j\)\) \le \n f(x^k_j,a^k_j)-f(x(t_j),a(t_j)) \en\\
& \le L_f\n x^k_j-x(t_j)\en +L_f\n a^k_j-a(t_j)\en,
\end{aligned}
\end{equation}
where $L_f$ is taken from {\rm (H1)}. As a consequence, we get 
\eq
\label{x-est}
\begin{aligned}
&\n x^k_{j+1} -x(t_{j+1})\en = \n x^k_j + h_kv^k_j-x(t_j)-\int^{t_{j+1}}_{t_j}\dot x(s)ds\en\\
&\le \n x^k_j-x(t_j)\en + \int^{t_{j+1}}_{t_j}\n v^k_j-\dot x(s)\en ds\\
&\le  \n x^k_j-x(t_j)\en + \int^{t_{j+1}}_{t_j}\n v^k_j-\dot x(t_j)\en ds + \int^{t_{j+1}}_{t_j}\n \dot x(t_j)-\dot x(s)\en ds\\
&\le \(1+h_kL_f\)\|x^k_j-x(t_j)\|+L_f \sum^{k-1}_{j=0}\|a^k_j-a(t_j)\|+\int^{t_{j+1}}_{t_j}f^x_j(s)ds\\
&=(1+h_k L_f)\|x^k_j-x(t_j)\|+L_f\int^{t_{j+1}}_{t_j} \n a^k(s)-a(s)\en ds +
\int^{t_{j+1}}_{t_j}f^x_j(s)ds\\
&\le(1+h_k L_f)\|x^k_j-x(t_j)\|+L_f\int^{t_{j+1}}_{t_j}\|a^k_j-a(t_j)\|ds+
\int^{t_{j+1}}_{t_j}f^x_j(s)ds\\
&\le(1+h_k L_f)\|x^k_j-x(t_j)\|+L_f\int^{t_{j+1}}_{t_j}\|a^k(s)-a(s)\|ds+\int^{t_{j+1}}_{t_j}\big[L_f f^a_j(s)+f^x(s)\big]ds
\end{aligned}
\eeq
where $f^x_j(s):=\|\dot{x}(t_j)-\dot{x}(s)\|$ and $f^a_j(s):=\|a(t_j)-a(s)\|$, $j=0,\ldots,k$. Let us further denote
$$
\left\{\begin{array}{ll}
\Lambda:=1+h_k L_f,\\
\gamma_j:=\|x^k_j-x(t_j)\|,\\
\lambda_j:= L_f\int^{t_{j+1}}_{t_j}\|a^k(s)-a(s)\|ds+\int^{t_{j+1}}_{t_j}\big[L_f f^a_j(s)-f^x_j(s)\big]ds
\end{array}\right.
$$
for $j=0,\ldots,k-1$. Then \eqref{x-est} reads as
\begin{equation*}
\gamma_{j+1}\le\Lambda\gamma_j+\lambda_j,\quad j=0,\ldots,k-1,
\end{equation*}
which in turn implies that 
\begin{equation*}
\gamma_j\le\Lambda^{j-1}\lambda_0+\Lambda^{j-2}\lambda_1+\ldots+\Lambda^0\lambda_{j-1},\;\mbox{ where }\;\Lambda^j:=(1+h_k L_f)^j\le e^{L_f\T}.
\end{equation*}
In this way, we arrive at the estimate
\begin{equation}\label{gamma-est}
\gamma_j\le e^{L_f\T}\sum_{j=0}^{k-1}\lm_j.
\end{equation}
Moreover, it follows from the definitions that
$$
\sum_{j=0}^{k-1}\lm_j=\sum^{k-1}_{j=0}\Big[L_f\int^{t_{j+1}}_{t_j}\|a^k(s)-a(s)\|ds+\int^{t_{j+1}}_{t_j}\big[L_f f^a_j(s)+f^{x}_j(s)\big]ds\Big],
$$
which readily yields the relationships
\eq\label{a-est}
\begin{aligned}
&\sum^{k-1}_{j=0}\int^{t_{j+1}}_{t_j} \n a^k(s)-a(s)\en ds= \sum^{k-1}_{j=0}\int^{t_{j+1}}_{t_j}  \n a(\tau_j)-a(s)\en ds\\
&\leq  \sum^{k-1}_{j=0}\int^{t_{j+1}}_{t_j}\[ \n a(\tau_j)-a(s^a_{\tau_j})\en+2^{-k}\] ds= \sum^{k-1}_{j=0}h_k\n a(\tau_j)-a(s^a_{\tau_j})\en +\T2^{-k}\\
&\leq h_k\var\(a(\cdot);[0,\T] \) + \T2^{-k}\leq h_k\mu + \T2^{-k}.
\end{aligned}
\eeq
To proceed further, pick $s^x_j$ and $s^a_j$ from the subintervals $[t_j,t_{j+1}]$ such that 
\eq
\label{BV}
\begin{cases}
\disp\sup_{s\in[t_j,t_{j+1}]}f^x_j(s)\le \n \dot x(t_j)-\dot x(s^x_j)\en + 2^{-k}\\
\disp\sup_{s\in[t_j,t_{j+1}]}f^a_j(s)\le \n a(t_j)- a(s^a_j)\en + 2^{-k}.
\end{cases}
\eeq
Using \eqref{e:var}, \eqref{BV} and employing the same arguments as above give us
\eq
\label{a-est1}
\begin{aligned}
\sum^{k-1}_{j=0}\int^{t_{j+1}}_{t_j} f^a_j(s)ds &\le h_k\sum^{k-1}_{j=0}\( \n a(t_j)- a(s^a_j)\en + \n a(s^a_j) - a(t_{j+1})\en+2^{-k}\)\\
&\le h_k\var\(a;[0,\T]\) +h_k k 2^{-k}\le h_k\mu+\T2^{-k},
\end{aligned}
\eeq
\eq
\label{x-est1}
\sum^{k-1}_{j=0}\int^{t_{j+1}}_{t_j} f^x_j(s)ds \le  h_k\mu+\T2^{-k}.
\eeq
Then it follows from ~\eqref{gamma-est},~\eqref{a-est},~\eqref{a-est1}, and~\eqref{x-est1} that 
\begin{equation}
\label{x-est1a}
\n x^k_j- x(t_j)\en =\gg_j \le \delta_k:=e^{L_f\T}\(2L_f+1\)\( h_k\mu+\T2^{-k}\),
\end{equation}
and therefore we arrive at the relationships
$$
\begin{aligned}
h_k\n v^k_j-\dot x(t_j)\en&\le h_k L_f\n x^k_j-x(t_j)\en +h_kL_f\n a^k_j-a(t_j)\en\\
&\le h_kL_f\delta_k+ L_f\int^{t_{j+1}}_{t_j}  \n a^k(s)-a(s)\en ds +L_f\int^{t_{j+1}}_{t_j}f^a_j(s)ds\\
&\le h_kL_f\delta_k+2L_f\(h_k\mu+\T 2^{-k}\)= L_f\[h_k\delta_k+2\(h_k\mu+\T 2^{-k}\) \]
\end{aligned}
$$
for $j=0,\ldots,k$, which allows us to estimate the error of using $x^k(t)$ to approximate $x(t)$ by
\begin{equation}
\label{x-est1b}
\begin{aligned}
\n x^k(t)-x(t)\en& \le \n x^k_j-x(t_j)\en + \int^t_{t_j}\n v^k_j-\dot x(s)\en ds\\
&\le \delta_k+h_k\n v^k_j-\dot x(t_j)\en+ \int^{t_{j+1}}_{t_j}f^x_j(s)ds\\
&\le \mu^x_k:=(h_k\mu+\T 2^{-k})(2L_f+1)\[e^{L_f\T}(1+L_fh_k)+1\]
\end{aligned}
\end{equation}
for any $t\in[t_j,t_{j+1})$ and all indices $j=0,\ldots,k-1$. This justifies the uniform convergence and hence the strong convergence in $L^2([0,\T];\R^n)$ of $\nn x^k(\cdot)\hnn$ to $x(\cdot)$ together with the claimed inclusion 
$$
(x^k(\T), \T)\in (\Xi_x+\mu^x_kB(0,1))\times\Xi_T.
$$
{\bf Step~3:} {\it Verifying the strong convergence of $\nn \dot x^k(\cdot)\hnn$ to $\dot x(\cdot)$ in $L^2\([0,\T];\R^n\)$.} We clearly have
$$\vspace*{-0.5ex}
\begin{aligned}
\int^{\T}_0 \n\dot x^k(s) -\dot x(s) \en^2 ds &=\sum^{k-1}_{j=0}\int^{t_{j+1}}_{t_j}\n v^k_j-\dot x(s)\en^2 ds \\ 
&\le 2\sum^{k-1}_{j=0}\[ h_k\n v^k_j-\dot x(t_j)\en^2 +  \int^{t_{j+1}}_{t_j}\[f^x_j(s)\]^2ds\],
\end{aligned}
$$
which readily yields the following inequalities:
$$
\begin{aligned}
h_k\n v^k_j-\dot x(t_j)\en^2 &\le L^2_fh_k\[\n x^k_j-x(t_j)\en +\n a^k_j-a(t_j)\en\]^2\\
&\displaystyle\le 2 L^2_f h_k\delta^2_k +2 L^2_f \int^{t_{j+1}}_{t_j} \[ \n a^k(s)-a(s)\en+f^a_j(s)\]^2 ds\\
&\displaystyle\le 2 L^2_f h_k\delta^2_k +4 L^2_f\int^{t_{j+1}}_{t_j}\n a^k(s)-a(s)\en^2 ds +4 L^2_f\int^{t_{j+1}}_{t_j}\[f^a_j(s)\]^2 ds,
\end{aligned}
$$
$$
\begin{aligned}
&\sum^{k-1}_{j=0}\int^{t_{j+1}}_{t_j}\[f^a_j(s)\]^2 ds\le \sum^{k-1}_{j=0}\int^{t_{j+1}}_{t_j}\[ f^a_j(s^a_j)+2^{-k}\]^2 ds = \sum^{k-1}_{j=0}h_k\[ f^a_j(s^a_j)+2^{-k}\]^2\\
&\le 2\sum^{k-1}_{j=0}h_k \nn\[f^a_j(s^a_j)\]^2+4^{-k} \hnn  \le 2 h_k\[\sum^{k-1}_{j=0}f^a_j(s^a_j)\]^2+2 h_k k 4^{-k}\\
&\le 2 h_k\[\sum^{k-1}_{j=0}\(\n a(t_j)- a(s^a_j)\en + \n a(s^a_j) - a(t_{j+1})\en\) \]^2+2\T4^{-k}
\le 2 h_k\var^2\(a(\cdot)\) + 2\T 4^{-k}\le 2 h_k\mu^2+2\T4^{-k}.
\end{aligned}
$$
Using the same arguments as above brings us to
$$
\sum^{k-1}_{j=0}\int^{t_{j+1}}_{t_j}\[f^x_j(s)\]^2 ds\le 2 h_k\mu^2+2\T 4^{-k},
$$
and therefore we arrive at the estimates 
$$
\begin{aligned}
&\int^{\T}_0 \n\dot x^k(s) -\dot x(s) \en^2 ds\le 4L^2_fh_kk\delta^2_k+8L^2_f\int^{\T}_0\n a^k(s)-a(s)\en^2 ds +4\(2h_k\mu^2+2\T 4^{-k}\)\\
&+8 L^2_f\(2h_k\mu^2+2\T 4^{-k}\)\le 4 L^2_f\Tilde\nu\delta^2_k+8 L^2_f\mu^a_k+\(8+16L^2_f\)\(h_k\mu^2+\T4^{-k}\).
\end{aligned}
$$
which justify the convergence of the derivatives $\nn \dot x^k(\cdot)\hnn$ to $\dot x(\cdot)$ in $L^2\([0,\T];\R^n\)$. Using finally the Newton-Leibniz formula and the Lebesque dominated convergence theorem, we conclude that the state sequence $\{x^k(\cdot)\}$ converges strongly to $x(\cdot)$ in $W^{1,2}([0,\T];\R^n)$ as $k\to\infty$.\\[0.5ex]
\noindent
{\bf Step~4:} {\it Verifying the strong convergence of $\nn  u^k(\cdot)\hnn$ to $u(\cdot)$ in $W^{1,2}\([0,\T];\R^n\)$.} Observe first that 
$$
\n u^k_j-u(t_j)\en = \n x^k_j-x(t_j) \en \le \delta_k,
$$
where $\delta_k$ is defined in \eqref{x-est1a}. It follows from the Lipschitz continuity of $\nu_i(\cdot)$ that 
$$
\nu_i(u^k_j) \leq \nu_i(u(t_j)) + L_\nu \n u^k_j-u(t_j)\en \leq L_\nu\delta_k,
$$
which justifies the inclusions $u^k(t_j) \in U_k$ in \eqref{e:ua-dc}. Picking any $t\in [0,\T]$ ensures that $t\in [t_j, t_{j+1})$ for some $j\in\{1,\ldots,k-1\}$. Then it follows from the definition of $u^k_j$ that
\begin{equation}
\label{u-est}	
\begin{aligned}
&\n u^k(t) - u(t)\en\le \n u^k_j-u(t)\en + \n u^k_{j+1}-u^k_j\en\\
&\le \n u^k_j -u(t_j)\en + \n u(t_j)-u(t)\en + \n x^k_{j+1}-x(t_{j+1})\en+\n x(t_j)-x^k_j\en+\n u(t_{j+1})-u(t_j)\en\\
&\leq 3\delta_k+ \n u(\al^k(t))-u(t)\en +  \n u(\be^k(t))-u(\al^k(t))\en:=\mu^u_k,
\end{aligned}
\end{equation}
where $\al^k(\cdot)$ and $\be^k(\cdot)$ are defined, respectively, by 
$$
\al^k(t) = 
\begin{cases}
t_j &\mbox{if  }t\in [t_j,t_{j+1}), \\
t_{k-1} &\mbox{if  } t=\T,	
\end{cases}
$$ 
$$
\be^k(t) = 
\begin{cases}
t_{j+1} &\mbox{if  }t\in [t_j,t_{j+1}), \\
\T &\mbox{if  } t=\T.	
\end{cases}
$$ 
It is clear that $\al^k(t) \to t$ and $\beta^k(t)\to t$ uniformly on $[0, \T]$, which yields the uniform convergence (and hence the $L^2$-strong convergence) of $\nn  u^k(\cdot)\hnn$ to $u(\cdot)$. It remains to justify the convergence of $\nn \dot u^k(\cdot)\hnn $ to $\dot u(\cdot)$ in $L^2([0,T];\R^n)$. To this end, pick any $t\in[t_j,t_{j+1})$ and use \eqref{e:var} to obtain 
$$
\begin{array}{ll}
&\n \dot u^k(t)-\dot u(t)\en = \n \frac{u^k_{j+1}-u^k_j}{h_k}-\dot u(t)\en \\
&\le\Big\|\frac{x^k_{j+1}-x^k_j}{h_k} -\frac{x(t_{j+1})-x(t_j)}{h_k}\Big\| +\Big\| \frac{u(t_{j+1})-u(t_j)}{h_k}-\dot u(t)\Big\|
=\Big\|\dot x^k(t)-\disp\int^{t_{j+1}}_{t_j}\Big\|\frac{\dot x(s)}{h_k}ds\Big\| +\disp\Big\|\int^{t_{j+1}}_{t_j}\frac{\dot u(s)}{h_k}ds-\dot u(t)\Big\|.
\end{array}
$$
We get therefore the relationships
$$
\begin{array}{ll}
&\n \dot x^k(t)-{\disp\int}^{t_{j+1}}_{t_j}\frac{\dot x(s)}{h_k}ds\en \le \n \dot x^k(t) - \dot x(t_j)\en + \n \dot x(t_j)-{\disp\int}^{t_{j+1}}_{t_j}\frac{\dot x(s)}{h_k}ds\en \\[2ex]
&\le  \n \dot x^k(t) - \dot x(t_j)\en +{\disp\int}^{t_{j+1}}_{t_j}\n \frac{\dot x(t_j)-\dot x(s) }{h_k}\en ds= \n v^k_j-\dot x(t_j)\en+{\disp\int}^{t_{j+1}}_{t_j} \frac{f^x(s)}{h_k} ds\\[2ex]
&\le \n v^k_j-\dot x(t_j)\en + {\disp\int}^{t_{j+1}}_{t_j}\frac{\n\dot x(s^x_j)-\dot x(t_j) \en+\n \dot x(t_{j+1})-\dot x(s^x_j)\en+2^{-k}}{h_k}ds\\[2ex]
&= \n v^k_j-\dot x(t_j)\en +\n\dot x(s^x_j)-\dot x(t_j) \en+\n \dot x(t_{j+1})-\dot x(s^x_j)\en+2^{-k}.
\end{array}
$$
The same arguments as above bring us to estimate
$$
\n \int^{t_{j+1}}_{t_j}\dfrac{\dot u(s)}{h_k}ds-\dot u(t)\en \le \n\dot u(s^u_j)-\dot u(t_j) \en+\n \dot u(t_{j+1})-\dot u(s^u_j)\en+2^{-k},
$$
where $s^u_j\in [t_j,t_{j+1}]$ is a number such that $\disp\sup_{s\in[t_j,t_{j+1}]} \n \dot u(t_j)-\dot u(s)\en \le \n \dot u(t_j)-\dot u(s^u_j)\en + 2^{-k}$. Therefore,
$$
\begin{aligned}
&\disp\int^{\T}_0 \n \dot u^k(t)-\dot u(t)\en^2 dt= \sum^{k-1}_{j=0} \int^{t_{j+1}}_{t_j}\n \dot u^k(t)-\dot u(t)\en^2 dt\\
&\le 2 \sum^{k-1}_{j=0} \int^{t_{j+1}}_{t_j}\[\n \dot x^k(t)-\int^{t_{j+1}}_{t_j}\dfrac{\dot x(s)}{h_k}ds\en ^2 +\n \int^{t_{j+1}}_{t_j}\dfrac{\dot u(s)}{h_k}ds-\dot u(t)\en^2\]dt\\
&\le 2 \sum^{k-1}_{j=0} \int^{t_{j+1}}_{t_j}\[ \n v^k_j-\dot x(t_j)\en +\n\dot x(s^x_j)-\dot x(t_j) \en+\n \dot x(t_{j+1})-\dot x(s^x_j)\en+2^{-k}\]^2\\
&\quad+2\sum^{k-1}_{j=0} \int^{t_{j+1}}_{t_j}\[  \n\dot u(s^u_j)-\dot u(t_j) \en+\n \dot u(t_{j+1})-\dot u(s^u_j)\en+2^{-k} \]^2\\
&\le 6\sum^{k-1}_{j=0} h_k \n v^k_j-\dot x(t_j)\en^2+6\sum^{k-1}_{j=0} h_k\[ \n\dot x(s^x_j)-\dot x(t_j) \en+\n \dot x(t_{j+1})-\dot x(s^x_j)\en\]^2\\
&\quad +6 k h_k 4^{-k} +4\sum^{k-1}_{j=0} h_k\[\n\dot u(s^u_j)-\dot u(t_j) \en+\n \dot u(t_{j+1})-\dot u(s^u_j)\en \]^2+4k h_k 4^{-k}\\
&\le 24 L^2_f\T \delta^2_k+48 L^2_f\mu^a_k+48\(h_k\mu^2+\T 4^{-k}\) + 6 h_k\var^2\(\dot x(\cdot);[0,\T] \) +4 h_k\var^2\(\dot u(\cdot);[0,\T] \)+10\T 4^{-k} \\
&\le 24 L^2_f\T \delta^2_k+48 L^2_f\mu^a_k+48\(h_k\mu^2+\T 4^{-k}\) + \mu^2+10\T 4^{-k},
\end{aligned}
$$
which justifies our claim.  The endpoint constraints \eqref{ec} follows from the constructions of the state-control sequence $\{x^k(\cdot),u^k(\cdot),a^k(\cdot)\}$.\\[0.5ex] 
\noindent 
{\bf Step~5:}  {\it Verifying the estimate \eqref{uk-est}.} We have from the above that
\begin{equation}
\label{var-uk}
\begin{array}{ll}
\disp\sum^{k-2}_{j=0}&\disp\l\|\frac{u^k(t_{j+2})-u^k(t_{j+1})}{h_k} - \frac{u^k(t_{j+1})-u^k(t_j)}{h_k}\r\|
\leq\disp\sum^{k-2}_{j=0}\l\| \frac{u(t_{j+2})-u(t_{j+1})}{h_k}-\frac{u(t_{j+1})-u(t_j)}{h_k}\r\|\\
&+\disp\sum^{k-2}_{j=0}\l\| \frac{u^k(t_{j+2})-u(t_{j+2})}{h_k}-\frac{u^k(t_{j+1})-u(t_{j+1})}{h_k}\r\| 
+\disp\sum^{k-2}_{j=0}\l\| -\frac{u^k(t_{j+1})-u(t_{j+1})}{h_k}+ \frac{u^k(t_j)-u(t_j)}{h_k}\r\|\\
&\leq \mu + 2\disp\sum^{k-1}_{j=0}\l\|\frac{x^k(t_{j+1})-x(t_{j+1})}{h_k} -\frac{x^k(t_j)-x(t_j)}{h_k} \r\|=\mu + 2\disp\sum^{k-1}_{j=0}\l\|\frac{x^k(t_{j+1})-x^k(t_{j})}{h_k}  -\frac{x(t_{j+1})-x(t_j)}{h_k} \r\|\\
&=\mu +2\disp\sum^{k-1}_{j=0}\n \dot x^k(t)-\disp\int^{t_{j+1}}_{t_j}\frac{\dot x(s)}{h_k}ds\en\leq \mu + 2\disp\sum^{k-1}_{j=0}\l[  \n v^k_j-\dot x(t_j)\en +\n\dot x(s^x_j)-\dot x(t_j) \en+\n \dot x(t_{j+1})-\dot x(s^x_j)\en+2^{-k}\r] \\
&\leq \mu + 2\var (\dot x(\cdot);[0,\T]) +2k2^{-k} +  2\disp\sum^{k-1}_{j=0}\n v^k_j-\dot x(t_j)\en.\numberthis 	
\end{array}
\end{equation}
On the other hand, it follows from \eqref{mu-k}, \eqref{dv1}, and \eqref{x-est1a} that 
$$
\begin{aligned} 
&\sum^{k-1}_{j=0}\n v^k_j-\dot x(t_j)\en\leq L_f\sum^{k-1}_{j=0}\n x^k_j-x(t_j)\en +L_f\sum^{k-1}_{j=0}\n a^k_j-a(t_j)\en\\
&\leq L_fk\delta_k + L_f\var (a(\cdot);[0,\T ])\leq L_fe^{L_f\T}(2L_f+1)(\T\mu +\T k2^{-k}) + L_f\mu\leq 2L_f\T\mu e^{L_f\T}(2L_f+1) + L_f\mu.
\end{aligned}
$$
Combining the latter with \eqref{var-uk} gives us 
$$
\mbox{var}\l(\dot u^k;[0,\T] \r)  \leq  3\mu+1+4L_f\T\mu e^{L_f\T}(2L_f+1) + 2L_f\mu \leq \Tilde\mu,
$$
where $\Tilde\mu$ is defined in \eqref{tmu}.
To verify the second estimate therein, we deduce from~\eqref{e:var} and~\eqref{x-est1a} that 
$$
\begin{array}{ll}
&\n\frac{u^k(t_1)-u^k(0)}{h_k}\en \leq \n \frac{u^k(t_1)-u(t_1)}{h_k}\en  +\n \frac{u^k(0)-u(0)}{h_k}\en + \n\frac{u(t_1)-u(0)}{h_k} \en\\
&\leq \n \frac{x^k(t_1)-x(t_1)}{h_k}\en  +\n \frac{x^k(0)-x(0)}{h_k}\en + \mu\leq 2e^{L_f\T }(2L_f+1)(\mu+k2^{-k}) + \mu\leq 2e^{L_f\T }(2L_f+1)(\mu+1) + \mu\leq \Tilde\mu.
\end{array}
$$
In the same way, we also get $\n\frac{u^k(T)-u^k(t^k_{k-1})}{h_k}\en\leq \Tilde\mu$ and thus complete the proof of the theorem. $\h$\vspace*{-0.2in}

\section{Existence of Optimal Solutions and Local Minimizers}\label{sec:exist}
\setcounter{equation}{0}\vspace*{-0.1in}

This section begins the study of the optimal control problem $(P)$ in \eqref{e:BP}--\eqref{A}. We first address the {\em existence of optimal solutions} to problem $(P)$. To proceed, let us consider the other form of the sweeping inclusion, which---in contrast to \eqref{e:SP}---doesn't explicitly contains measurable controls $a(\cdot)\in A$:
\eq
\label{e:SP2}
\dot x(t) \in -N_{C_{u(t)}}\big(x(t)\big) + f\big(x(t), A\big) \quad\mbox{ a.e. }\;t \in [0,T], 
\eeq
where the image set $f(x, A)$ is defined by 
$$
f(x,A): =\big\{v\in \R^n\;\big|\;v= f(x,a)\;\mbox{ for some }\;a \in A\big\}, \; x\in \R^n. 
$$
The following proposition shows that the above forms are equivalent.\vspace*{-0.1in}

\begin{proposition}\label{prop:equi} The sweeping inclusions in \eqref{e:SP} and \eqref{e:SP2} are equivalent to each other.
\end{proposition}\vspace*{-0.12in}
{\bf Proof}. It is obvious that the sweeping inclusion in form \eqref{e:SP} yields that in \eqref{e:SP2}. To verify the converse implication, define the set-valued mapping $S:[0, T]\tto \R^n$ by 
\eq
\label{e:SVM}
S(t): = \l\{ a\in A\;\big|\;-\dot x(t) \in  F\big(x(t),u(t),a)\big) = N_{C_{u(t)}}\big(x(t)\big) + f\big(x(t); a\big)\r\},
\eeq
which is closed-valued for a.e. $t\in[0,T]$. We claim that the mapping $S$ is measurable on $[0,T]$. Indeed, employing the classical Luzin theorem from real analysis to the measurable function $-\dot\ox(\cdot)$ gives us a closed set $T_{\ve}\subset[0,T]$ with $\mbox{mes}([0,T]\backslash T_\ve)<\ve$, where $\mbox{mes}(\cdot)$ denotes the Lebesgue measure, such that $-\dot\ox(\cdot)$ is continuous on $T_\ve$ for any given $\ve>0$. It follows from the continuity of $f(x,a)$ together with the closed-graph property of the normal cone \eqref{e:SVM} that the restricted mapping $S:T_\ve\tto \R^n$ in \eqref{e:SVM} is of closed graph. Applying \cite[Theorem~14.10]{rw} to the set-valued mapping $S$ ensures its measurability on $[0,T]$. Hence we can find a measurable control $a(\cdot)\in \cA$ such that the triple $(x(\cdot),u(\cdot),a(\cdot))$ is feasible to the differential inclusion \eqref{e:SP1} by the measurable selection theorem from \cite[Corollary~14.6]{rw}. This justifies the equivalence between the inclusions in \eqref{e:SP} and \eqref{e:SP2}. $\h$ 

The next theorem provides conditions ensuring the {\em existence of optimal solutions} to  the sweeping optimal control problem $(P)$ defined in Section~\ref{S1}.\vspace*{-0.15in}

\begin{theorem}\label{Th2} Fix $\ve>0$ and suppose that along some minimizing sequence $\{\ox^k(\cdot),\ou^k(\cdot), \oa^k(\cdot), \T^k\}$ in problem $(P)$, the assumptions in {\bf(H1)}--{\bf(H5)} are satisfied with $T=\T^k+\ve$ as $k\in\N$. Suppose in addition that  $\{\ou_k(\cdot)\}$ is bounded in $W^{1,2}([0,\T^k+\ve];\R^n)$, that $\{\oa^k(0), \T^k\}$ is bounded, that  
$\var(\oa^k(\cdot);[0,\T^k+\ve])\leq \mu$,  and that the running cost $\ell$ in \eqref{e:BP} is convex with respect to the velocity variables $(\dot x, \dot u)$ while $\ell_t(t,\cdot)$ is majorized by a summable function. Then problem $(P)$ admits an optimal solution $(\ox(\cdot),\ou(\cdot),\oa(\cdot),\T)$ in  $W^{1,2}([0,\T];\R^n)$ $\times$ $ W^{1,2}([0,\T];\R^n)$ $\times L^2([0,\T];\R^d)\times [0,\infty)$. 
\end{theorem}\vspace*{-0.12in}
{\bf Proof}.
Observe by Proposition~\ref{Prop1} that the set of feasible solutions to problem $(P)$ is nonempty, which enables us to select the minimizing sequence  
$$
\l(\ox^k(\cdot),\ou^k(\cdot), \oa^k(\cdot), \T^k\r)\in W^{1,2}([0,\T^k];\R^{2n})\times L^2([0,\T^k];\R^d) \times [0,\infty)
$$ 
in $(P)$. The boundedness of $\{\ou_k(\cdot)\}$ in $W^{1,2}([0,\T^k+\ve];\R^n)$ yields the boundedness of  $\{\dot\ou^k(\cdot)\}$ in $L^2([0,\T^k+\ve];\R^n)$. Hence the Banach-Alaoglu theorem ensures the existence of a function $\vartheta^u(\cdot)$ such that $\dot\ou^k(\cdot)\to \vartheta^u(\cdot)$ as $k\to\infty$ weakly in {$L^2([0,\T^k+\ve];\R^n)$. It follows from the boundedness of $\{\oa^k(0)\}$ and the uniform bounded variation property of $\{\oa^k(\cdot)\}$ that $\{\oa^k(\cdot)\}$ is bounded on $[0, \T^k+\ve]$. Employing now Helly's selection theorem gives us a function of bounded variation $\oa(\cdot)$ such that $\oa^k(t) \to \oa(t)$ as $k\to\infty$ for all $t\in [0,\T]$, where $\T=\lim_{k\to\infty}\T^k+\varepsilon$ along some subsequence. We get by the compactness of  the set $A$ that $\oa(t)\in A$ for all $t\in [0, \T]$. Since $\{\ou^k(0)\}$ belongs to the compact set $U$, there exists a subsequence of $\{\ou^k(0)\}$, still denoted by $\{\ou^k(0)\}$, such that $\{\ou^k(0)\}$ converges to some $u_0\in\R^n$. Define further the function $\ou(\cdot)\in W^{1,2}([0,\T];\R^n)$ as follows
\eq 
\label{e:ubar}
\ou(t): = u_0 + \int^t_0\vartheta^u(s)ds\;\mbox{ for all }\;t\in [0, \T]
\eeq 
and conclude that $\dot\ou^k(\cdot) \to \dot\ou(\cdot)$ as $k\to\infty$ weakly in $L^2([0, \T^k+\varepsilon];\R^n)$. Hence we get the pointwise convergence $\ou^k(t) \to \ou(t)$,  which implies that $\ou(t) \in U$ for all $t\in [0, \T]$. Next we deduce from estimate \eqref{e:boundedness} in Proposition~\ref{Prop1} for $T=\T^k+\ve$ and from the boundedness of $\{\dot\ou^k(\cdot)\}$ in {$L^2([0,\T^k+\ve];\R^n)$} that $\{\dot\ox^k(\cdot)\}$ is bounded in {$L^2([0,\T^k+\ve];\R^n)$}. This justifies the weak convergence of $\{\dot\ox^k(\cdot)\}$ to $\dot\ox(\cdot)$ {in $L^2([0,\T];\R^n)$} and ensures that $\ox^k(t)\to \ox(t)$ as $k\to\infty$, where $\ox(\cdot) \in W^{1,2}([0,\T];\R^n)$ is defined by
\eq 
\label{e:xbar}
\ox(t): = x_0 + \int^t_0\vartheta^x(s)ds\;\mbox{ for all }\;t\in [0, \T]
\eeq 
with $\vartheta^x(t) = \lim_{k\to\infty}\dot\ox^k(t)$ along some subsequence. To verify the feasibility of $\l(\ox(\cdot), \ou(\cdot), \oa(\cdot), \T\r)$ to the sweeping differential inclusion \eqref{e:SP1}, we elaborate the arguments similar to those developed in the proof of \cite[Theorem~4.1]{cm3} for $T=\T$. Using finally the Lebesgue dominated convergence theorem together with the imposed convexity of the integrand $\ell$ with respect to $(\dot x, \dot u)$ brings us to
$$
J[\ox, \ou, \oa, \T] \leq \liminf_{k\to\infty}J\l[\ox^k, \ou^k, \oa^k, \T^k+\varepsilon\r]
$$
due to the lower semicontinuity of integral functionals with respect to the weak topology in $L^2$. This  justifies the optimality of $(\ox(\cdot), \ou(\cdot), \oa(\cdot), \T)$ in problem $(P)$ and thus completes the proof. $\h$
 
As seen from Theorem~\ref{Th2}, the existence of optimal solutions to problem $(P)$ requires the {\em convexity} of the running cost with respect to velocity variables. On the other hand, it has been well recognized in variational and control theories, starting with the classical Bogolyubov-Young theorem in the calculus of variations, that it is possible to {\em relax} a given nonconvex problem to a certain {\em convexified} one, which admits an optimal solution that can be approximated by feasible solutions of the original problems and keeps the same value of the cost functional. For various controlled sweeping processes, such a relaxation procedure was implemented in, e.g., \cite{CMN20,ET,tolstonogov}.\vspace*{-0.05in}

In the setting of problem $(P)$, we proceed as follows. Let $\hat \ell_F(t, x, u, a, \dot x, \dot u)$ be the convexification 
of the running cost in \eqref{e:BP} (i.e., the largest l.s.c. convex function majorized by $\ell(t,x,u,a,\dot x, \dot u)$)  on the set $F(x,u,a)$ in \eqref{e:F-rep}. Along with problem $(P)$, consider the  {\it relaxed optimal control problem $(R)$} defined by 
\eq
\label{Relax}
\mbox{minimize } \widehat J[x,u,a,T]: = \ph(x(T),T) + \int^T_0\hat\ell_F(t,x(t),u(t),a(t),\dot x(t),\dot u(t))dt
\eeq
over $(x(\cdot), u(\cdot), a(\cdot),T)\in  W^{1,2}([0,T];\R^n\times\R^n)\times L^2([0,T];\R^d)\times [0,\infty)$ satisfying  
\eq 
\label{e:RSP}
\begin{cases}
-\dot x(t) \in N_{C_{u(t)}}(x(t)) + \co f(x(t), A) \mbox{ a.e. } t \in [0,T],\\
x(0) =   x_0 \in C_u(0),\\
(u(t),a(t)) \in U\times A \mbox{ for all } t\in [0,T].
\end{cases}
\eeq
Note that we don't need to convexify the normal cone $N_{C_{u(t)}}(x(t))$ in \eqref{e:RSP}, since the assumptions imposed on the moving set guarantee the cone convexity.\vspace*{-0.05in}

Now we describe two notions of {\em local minimizers} in $(P)$, which are of our interest in this paper.\vspace{-0.1in}

\begin{definition}\label{D5.1} We say that the quadruple $(\ox,\ou,\oa,\T)$ is:

{\bf(i)} An {\sc intermediate local minimizer} $($i.l.m.$)$ in problem $(P)$ if there exists $\ve>0$ such that
$$
J[\ox,\ou,\oa,\T]\le J[x,u,a,T]
$$
for all feasible solutions $(x,u,a,T)$ to $(P)$ that satisfy the localization conditions
\eq\label{nb} 
\begin{cases}
\disp\int^{\T}_0\n (\ox(t),\ou(t),\oa(t)) - (x(t),u(t),a(t))\en^2dt < \ve, \\[1ex]
\disp\int^{\T}_0\n( \dot\ox(t), \dot\ou(t)) - (\dot x(t), \dot u(t)) \en^2 dt  <\ve, \\[1ex]
|T-\T|<\ve,  
\end{cases}
\eeq 
\vspace*{-0.2in}
 
{\bf(ii)} A  {\sc relaxed intermediate local minimizer} $($r.i.l.m.$)$ in $(P)$ if it is feasible to $(P)$, provides an i.l.m.\ for the relaxed problem $(R)$ with $\widehat J[\ox,\ou,\oa,\T]=J[\ox,\ou,\oa,\T]$, and there exists $\ve>0$ such that the localization conditions in \eqref{nb} are satisfied.  
\end{definition}\vspace*{-0.15in} 

It follows from Definition~\ref{D5.1}(ii) and the form of the cost functional \eqref{Relax} that the r.i.l.m. $(\ox,\ou,\oa,\T)$ therein satisfies the relaxed system \eqref{e:RSP} with the convexified sweeping inclusion. The notions of intermediate local minimizers and their relaxed counterpart, introduced in \cite{m95} for Lipschitzian differential inclusions (see also \cite{m96} for such systems with free time), occupy an intermediate position between weak and strong minima in the classical calculus of variations and optimal control. Versions of these notions for various controlled sweeping processes can be found in \cite{CCMN21,chhm,CMN20,cmnn24}.\vspace*{-0.05in} 

While there is obviously no difference between i.l.m.\ and r.i.l.m.\ for problems convex in velocity variables, such a {\em relaxation stability} holds for large classes of nonconvex control problems of Lipschitzian and sweeping types. This is due to the {\em hidden convexity} of continuous-time differential systems; see more discussions in \cite{CMN20,ET,m-book1,tolstonogov} and the references therein.\vspace*{-0.25in}

\section{Discrete Approximations of Local Minimizers}\label{sec:disc-local}
\setcounter{equation}{0}\vspace*{-0.15in}

In this section, we construct and justify well-posedness of discrete approximations of a designated r.i.l.m.\ for the sweeping control problem $(P)$ under consideration. Let $(\ox(\cdot),\ou(\cdot),\oa(\cdot), \T)$ be a r.i.l.m.\ in $(P)$. For each $k\in \N$, construct a sequence of discrete sweeping control problems $(P_k)$  on the discrete mesh {in~\eqref{e:DP} as} defined  in the following way: minimize the cost function
\eq
\label{DP}
\begin{array}{ll}
&J_k[x^k,u^k,a^k, T^k]:= \ph\(x^k_{k},T^k\) + \sum\limits^{k-1}_{j=0} h_k\ell\(t^k_j, x^k_j,u^k_j,a^k_j, \frac{x^k_{j+1}-x^k_j}{h_k}, \frac{u^k_{j+1}-u^k_j}{h_k}\)\\[2ex]
&+\frac{1}{2}(T^k-\T)^2
+\frac{1}{2}\sum\limits^{k-1}_{j=0}{\disp\int}^{t^k_{j+1}}_{t^k_j}\n \l(a^k_j, \frac{x^k_{j+1}-x^k_j}{h_k}, \frac{u^k_{j+1}-u^k_j}{h_k}\r) - \l(\oa(t),\dot \ox(t), \dot \ou(t)  \r) \en^2 dt \\[2ex]
&+\dist^2\l(\l\| \frac{ u^k_1-u^k_0}{h_k}\r\|;(-\infty,\Tilde\mu ]\r) + \dist^2\l(\sum\limits^{k-1}_{j=1}\l\|\frac{u^k_{j+1}-2u^k_j+u^k_{j-1}}{h_k} \r\|;(-\infty,\Tilde\mu)\r)
\end{array}
\eeq 
over elements $(x^k_0,x^k_1,\ldots,x^k_{k},u^k_0,\ldots, u^k_k,a^k_0,\ldots,a^k_{k-1},T^k)$ satisfying the constraints
\begin{equation}
\label{con1}
x^k_{j+1}\in x^k_j-h_kF\(x^k_j,u^k_j,a^k_j\)\;\;\mbox{for}\;\;j=0,\ldots, k-1,
\end{equation}
\begin{equation}
\label{con2}
\(x^k_0, u^k_0 \) = \(x_0, \ou(0)\),
\end{equation}
\begin{equation}
\label{con2a}
\left (x^k_{k}, T^k \right) \in \Xi^k_x\times \Xi^k_T\colon =  (\Xi_x+\mu^x_kB(0,1))\times (\Xi_T+\mu^x_k),
\end{equation}
\eq 
\label{con3}
T^k \leq \T +\ve,
\eeq
\eq
\label{con3a}
u^k_j \in U_k \mbox{ for } j=0, \ldots, k-1,
\eeq
\begin{equation}
\label{con4}
a^k_j \in  A\;\;\mbox{for}\;\;j=0,\ldots,k-1,
\end{equation}
\begin{equation}
\label{con5}
\begin{array}{ll}
{\disp\sum^{k-1}_{j=0}\int^{t^k_{j+1}}_{t^k_j}}\n \(x^k_j,u^k_j,a^k_j\)-\(\ox(t),\ou(t),\oa(t)\)\en^2 dt  \le \frac{\ve}{2},
\end{array}
\end{equation}
\begin{equation}
\label{con6}
\begin{array}{ll}
{\disp\sum\limits^{k-1}_{j=0}\int^{t^k_{j+1}}_{t^k_j}}\n \l( \frac{x^k_{j+1}-x^k_j}{h_k}, \frac{u^k_{j+1}-u^k_j}{h_k}\r) - \l(\dot \ox(t), \dot \ou(t) \r) \en^2 dt \leq \frac{\ve}{2},
\end{array}
\end{equation}
\begin{equation}
\label{con7}
\begin{array}{ll}
\max\l\{\l\| \frac{ u^k_1-u^k_0}{h_k}\r\|,\l\|\frac{u^k_k-u^k_{k-1}}{h_k}\r\|  \r\}\leq \Tilde\mu+1 \mbox{ and }\sum\limits^{k-1}_{j=1}\l\|\frac{u^k_{j+1}-2u^k_j+u^k_{j-1}}{h_k} \r\| \leq \Tilde\mu+1 
\end{array}
\end{equation}
where $\delta_k, \mu^x_k$, and $\Tilde\mu$ are given in Theorem~\ref{Th1}, and where $\ve$ is taken from Definition~\ref{D5.1}.\vspace*{-0.05in} 

To proceed further, we imposed the following assumption on the endpoint constraint set $\Xi_x\times\Xi_T$:

{\bf(H6)} The set $\Xi_x\times\Xi_T$ is closed around $(\ox(\T),\T)$. 

\noindent To study the relationships between problems $(P_k)$ and $(P)$, we first verify the {\em existence of optimal solutions} to $(P_k)$ for large $k$, which is a must issue for the method of discrete approximations.\vspace*{-0.1in}

\begin{proposition}
\label{P5.2}
In addition to the assumptions of Theorem~{\rm\ref{Th1}}, suppose that {\rm(H5)} and {\rm(H6)} are satisfied along the r.i.l.m. $(\ox,\ou,\oa,\T)$. Then each problem $(P_k)$ defined in \eqref{DP}--\eqref{con7} admits an optimal solution $\l(\ox^k,\ou^k,\oa^k, \T^k\r)$ whenever $k$ is sufficiently large. 
\end{proposition}\vspace*{-0.15in}
\begin{proof} The existence of optimal solutions to finite-dimensional problems $(P_k)$ follows directly from the classical Weierstrass existence theorem provided that the set of feasible solutions to $(P_k)$ is nonempty, bounded. and closed. The nonemptiness of the feasible sets for all large $k$ follows from Theorem~\ref{Th1} applied to the designated local minimizers of $(P)$. Furthermore, the imposed constraints \eqref{con4}--\eqref{con6} these sets are bounded for all $k$. To justify the closedness, fix $k$ and take a sequence 
$$
z^m: = (x_0^m, \ldots, x^m_{k}, u^m_0, \ldots,  u^m_k, a^m_0, \ldots, a^m_{k-1}, T^m_k) 
$$
of feasible solutions for $(P_k)$ converging to $z: = (x_0, \ldots, x_{k}, u_0, \ldots, u_k, a_0, \ldots, a_{k-1}, T^k)$ and check that $z$ is feasible to $(P_k)$. Indeed, we have $g_i(x_j-u_j)=\disp\lim_{m\to\infty}g_i(x^m_j-u^m_j)\ge 0$ for all $i=1,\ldots, m$ and $j=0,\ldots, k-1$ which implies that $x_j - u_j\in C$ for all $j=0,\ldots, k-1$. It is not hard to see that $I(x^m_j-u^m_j) \subset I(x_j-u_j)$ for $m\in \N$ sufficiently large. This allows us to deduce from \eqref{con1} and \eqref{F} the inclusion
$$
\dfrac{x^m_{j+1}-x^m_j}{-h_k} - f(x^m_j, a^m_j) \in N_C(x^m_j-u^m_j).
$$
Taking into account the limiting conditions
$$
\dfrac{x^m_{j+1}-x^m_j}{-h_k} - f(x^m_j, a^m_j)  \to  \dfrac{x_{j+1}-x_j}{-h_k} - f(x_j, a_j) \;\; \mbox{and}\;\; x^m_j-u^m_j \to x_j-u_j
$$
as $m\to\infty$ and the robustness of the limiting normal cone, the latter inclusion readily implies that
$$
\dfrac{x_{j+1}-x_j}{-h_k} - f(x_j, a_j) \in \Limsup_{x-u\to x_j-u_j}N_C(x-u)= N_C(x_j-u_j)
$$
ensuring in this way that $x_{j+1}-x_j\in -h_kF(x_j,u_j,a_j)$ for all $j=0,\ldots, k-1$. It is obvious that $z$ also satisfies the constraints \eqref{con2}--\eqref{con6}. This verifies that the set of feasible solutions to $(P_k)$ is closed and thus completes the proof of the proposition.  
\end{proof}

The next theorem justifies the {\em strong convergence} (in the corresponding spaces) of the optimal solutions $\l(\ox^k,\ou^k,\oa^k, \T^k\r)$ to $(P_k)$, properly extended to the continuous-time interval $[0,\T]$, to the designated optimal solution $(\ox,\ou,\oa,\T)$ of the sweeping control problem $(P)$. This ensures that $\l(\ox^k,\ou^k,\oa^k, \T^k\r)$ can be treated as {\em approximately optimal/suboptimal} solutions to $(P)$, and hence the necessary optimality conditions for $\l(\ox^k,\ou^k,\oa^k, \T^k\r)$ can be treated as {\em almost optimality conditions} for  $(\ox,\ou,\oa,\T)$.\vspace*{-0.12in}

\begin{theorem}\label{Th2a}
Let $(\ox(\cdot),\ou(\cdot),\oa(\cdot),\T)$ be an r.i.l.m.\ for problem $(P)$, and let all the assumptions of Proposition~\ref{P5.2} be satisfied for this quadruple. Suppose in addition that the terminal cost $\ph$ is continuous around $\ox(\T)$, that the running cost $\ell$ is continuous at $(t, \ox(t), \ou(t),\oa(t),\dot\ox(t),\dot\ou(t))$ and that $\ell(\cdot, x, u, a, \dot x, \dot u)$ is uniformly majorized around $(\ox(\cdot),\ou(\cdot),\oa(\cdot))$ by a summable function on $[0,\T+\ve]$. Then any sequence of optimal solutions $\l(\ox^k(\cdot),\ou^k(\cdot),\oa^k(\cdot), \T^k\r)$ of $(P_k)$, where $\l(\ox^k(\cdot),\ou^k(\cdot)\r )$ and $\oa^k(\cdot))$ are piecewise linearly and piecewise constantly extended to the whole interval $[0, \T^k]$, converges to $(\ox(\cdot),\ou(\cdot),\oa(\cdot),\T)$ in the following senses:
\eq 
\label{e:conv1}
\T^k \to \T\;\mbox{ as }\;k \to\infty,
\eeq 
\eq 
\label{e:conv2}
\max_{t\in [0,\T^k]}\l\{\n(\ox^k(t),\ou^k(t)) - (\ox(t),\ou(t))\en \r\}\to 0\;\mbox{ as }\;k \to\infty,
\eeq 
\eq 
\label{e:conv3}
\int^{\T^k}_0\n \l( \dot\ox^k(t), \dot\ou^k(t),\oa^k(t)\r) - \l(\dot \ox(t), \dot\ou(t),\oa(t)\r) \en^2 dt \to 0\;\mbox{ as }\;k \to \infty
\eeq 
with the fulfillment of the estimates
\eq
\label{e:conv4}
\begin{array}{ll}
\max\l\{\l\| \frac{ \ou^k_1-\ou^k_0}{h_k}\r\|, \l\|\frac{\ou^k_k-\ou^k_{k-1}}{h_k}\r\| \r\}\leq \Tilde\mu\;\mbox{  and  }\; {\disp\limsup_{k\to\infty }\sum^{k-1}_{j=1}}\l\|\frac{\ou^k_{j+1}-2\ou^k_j+\ou^k_{j-1}}{h_k} \r\|  \leq \Tilde\mu,
\end{array}
\eeq
where the number $\Tilde\mu>0$ is taken from 
Theorem~{\rm\ref{Th1}}. 
\end{theorem}\vspace*{-0.1in}
{\bf Proof}.
Take any sequence $\l(\ox^k(\cdot),\ou^k(\cdot),\oa^k(\cdot), \T^k\r)$ of extended optimal solutions of $(P_k)$, where $\{(\ox_k(t^k_j),$ $\ou^k(t^k_j))\;\big|\; j=0,\ldots, k\}$ and $\{\oa^k(t^k_j)\big| \; j=0,\ldots, k\}$ are piecewise linear and piecewise constant on $[0,\T^k]$, respectively, and where the mesh points $t^k_j$ are taken from $\Delta_k(\T^k)$ defined in \eqref{e:DP}. Let us show that
\eq 
\label{e:sigma}
\begin{array}{ll}
\disp\lim_{k\to\infty}&\bigg(\sigma_k:= \disp\int^{\T^k}_0\n \l( \dot\ox^k(t), \dot\ou^k(t),\oa^k(t)\r) - \l(\dot \ox(t), \dot\ou(t),\oa(t)\r) \en^2 dt+\l|\T^k-\T\r| \\[1ex]
&+\dist^2\l(\max\l\{\l\| \frac{ \ou^k_1-\ou^k_0}{h_k}\r\|, \l\|\frac{\ou^k_k-\ou^k_{k-1}}{h_k}\r\| \r\};(-\infty,\Tilde\mu ]\r) + \mbox{dist}^2\bigg(\sum\limits^{k-1}_{j=1}\l\|\frac{\ou^k_{j+1}-2\ou^k_j+\ou^k_{j-1}}{h_k} \r\| ; (-\infty, \Tilde\mu] \bigg)\bigg)=0,
\end{array}
\eeq 
which will provide the convergence of $\l(\ox^k(\cdot),\ou^k(\cdot),\oa^k(\cdot), \T^k\r)$ to $(\ox(\cdot),\ou(\cdot),\oa(\cdot),\T)$ claimed in \eqref{e:conv1}--\eqref{e:conv4}. To justify this claim, suppose on the contrary that the limit in \eqref{e:sigma}, along some subsequence (without relabeling), equals to $\sigma>0$. Using the boundedness of $\T^k$ due to \eqref{con3}, we find a real number $\tT \leq \T+\ve$ for which $\lim_{k\to\infty}\T_k = \tT$ along a subsequence (without relabeling). Moreover, it follows from the weak sequential compactness of the unit ball in $L^2$ that there exist functions $(v^x(\cdot), v^u(\cdot), \ta(\cdot)) \in L^2([0,\tT];\R^{2n+d})$ for which the triples $\l(\dot\ox^k(\cdot), \dot\ou^k(\cdot), \oa^k(\cdot)\r)$ converge weakly to $(v^x(\cdot), v^u(\cdot), \ta(\cdot))$ in $L^2([0,\tT];\R^{2n+d})$. Define further the state-control  pair $(\tx(\cdot), \tu(\cdot))\in W^{1,2}([0,\tT];\R^{2n})$ by 
$$
\l(\tx(t),\tu(t)\r) \colon=(x_0, \ou(0)) + \int^t_0(\vartheta^x(s), \vartheta^u(s))ds \;\;\mbox{ for all }\; t\in [0,\tT].
$$
Then $\l(\dot\ox^k(\cdot), \dot\ou^k(\cdot)\r)$ converges weakly to $(v^x(\cdot), v^u(\cdot)) = \l(\dot\tx(\cdot),\dot\tu(\cdot)\r)$ in $L^2([0,\tT];\R^{2n})$ and $(\ox^k(t),$ $\ou^k(t))$ converge to $(\tx(t),\tu(t))$ for a.e. $t\in [0,\tT]$. As a consequence, we get $\tx(t)-\tu(t)=\lim_{k\to\infty}(\ox^k(t)-\ou^k(t))\in C$ from \eqref{e:MC}{ and $(\tu(t),\ta(t)) \in U\times\co A$ for $t\in [0,\tT]$. Elaborating the arguments similar to those developed in the proof of \cite[Theorem~4.1]{cm3} for $T=\tT$ gives us 
$$
-\dot\tx(t) - f(\tx(t),\ta(t)) \in N_{C_{\tu}(t)}(\tx(t)) \;\;\mbox{ for a.e. }\;t\in [0,\tT],
$$
which justifies the fulfillment of the convexified sweeping differential inclusion in \eqref{e:RSP} for the limiting process. The convexity of the norm function and hence its lower semicontinuity in the $L^2$-weak topology allows us to arrive at the following estimates by passing \eqref{con5} and \eqref{con6} to the limit as $k\to\infty$:
$$
\begin{cases}
\l|\tT - \T \r| \leq \ve \\
\disp\int^{\Tilde T}_0\n (\ox(t),\ou(t),\oa(t)) - (\tx(t),\tu(t),\ta(t))\en^2dt  \\
\leq \disp\liminf_{k\to\infty}\sum^{k-1}_{j=0}\int^{t^k_{j+1}}_{t^k_j}\n 
\(\ox^k(t^k_j),\ou^k(t^k_j),\oa^k(t^k_j)\)-
\(\ox(t),\ou(t),\oa(t)\)\en^2 dt \leq \frac{\ve}{2},\\
\disp\int^{\tT}_0\n \l( \dot\tx(t), \dot\tu(t)\r) - \l(\dot \ox(t), \dot\ou(t)\r) \en^2 dt\\
\le\disp\liminf_{k\to\infty}\sum^{k-1}_{j=0}\int^{t^k_{j+1}}_{t^k_j}\n \l( \frac{\ox^k(t^k_{j+1})-\ox^k(t^k_j)}{h_k}, \frac{\ou^k(t^k_{j+1})-\ou^k(t^k_j)}{h_k}\r) - \l(\dot \ox(t), \dot \ou(t) \r) \en^2 dt \leq \frac{\ve}{2}.
\end{cases}
$$
As a result, the quadruple $\l(\tx(\cdot),\tu(\cdot),\ta(\cdot),\tT\r)$ belongs to the given $\ve$-neighborhood of the r.i.l.m. $\big(\ox(\cdot),\ou(\cdot),$ $\oa(\cdot),\T\big)$ in the space of $W^{1,2}([0,\T];\R^{2n})$ $\times L^2([0,\T];\R^d)\times [0,\infty)$. Applying now Theorem~\ref{Th1} to the r.i.l.m. $\l(\ox(\cdot),\ou(\cdot),\oa(\cdot),\T\r)$ enables us to select a sequence $\l(x^k(\cdot),u^k(\cdot),a^k(\cdot),T^k=\T\r)$ of the extended feasible solutions to $(P_k)$ for which $x^k(\cdot), u^k(\cdot)$ and $a^k(\cdot)$ strongly approximate $\ox(\cdot), \ou(\cdot)$ and $\oa(\cdot)$ in $W^{1,2}([0,\T];\R^{2n})$ and $L^2([0,\T];\R^d)$, respectively. 
The imposed convexity of $\hat\ell_F$ and the optimality of $\l(\ox^k(\cdot),\ou^k(\cdot),\oa^k(\cdot), \T^k\r)$ for problem $(P_k)$ yields the relationships
$$
\begin{array}{ll}
&\hat J[\tx, \tu, \tu, \tT] + \frac{\sigma}{2} = \ph(\tx(\tT),\tT) + \disp\int^{\tT}_0\widehat\ell_F(t,\tx(t),\tu(t),\ta(t),\dot\tx(t),\dot\tu(t))dt + \frac{\sigma}{2}\\
&\le\disp\liminf_{k\to\infty}\l[\ph\l(\ox^k(\T^k),\T^k \r) + \disp\sum^{k-1}_{j=0}{h_k}\ell\(t^k_j, \ox^k_j,\ou^k_j,\oa^k_j, \frac{\ox^k_{j+1}-\ox^k_j}{h_k}, \frac{\ou^k_{j+1}-\ou^k_j}{h_k}\)\ + \frac{\sigma}{2}\r]\\
&= \disp\liminf_{k\to\infty}J_k[\ox^k,\ou^k,\oa^k,\T^k] \leq \liminf_{k\to\infty}J_k[x^k,u^k,a^k, T^k].
\end{array}
$$
Using next the strong convergence of $\l(x^k(\cdot),u^k(\cdot),a^k(\cdot),T^k=\T\r)$ to $\l(\ox(\cdot),\ou(\cdot),\oa(\cdot),\T\r)$ in the space $W^{1,2}([0,\T];\R^{2n})\times L^2([0,\T];\R^d)\times [0,\infty)$ from Theorem~\ref{Th1} and the imposed continuity of $\ph$ and $\ell$ gives us $J[x^k,u^k,a^k,$ $ T^k=\T]\to J[\ox, \ou, \oa, \T] $ as $k\to\infty$. Therefore,  we arrive at
$$
\begin{aligned} 
\hat J[\tx, \tu, \tu, \tT] &< \hat J[\tx, \tu, \tu, \tT] + \dfrac{\sigma}{2} \leq \liminf_{k\to\infty}J_k[x^k,u^k,a^k, T^k]
\\
&=  J[\ox, \ou, \oa, \T]  =  \hat J[\ox, \ou, \oa, \T]  = \inf(R),
\end{aligned}
$$
which is not possible since $(\ox(\cdot) \ou(\cdot), \oa(\cdot), \T)$ is an r.i.l.m.\ for $(P)$. This justifies the claimed limiting condition \eqref{e:sigma} and thus completes the proof of the theorem. $\h$\vspace*{-0.2in}

\section{Second-Order Subdifferential Evaluations}\label{tools}
\setcounter{equation}{0}\vspace*{-0.1in}

In our derivation of necessary optimality conditions for the discrete-time sweeping control problems $(P_k)$ and eventually for their continuous-time counterpart $(P)$, {\em coderivative evaluations} of velocity mapping $F$ from \eqref{F} via the problem data play a significant role. It follows from the structure of \eqref{F} that the major part of $F$ is the normal cone mapping for which the coderivative reduces to the {\em second-order subdifferential} \eqref{2nd} of the set indicator function $\delta_\O$ equal to zero on the set in question and $\infty$ otherwise.\vspace*{-0.05in}

The results presented below are mostly based on the coderivative evaluations provided in \cite{hos} with further developments given in \cite{cm3,m24}. Recall that a set-valued mapping $M\colon\R^s\tto\R^q$ is {\em calm} at $(\bar\vt,\bar q)\in\gph M$ if there exist positive numbers $\mu$ and $\eta$ such that
\begin{equation}\label{e:calm}
M(\vt)\cap(\bar q+\eta\B)\subset M(\bar\vt)+\mu\|\vt-\bar\vt\|B(0,1)\;\mbox{ whenever }\;\vt\in\bar\vt+\ell B(0,1).
\end{equation}
This is a weak ``one-point" stability property, which is readily guaranteed by robust Lipschitzian behavior of  multifunctions known as Lipschitz-like or Aubin property; see, e.g., \cite{m-book,rw}.\vspace*{-0.05in}

The first statement below, inspired by \cite[Theorem~3.3]{hos} (see \cite[Theorem~6.1]{cm3} for the detailed proof), provides a constructive upper estimate and an exact calculation of $\partial^2\dd_\O=D^*N_\O$ for the case where 
\begin{equation}\label{omega}
\O:= C=\bigcap^m_{i=1}\big\{x\in\R^n\;\big| \;g_i(x)\ge0\big\}.
\end{equation}
Here we use the standard notion
\begin{eqnarray}\label{orthant}
\R^m_-:=\big\{(y_1,\ldots,y_m)\in\R^m\;\big|\;y_i\le 0\;\mbox{ for all }\;i=1,\ldots,m\big\}
\end{eqnarray}
and recall that the {\em positive linear independence condition} (PLICQ) for $\nabla g_1,\ldots,\nabla g_m$ holds at $\ox$ if
\begin{equation*}
\disp\Big[\sum_{i=1}^m\lambda_i\nabla g_i(\ox)=0,\;\lm_i\ge 0\Big]\Longrightarrow\big[\lambda_i=0,\;i=1,\ldots,m\big],
\end{equation*}\vspace*{-0.25in}

\begin{proposition}\label{Th:co-nc} Consider the set $\O$ from \eqref{omega}, where $g:=(g_!,\ldots,g_m)$ is $C^2$-smooth  around $\ox\in\O$, and suppose that the following two constraint qualifications are satisfied:\vspace*{-0.05in}

{\bf(i)} The vectors $\nabla g_1(\ox),\ldots$ $\nabla g_m(\ox)$ are positively linearly independent.

{\bf(ii)} Given a normal $\ov\in N_\O(\ox)$, the multifunction $M\colon\R^{2m}\tto\R^{n+m}$ defined by
\begin{eqnarray}\label{calm}
M(\vTh):=\big\{(x,\lm)\;\big|\;\big(-g(x),\lm\big)+\vTh\in\gph N_{\R^m_-}\big\}
\end{eqnarray}
is calm \eqref{e:calm} at $(0,\ox,\bar\lm)$ for all $\bar\lm=(\bar\lm_1,\ldots,\bar\lm_m)\ge 0$ satisfying the equation $-\nabla g(\ox)^*\bar\lm=\bar v$.\\\vspace*{0.03in}
Then we have the second-order upper estimate
\begin{equation}
\label{e:cd1}
D^*N_\O(\ox,\bar v)(u)\subset\bigcup_{\bar\lm\ge 0,-\nabla g(\ox)\bar\lm=\bar v}\left\{\left(-\sum^m_{i=1}\bar\lm_i\nabla^2g_i(\ox)\right)u-\nabla g(\ox)^*D^*N_{\R^m_-}\big(-g(\ox),\bar\lm\big)\big(-\nabla g(\ox)u\big)\right\}.
\end{equation}
Moreover, \eqref{e:cd1} holds as an equality if the Jacobian $\nabla g(\ox)$ is of full rank when $\bar\lm\ge 0$ is a unique solution to $-\nabla g(\ox)^*\bar\lm=\bar v$ and both assumptions in {\rm(i)} and {\rm(ii)} are satisfied. If in the latter case $\O=\R^m_-$, then we get
\begin{equation}\label{e:cor-nc-o}
D^*N_{\R^m_-}(x,v)(y)=\left\{\begin{array}{ll}
\emp\mbox{ if }\;\mbox{ there is }\;i\;\mbox{ with }\;v_iy_i\not=0,\\
\{\gg\;|\;\gg_i=0\;\mbox{ for all }\;i\in I_1(y),\;\gg_i\ge 0\;\mbox{ for all }\;i\in I_2(y)\big\}\;\mbox{otherwise}
\end{array}\right.
\end{equation}
whenever $(x,v)\in\gph N_{\R^m_-}$ with the index subsets in \eqref{e:cor-nc-o} defined by
\begin{equation}\label{e:indexes}
I_1(y):=\big\{i\;\big|\;x_i<0\big\}\cup\big
\{i\;\big|\;v_i=0,\;y_i<0\big\},\quad I_2(y):=\big
\{i\;\big|\;x_i=0,\;v_i=0,\;y_i>0\big\}.
\end{equation}
\end{proposition}\vspace*{-0.05in}

The importance of the formula \eqref{e:cd1} is {\em transforming} the evaluation of the coderivative for normal cone mappings for general sets given in \eqref{omega} to the case of orthants \eqref{orthant}. The next proposition provide, suppose the crucial second-order computations of the coderivative of the mapping \eqref{F} associated with the weeping dynamics in terms of its given data; see \cite[Theorem~6.2]{cm3} for the detailed proof and discussions.\vspace*{-0.1in}

\begin{proposition}\label{Thm_Co}
Consider the set-valued mapping $F$ in \eqref{F} associated with the sweeping process \eqref{e:SP}, where the nonconvex set $C$ is taken from \eqref{e:MS}, and where the perturbation mapping that $f$ is $C^1$-smooth. Given $x,u\in\R^n$ with $x-u\in C$ as well as $w+f(x,a)\in N_C(x-u)$ and $a\in\R^d$, assume that the vectors $\nabla g_1(x-u),\ldots,\nabla g_m(x-u)$ are positively linearly independent and that the multifunction $M$ from \eqref{calm} is calm \eqref{e:calm} at $(0,x-u,\lm)$ for all $\lm=(\lm_1,\ldots,\lm_m)\ge 0$ satisfying the equation $-\nabla g(x-u)^*\lm=w+f(x,a)$. Then we have the coderivative upper estimate
\begin{equation}\label{e:co-upper-est1}
\begin{aligned}
D^*F(x,u,a,w)(y)\subset\disp\bigcup_{\lm\ge 0,-\nabla g(x-u)\lm= w+f(x,a)}\bigg\{\bigg(-\nabla_xf(x,a)^*y-\bigg(\sum^m_{i=1}\lm\nabla^2g_i(x-u)\bigg)y-\\
\nabla g(x-u)^*\gg,\disp\bigg(\sum^m_{i=1}\lm\nabla^2g_i(x-u)\bigg)y+\nabla g(x-u)^*\gg,-\nabla_af(x,a)^*y\bigg)\bigg\}
\end{aligned}	
\end{equation}
for all $y\in\dom D^*N_C\big(x-u,w+f(x,a)\big)$, where the coderivative domain satisfies the inclusion
\begin{equation}\label{e:dom-co-est}
\begin{aligned}
\dom D^*N_C(x-u,w+f(x,a))\subset&\big\{y|\;\exists\,\lm\ge0\;\mbox{ such that }\;-\nabla g(x-u)\lm=w+f(x,a),\\
&\lm_i\la\nabla g_i(x-u),y\ra=0\;\mbox{ for }\;i=1,\ldots,m\big\}
\end{aligned}
\end{equation}
with $\gg_i=0$ if either $g_i(x-u)>0$ or $\lm_i=0$, and with $\la\nabla g_i(x-u),y\ra>0$. We also get
$\gg_i\ge 0$ if $g_i(x-u)=0,\lm_i=0$ and that $\la\nabla g_i(x-u),y\ra<0$.
Moreover, replacing the calmness of \eqref{calm} by the stronger full rank assumption on the Jacobian matrix $\nabla g(x-u)$ ensures that both inclusions \eqref{e:dom-co-est} and \eqref{e:dom-co-est}
hold as equalities\} that the equalities with the collection of nonnegative multipliers $\lm=(\lm_1,\ldots,\lm_m)\ge 0$ uniquely determined by the equation $-\nabla g(x-u)^*\lm=w+f(x,a)$.
\end{proposition}\vspace*{-0.25in}

\section{Necessary Optimality Conditions via Discrete Approximations}\label{nec}
\setcounter{equation}{0}\vspace*{-0.1in}

In this section, we establish necessary optimality conditions for each discrete-time problem $(P_k)$ as $k\in\N$. Through this process and by applying Theorem~\ref{Th2a}, we derive suboptimal conditions for the continuous-time sweeping optimal control problem $(P)$ with any desired level of accuracy.\vspace*{-0.1in}

\begin{theorem}\label{Th8} Fix any $k\in\N$ and let
\begin{equation*}
\(\ox^k,\ou^k,\oa^k,\oT^k\)=\(\ox_0^k,\ldots,\ox_{k}^k,\ou_0^k,\ldots,\ou_{k}^k,\oa_0^k,\ldots,\oa_{k-1}^k,\oT^k\)
\end{equation*}
be an r.i.l.m.\ in  $(P_k)$ along which the assumptions of Theorem~{\rm\ref{Th2a}} together with general qualification 
conditions of Proposition~{\rm\ref{Th:co-nc}} are fulfilled. Assume in addition that the cost functions $\ph$ and $\ell_t:=\ell(t,\cdot,\cdot)$ are locally Lipschitzian around the optimal points for any $t\in\Delta_k(\T^k)$ defined in \eqref{e:DP}, and that $\ell_t$ does not depend on the time variable $t$. Then there exist dual elements
$\al^k=(\al^k_1,\ldots,\al^k_m)\in\R^m_+$, $\psi^{uk}_j =\( \psi^{uk}_{j1}, \ldots, \psi^{uk}_{js}\)\in \R^s_+,\; \psi^{ak}_j \in \R^d$ for $j=0,\ldots,k-1$, and $p^k_j=(p^{xk}_j,p^{uk}_j)\in\R^n\times\R^n$ for $j=0,\ldots,k$ satisfying the conditions
\begin{equation}
\label{e:nontri-da}
\lm^k+\|\al^k\|+\sum^{k}_{j=0}\|p^{xk}_j\|+\|p^{uk}_0\| + \sum^{k-1}_{j=0}\(\n\psi^{uk}_j\en +  \n \psi^{ak}_j\en\)\not=0,
\end{equation}
\begin{equation}\label{e:dac1}
\al^k_ig_i(\ox^k_{k}-\ou^k_{k})=0,\;\;i=1,\ldots,m,
\end{equation}
\begin{equation}\label{e:dac4}
p^{uk}_{k}=-\sum^m_{i=1}\al^k_i\nabla g_i(\ox^k_{k}-\ou^k_{k}),
\end{equation}
\begin{equation}\label{e:dac5}
\begin{array}{ll}
p^{uk}_{j+1}=\lm^k\(v^{uk}_j+\frac{\th^{Uk}_j}{h_k}\),\;\;j=0,\ldots,k-1,
\end{array}
\end{equation}
\begin{equation}
\label{e:dac6}
\begin{array}{ll}
&\bigg(\frac{p^{xk}_{j+1}-p^{xk}_j}{h_k}-\lm^kw^{xk}_j,\frac{p^{uk}_{j+1}-p^{uk}_j}{h_k}-\lm^kw^{uk}_j,-\lm^kw^{ak}_j -\frac{1}{h_k}\lm^k\th^{ak}_j, -p^{xk}_{j+1}+\lm^k\Big(v^{xk}_j+ \frac{\th^{Xk}_j}{h_k}\Big)\bigg) \\
&\in N\bigg(\bigg(\ox^k_j,\ou^k_j,\oa^k_j,\frac{\ox^k_{j+1}-\ox^k_j}{-h_k}\bigg);\gph F\bigg) + \bigg(0,\disp\sum^s_{i=1}\frac{\psi^{uk}_{ji}}{h_k}\nabla\nu_i(\ou^k_j), \frac{\psi^{ak}_j}{h_k},0\bigg),
\end{array}
\end{equation}
\begin{equation}
\label{e:dac7}
\begin{aligned}
&\(-p^{xk}_k  +\sum^m_{i=1}\al^k_i\nabla g_i(\ox^k_{k}-\ou^k_{k}), \H^k +\lm^k(\T-\T^k)+ \lm^k\varrho^k \)\\
&\in  \lm^k \partial\varphi\(\ox^k_{k}, \T^k\) +  N\Big(\big(\ox_{k},\T\big); \Xi^k_x\times \Xi^k_T \Big)
\end{aligned}
\end{equation}
for $j=0,\ldots,k-1$ with the triples
\begin{equation}\label{e:dac9}
\begin{array}{ll}
(\th^{Xk}_j,\th^{Uk}_j, \th^{ak}_j):={\disp\int}^{t^k_{j+1}}_{t^k_j}\left(\frac{\ox^k_{j+1}-\ox^k_j}{h_k}-\dot{\ox}(t),\frac{\ou^k_{j+1}-\ou^k_j}{h_k}-\dot{\ou}(t), \oa^k_j-\oa(t)\right)dt
\end{array}
\end{equation}
along with the collection of subgradients
\begin{equation}
\label{e:dac9a}
\begin{array}{ll}
(w^{xk}_j,w^{uk}_j,w^{ak}_j,v^{xk}_j,v^{uk}_j)\in\partial_{x,u,a,X,U}\ell\left(t^k_j,\ox^k_j,\ou^k_j, \oa^k_j,\frac{\ox^k_{j+1}-\ox^k_j}{h_k},\frac{\ou^k_{j+1}-\ou^k_j}{h_k}\right),\;\;j=0,\ldots,k-1,
\end{array}
\end{equation}
where the sequence $\{\delta_k\}\dn0$ as $k\to\infty$ is taken from in
Theorem~{\rm\ref{Th1}}, and where
$$
h_k\colon = \T^k/k;\;\; t^k_j\colon = jh_k  \mbox{ for } j=0,\ldots, k,
$$
\begin{equation}
\label{e:Hk}
\begin{array}{ll}
\H^k: =\frac{1}{k}{\disp\sum}^{k-1}_{j=0}\[\la p^{xk}_{j+1}, \frac{\ox^k_{j+1}-\ox^k_j}{h_k}\ra +\la p^{uk}_{j+1}, \frac{\ou^k_{j+1}-\ou^k_j}{h_k}\ra -\lm^k\ell\( t^k_j,\ox^k_j,\ou^k_j,\oa^k_j,\frac{\ox^k_{j+1}-\ox^k_j}{h_k}, \frac{\ou^k_{j+1}-\ou^k_j}{h_k} \)  \],
\end{array}
\end{equation}
\begin{equation}
\label{e:varrho}
\begin{array}{ll}
\varrho^k&: =-  {\disp\sum}^{k-1}_{j=0}\bigg[\frac{j+1}{k} \bigg\|\(\oa^k_j,  \frac{\ox^k_{j+1}-\ox^k_j}{h_k},\frac{\ou^k_{j+1}-\ou^k_j}{h_k} \)- \bigg(\oa\(t^k_{j+1}\),\dot{\ox}\(t^k_{j+1}\), \dot{\ou}\(t^k_{j+1}\) \bigg)\bigg\|^2 \\
&-\frac{j}{k}\n\(\oa^k_j,  \frac{\ox^k_{j+1}-\ox^k_j}{h_k},\frac{\ou^k_{j+1}-\ou^k_j}{h_k} \)-\(\oa\(t^k_j\),\dot{\ox}\(t^k_j\),\dot{\ou}\(t^k_j\) \)\en^2\bigg],
\end{array}
\end{equation}
\begin{equation}\label{psi}
\psi^{uk}_{ij}\in N\(\nu_i(\ou^k_j);(-\infty,L_\nu\delta_k]\),\;\;i=1, \ldots, s, \; j=0\ldots,k;
\end{equation}
\begin{equation}\label{rho}
\psi^{ak}_j\in N\(\oa^k_j;A\),\;\;j=0\ldots,k-1.
\end{equation}
\end{theorem}
{\bf Proof}. Fix $k\in \N$, and let $T^k\geq 0$ be a given constant. Define the partition of the interval $[0,T^k]$ by
$$
\Delta_k(T^k):=\left\{0=t^k_0<t^k_1<\ldots<t^k_{k-1}<t^k_{k}=T^k\right\}
$$
with the uniform stepsize $h_k:=t^k_{j+1}-t^k_j = \frac{T^k}{k}$ for all $j=0,\ldots,k-1$. Denote
$$
y^k:=(x^k_0,\ldots,x^k_{k},u^k_0,\ldots,u^k_{k},a^k_0,\ldots,a^k_{k-1} ,X^k_0,\ldots,X^k_{k-1},U^k_0,\ldots,U^k_{k-1},T^k).
$$
Now we reformulate $(P_k)$ as an equivalent mathematical programming problem $(MP)$ with  the variable $y^k$ and the fixed initial point $x_0$ by
\begin{eqnarray*}
\begin{array}{ll}
\mbox{minimize }\;\varphi_{0}[y^k]&:=\varphi\(x^k_{k},T^k\)+\disp\frac{T^k}{k}{\sum\limits}^{k-1}_{j=0} \ell(t^k_j,x^k_j,u^k_j,a^k_j,X^k_j,U^k_j)+\frac{1}{2}(T^k-\T)^2\\
&+\disp\frac{1}{2}\sum\limits^{k-1}_{j=0}{\disp\int^{(j+1)T^k/k}_{jT^k/k}}\n \l( a^k_j,X^k_j,U^k_j\r) - \l(\oa^k(t),\dot \ox^k(t), \dot \ou^k(t)\r) \en^2 dt\\
&+\dist^2\l(\l\| \frac{ u^k_1-u^k_0}{h_k}\r\|;(-\infty,\Tilde\mu ]\r) + \disp\dist^2\Big(\disp\disp\sum^{k-2}_{j=0}\Big\|U^k_{j+1}-U^k_j\Big\|; (-\infty,\Tilde \mu]\Big)
\end{array}
\end{eqnarray*}
constrained by the following equality, inequality, and geometric conditions:
\begin{equation*}
\begin{array}{ll}
b^x_j(y^k):= x^k_{j+1}-x^k_j-(T^k/k) X^k_j=0\;\mbox{ for }\;j=0,\ldots,k-1,\\
b^u_j(y^k):=u^k_{j+1}-u^k_j- (T^k/k) U^k_j=0\;\mbox{ for }\;j=0,\ldots,k-1,\\
c_i(y^k):=-g_i(x^k_{k}-u^k_{k})\le 0\;\mbox{ for }\;i=1,\ldots,m,\\
\phi_{j}(y^k):=\sum\limits^{k-1}_{j=0}{\disp\int^{t^k_{j+1}}_{t^k_j}}\left\|(x^k_j,u^k_j,a^k_j)-(\ox^k(t^k_j),\ou^k(t^k_j),\oa^k(t^k_j))\right\|^2 -\frac{\ve}{2}\le 0,\\
\phi_{k+1}(y^k):=\sum\limits^{k-1}_{j=0}{\disp\int^{t^k_{j+1}}_{t^k_j}}\n \l( a^k_j,X^k_j,U^k_j\r) - \l(\oa^k(t),\dot \ox^k(t), \dot \ou^k(t)\r) \en^2 dt -\frac{\ve}{2}\le 0,\\
\phi_{k+2}(y):=\sum\limits^{k-2}_{j=0}\l\| U^k_{j+1}-U^k_j\r\| -\Tilde\mu-1\leq 0,\\
\phi_{k+3}(y):= \l\| u^k_1-u^k_0\r\|\leq (\Tilde\mu +1)h_k, 
\\
y^k\in\O^X_{j}:=\big\{(x^k_0, \ldots, U^k_{k-1},T^k)\big{|}\;-X^k_{j}\in F(x^k_{j},u^k_{j},a^k_{j})\big\}\;\mbox{ for }\;j=0,\ldots,k-1,\\
y^k\in\O^{ua}_{j}:=\big\{(x^k_0, \ldots, U^k_{k-1},T^k)\;\big|\;\(u^k_j,a^k_j\)\in U_k \times A\big\}, \; \mbox{ for }\; j=0,\ldots,k-1, \\
y^k\in\O_{bd}:= \big\{(x^k_0, \ldots, U^k_{k-1},T^k)\;\big|\;\big((x^k_0,u^k_0\big)=\big(x_0,\ou(0)\big)\;\mbox{ and }\;(x^k_{k}, T_k) \in \Xi^k_x \times \Xi^k_T\big\}.
\end{array}
\end{equation*}
Next we apply the necessary optimality conditions from \cite[Theorem~3.5(ii)]{m18} to any local optimal solution  $\oy^k$ of problem $(MP)$ by taking into account that all the inequality constraints in $(MP)$ associated with the functions $\phi_j$ as $j=0,\ldots,k+1$ become {\em inactive} for sufficiently large $k\in\N$. As a result, the corresponding multipliers vanish from the optimality conditions. This leads us to dual elements $\lm^k \ge 0,\;\alpha^k=(\al^k_1,\ldots,\al^k_m)\in\R^m_+,\;p^k_j=(p^{xk}_j,p^{uk}_j)\in\R^{2n}$, $\psi^k_j =(\psi^{uk}_j,\psi^{ak}_j)\in \R^{n}\times\R^d$ as $j=0,\ldots,k$, and
$$
y^\ast_j=(x^\ast_{0j},\ldots,x^\ast_{kj},u^\ast_{0j},\ldots,u^\ast_{kj},a^\ast_{0j},\ldots,a^\ast_{(k-1)j},X^\ast_{0j},\ldots,X^\ast_{(k-1)j},U^\ast_{0j},\ldots,U^\ast_{(k-1)j}, T^\ast_j),
$$
which are not zero simultaneously while satisfying the conditions in \eqref{e:dac4} and the inclusions
\begin{equation}\label{69}
y^\ast_j\in\left\{
\begin{array}{llll}
N_{\O^X_j}(\oy^k)+N_{\O^{ua}_j}(\oy^k)&\textrm{ if }&j\in\{0,\ldots, k-1\}\\
N_{\O_{bd}}(\oy^k)&\textrm{ if }& j=k
\end{array}
\right.,
\end{equation}
\begin{equation}\label{e:dac11}
-y^\ast_0-\ldots-y^\ast_{k}\in\lm^k\partial\varphi_{0}(\oy^k)+\sum^m_{i=1}\alpha^k_i\nabla c_i(\oy^k)+\sum^{k-1}_{j=0}\nabla b^x_j(\oy^k)^\ast p^{xk}_{j+1}+\sum^{k-1}_{j=0}\nabla b^u_j(\oy^k)^\ast p^{uk}_{j+1},
\end{equation}
\begin{equation*}
\al^k_i c^k_i(\oy^k)=0\;\mbox{ for }\;i=1,\ldots,m.
\end{equation*}
Observe that the first line in $(\ref{69})$ is derived by applying the normal cone intersection rule from \cite[Theorem~2.16]{m18} to $\oy\in\O^X_j\cap\O^{ua}_j$ for
$j=0,\ldots,k-1$ provided that the qualification condition
\begin{equation}\label{qc}
N_{\O^X_j}(\oy^k)\cap\big(-N_{\O^{ua}_j}(\oy^k)\big)=\{0\},\quad j=0,\ldots,k-1,
\end{equation}
holds. To verify the latter, pick any vector $y^\ast_j\in
N_{\O^X_j}(\oy^k)\cap(-N_{\O^{ua}_j}(\oy^k))$ and get the inclusions
\begin{equation}\label{qc1}
(x^\ast_{jj},u^\ast_{jj},a^\ast_{jj},-X^\ast_{jj})\in
N_{\gph F}\(\ox^k_j,\ou^k_j,\oa^k_j,-\frac{\ox^k_{j+1}-\ox^k_j}{h_k}\),
\end{equation}
$$-u^\ast_{jj}\in N_{U_k}(\ou^k_j),\,-a^\ast_{jj}\in N_A(\oa^k_j),
$$
while the remaining components of $y^\ast_j$ are zero. It follows from the second inclusion in \eqref{qc1} that
\begin{equation*}
x^\ast_{jj}=0\;\mbox{ and }\;X^\ast_{jj}=0.
\end{equation*}
Substituting this into the first inclusion in \eqref{qc1} and applying the coderivative definition \eqref{e:cor} give us
\begin{equation*}
(0,u^\ast_{jj},a^\ast_{jj})\in
D^\ast F\Big(\ox^k_j,\ou^k_j,\oa^k_j,-\frac{\ox^k_{j+1}-\ox^k_j}{h_k}
\Big)(0),\quad j=0,\ldots,k-1.
\end{equation*}
Directly from the coderivative estimate \eqref{e:co-upper-est1} for
$F$ in \eqref{F}, we deduce under the imposed PLICQ assumption that $u^\ast_{jj}=0,\,a^\ast_{jj}=0$ for
all $j=0,\ldots,k-1$. As a result, for such indices $j$, we have $y^\ast_j=0$, confirming that the qualification condition \eqref{qc} holds. Consequently, the inclusions in \eqref{69} are equivalent to
\begin{equation}\label{e:5.18*}
\begin{cases}
(x^\ast_{jj},u^\ast_{jj}-\Psi^{uk}_{j},a^\ast_{jj}-\psi^{ak}_j,\,-X^\ast_{jj})\in
N_{\gph F}\Big( \ox^k_j,\ou^k_j,\oa^k_j,-\frac{\ox^k_{j+1}-\ox^k_j}{h_k}
\Big)  \textrm{ for }\;j=0,\ldots,k-1,\\
\big(x^\ast_{kk}, T^\ast_{k} \big) \in {N_{\Xi^k_x \times \Xi^k_T }\(\ox^k_{k},\T\)}, \\
\big(x^\ast_{ij}, u^\ast_{ij}, a^\ast_{ij}, X^\ast_{ij}, T^\ast_j \big) = (0, 0, 0, 0, 0) \mbox{ if } i\not=j \in\{0,\ldots,k-1\}, \\
\big(x^\ast_{i k}, u^\ast_{j k}, a^\ast_{j k} \big) = (0,0,0,0,0) \mbox{ for } i = 1,\ldots, k-1,\; j=1, \ldots, k, \\
\big( X^\ast_{i k}, U^\ast_{i k}\big)  = (0,0) \mbox{ for } i = 0,\ldots, k-1, \\
a^\ast_{0k} = 0,
\end{cases}
\end{equation}
with $\Psi_{j}^{uk}\in N_{U_k}(\ou^k_j)$ and $\psi^{ak}_j$ taken from \eqref{rho}. Using assumption (H3.2) and the calculus rules for the normal cones from \cite[Theorem~1.17]{m-book} allows us to represent $N_{U_k}(\ou^k_j)$ as follows:
$$
N_{U_k}(\ou^k_j) = \sum^s_{i=1}\nabla\nu_i\(\ou^k_j\)^\ast N_{(\infty,L_\nu\delta_k]}\(\nu_i(\ou^k_j)\) \mbox{ for } j=0,\ldots, k-1,
$$
where the functions $\nu_i$ are given in \eqref{U}. As a consequence, we get
$$
\begin{cases}
\Psi^{uk}_j= \disp\sum^s_{i=1}\psi^{uk}_{ij}\nabla\nu_i\(\ou^k_j \) \mbox{ for } j=0,\ldots, k-1;\\
u^\ast_{kk} = \disp\sum^s_{i=1}\psi^{uk}_{ik}\nabla\nu_i\(\ou^k_k \)
\end{cases}
$$
with the nonnegative numbers $\psi^{uk}_{ij}$ taken from \eqref{psi}. In a similar manner, the vectors $\big(x^\ast_{k}, u^\ast_{0 k},x^\ast_{kk}, T^\ast_{k}\big)$, arising from the normal cone $N_{\O_{bd}}(\oy^k)$, are the only nonzero components of $y^\ast_{k}$. Consequently, this yields
\begin{equation}\label{e:dac13}
\begin{aligned}
&-y^\ast_0-y^\ast_1-\ldots-y^\ast_{k}=\big(\underbrace{-x^\ast_{0k}-x^\ast_{00},-x^\ast_{11},\ldots,-x^\ast_{k-1,k-1}, -x^\ast_{kk} }_{x-\mbox{components}}, \\
&\underbrace{-u^\ast_{0k}-u^\ast_{00},\ldots,  -u^\ast_{k-1,k-1}, 0}_{u-\mbox{components}}, 
\underbrace{-a^\ast_{00},-a^\ast_{11},\ldots,-a^\ast_{k-1,k-1},0,}_{a-\mbox{components}} \\
&\underbrace{-X^\ast_{00},\ldots,-X^\ast_{k-1,k-1}}_{X-\mbox{components}},\underbrace{0,\ldots,0}_{U-\mbox{components}}, -T^\ast_{k}\big).
\end{aligned}
\end{equation}
Next we calculate the expression on the right-hand side of \eqref{e:dac11}. It follows that
\begin{eqnarray*}
\begin{aligned}
\left(\sum^m_{i=1}\alpha^k_i\nabla c_i(\oy^k)\right)_{(x^k_{k},u^k_{k},a^k_{k})}=&\left(-\sum^m_{i=1}\alpha_i^k\nabla g_i(\ox^k_{k}-\ou^k_{k}), \sum^m_{i=1}\alpha_i^k\nabla g_i(\ox^k_{k}-\ou^k_{k}), 0\right),\\
\left(\sum^{k-1}_{j=0}\big(\nabla b_j(\oy^k)\big)^\ast p^k_{j+1}\right)_{(x^k_j,u^k_j)}=&
\begin{cases}
-p^k_1 & \mbox{ if  }\; j=0,\\
p^k_j-p^k_{j+1} & \mbox{ if  }\; j=1,\ldots,k-1,\\
p^k_{k} & \mbox{ if  }\; j=k,
\end{cases}\\
\left(\sum^{k-1}_{j=0}\big(\nabla b_j(\oy^k)\big)^*p^k_{j+1}\right)_{(X,U)}=&\big(-h_kp^{xk}_1,\ldots,-h_kp^{xk}_{k},-h_kp^{uk}_1,\ldots,-h_kp^{uk}_{k}\big),
\end{aligned}
\end{eqnarray*}
$$
\begin{array}{ll}
\sum\limits^{k-1}_{j=0}\(\nabla_{\T^k} b^x_j(\oy^k)\)^*p^{xk}_{j+1} &= -\sum\limits^{k-1}_{j=0}\frac{1}{k} \la X^k_j,p^{xk}_{j+1}\ra, \\
\sum\limits^{k-1}_{j=0}\(\nabla_{\T^k} b^u_j(\oy^k)\)^*p^{uk}_{j+1} &= -\sum\limits^{k-1}_{j=0}\frac{1}{k} \la U^k_j,p^{uk}_{j+1}\ra.
\end{array}
$$
Applying the subdifferential sum rule from \cite[Theorem~2.19]{m18} brings us to the inclusion
$$
\begin{array}{ll}
&\partial\varphi_0(\oy^k)\subset \partial\varphi\(\ox^k_{k}, \T^k\)+ \frac{\T^k}{k}\sum\limits^{k-1}_{j=0} \partial\ell\( \frac{jT^k}{k} ,\ox^k_j,\ou^k_j,\oa^k_j,\bar{X}^k_j,
\bar{U}^k_j\)\\[2ex]
&+\sum\limits^{k-1}_{j=0}\nabla\zeta _j(\oy^k)+ (0,\ldots,0, \T^k-\T) + \partial\sigma(\oy^k),
\end{array}
$$ 
where $\zeta _j(\cdot)$ and $\sigma(\cdot)$ are real-valued functions defined by
$$
\zeta _j(y^k):=\frac{1}{2}\int^{(j+1)\T^k/k}_{j\T^k/k}\left\|(a^k_j, X^k_j,U^k_j )-(\oa(t),\dot{\ox}(t),\dot{\ou}(t) )\right\|^2dt,
$$
$$
\sigma(y^k)\colon = \dist^2\l(\l\| \dfrac{u^k_1-u^k_0}{h_k}\r\|;(-\infty,\Tilde\mu ]\r) + \mbox{dist}^2\bigg(\sum^{k-1}_{j=1}\l\|U^k_{j+1}-U^k_j\r\| ; (-\infty, \Tilde\mu] \bigg).
$$
The differentiability of $\psi(x)\colon = \dist^2(x;(-\infty,\Tilde\mu])$ with the gradient $\nabla\psi (x)=0$ whenever $x\leq \Tilde\mu$ together with \eqref{e:conv4} implies that $\partial\sigma(\oy^k)=\{0\}$. Differentiating $\zeta _j(\cdot)$ with respect to $a^k_j, X^k_j, U^k_j$ and $T^k$ gives us
\begin{equation}
\label{e:grad}
\begin{cases}
\nabla_{a^k_j,X^k_j,U^k_j}\zeta_j(\oy^k)=(\theta^{ak}_j,\theta^{Xk}_j,\theta^{Uk}_j), \\[2ex]
\sum\limits^{k-1}_{j=0}\nabla_{T^k}\zeta_j(\oy^k) = \zeta_{T^k}:= \sum\limits^{k-1}_{j=0}\bigg[\frac{j+1}{k} \bigg\|\(a^k_j, X^k_j,U^k_j \)- \bigg(\oa\(\frac{(j+1)\T^k}{k}\),\dot{\ox}\(\frac{(j+1)\T^k}{k}\),
\dot{\ou}\(\frac{(j+1)\T^k}{k}\) \bigg)\bigg\|^2 \\[2ex]
- \frac{j}{k}\left\|\(a^k_j, X^k_j,U^k_j \)-\(\oa\(j\T^k/k\),\dot{\ox}\(j\T^k/k\),\dot{\ou}\(j\T^k/k\) \)\right\|^2\bigg],
\end{cases}
\end{equation}
where the triples $(\theta^{Xk}_j,\theta^{Uk}_j,\theta^{ak}_j)$ are given in \eqref{e:dac9}. Combining this with the above inclusion we obtain the following upper estimate for $\lm^k\partial\varphi_0(\oy^k)$ in \eqref{e:dac11}:
$$
\begin{array}{ll}
\lm^k \bigg(\underbrace{h_k w^{xk}_0, h_kw^{xk}_1,\ldots,h_kw^{xk}_{k-1},\vartheta^{xk}_{k}}_{x-\mbox{components}}, 
\underbrace{h_kw^{uk}_0, h_kw^{uk}_1, \ldots ,h_kw^{uk}_{k-1}, 0,}_{u-\mbox{components}} \\
\underbrace{\th^{ak}_0 + h_kw^{ak}_0, \th^{ak}_1 + h_kw^{ak}_1, \ldots, \th^{ak}_{k-1}+h_kw^{ak}_{k-1}, 0,}_{a-\mbox{components}} \\
 \underbrace{\theta^{Xk}_0+h_kv^{xk}_0, \theta^{Xk}_1+h_kv^{xk}_1, \ldots, \theta^{Xk}_{k-1}+h_kv^{xk}_{k-1}}_{X-\mbox{components}},\\
\underbrace{\theta^{Uk}_0+h_kv^{uk}_0, \theta^{Uk}_1+h_kv^{uk}_1, \ldots,\theta^{Uk}_{k-1}+h_kv^{uk}_{k-1}}_{U-\mbox{components}}, \\\zeta_{T^k}+ \vartheta^{Tk}_k+\T^k-\T +  \frac{1}{k}\sum\limits^{k-1}_{j=0} \ell(t^k_j,x^k_j,u^k_j,a^k_j,X^k_j,U^k_j) \bigg),
\end{array}
$$
where $\(\vartheta^{xk}_{k}, \vartheta^{Tk}_k \)\in \partial\varphi\(\ox^k_{k}, \T^k\)$ and the components of $(w^{xk},w^{uk},w^{ak},v^{xk},v^{uk})$ satisfying \eqref{e:dac7}. Incorporating \eqref{e:dac13} together with \eqref{69} ensures that \vspace*{-0.05in}
\begin{align}
-x^\ast_{0k}-x^\ast_{00 }&=\lambda^k h_k w^{xk}_0-p^{xk}_1, \label{e:dx1}\\
-x^\ast_{jj} &=\lambda^k h_k w^{xk}_j+p^{xk}_j-p^{xk}_{j+1}\;\mbox{ for }\;j=1,\ldots, k-1, \label{e:dx2} \\
-x^\ast_{kk}&=\lambda^k\vartheta^{xk}_{k} - \sum^m_{i=1}\alpha_i^k\nabla g_i(\ox^k_{k}-\ou^k_{k}) + p^{xk}_{k}, \label{e:dx3}\\
-u^\ast_{0k}-u^\ast_{00} &= \lambda^k h_k w^{uk}_0 - p^{uk}_1, \label{e:du1}\\
-u^\ast_{jj} &=\lambda^k h_kw^{uk}_j+p^{uk}_j-p^{uk}_{j+1}\;\mbox{ for }\;j=1,\ldots,k-1,\label{e:du2}\\ 
0 &= \sum^m_{i=1}\alpha^k_i\nabla g_i(\ox^k_{k}-\ou^k_{k}) +p^{uk}_{k},\label{e:du3}\\
-a^\ast_{00} &=\lambda^k (h_k w^{ak}_0 + \th^{ak}_0), \label{e:da1}\\
-a^\ast_{jj} &=\lambda^k(h_k w^{ak}_j+ \th^{ak}_j) \;\mbox{ for }\;j=1,\ldots, k-1, \label{e:da2}\\[2ex]
-X^\ast_{jj} &=\lambda^k(h_k v^{xk}_j+\theta^{Xk}_j)-h_k p^{xk}_{j+1}\;\mbox{ for }\;j=0,\ldots,k-1, \label{e:dX} \\
0&=\lambda^k(h_k v^{uk}_j+\theta^{Uk}_j)-h_k p^{uk}_{j+1}\;\mbox{ for }\;j=0,\ldots,k-1,\label{e:dU}\\
-T^\ast_{k} &= \lm^k\( \zeta_{T^k}+ \vartheta^{Tk}_k + \T^k-\T+ \frac{1}{k}\sum^{k-1}_{j=0} \ell(t^k_j,x^k_j,u^k_j,a^k_j,X^k_j,U^k_j) \)\label{e:dT}\\
& -\sum^{k-1}_{j=0}\frac{1}{k} \la X^k_j,p^{xk}_{j+1}\ra -\sum^{k-1}_{j=0}\dfrac{1}{k} \la U^k_j,p^{uk}_{j+1}\ra\nonumber,
\end{align}
where $\zeta_{T^k}$ is defined in \eqref{e:grad}.
This enables us to justify all the necessary optimality conditions claimed in this theorem. It is obvious that \eqref{e:dac4} follows from \eqref{e:dx3} and \eqref{e:du3} and that \eqref{e:dac5} follows from \eqref{e:dU}.\vspace*{-0.05in}

Next we extend each vector $p^k = (p^{xk},p^{uk})$ by adding the zero component $p^k_0 = \left(x^\ast_{0k}, u^\ast_{0k}\right)$. Rewrite equations \eqref{e:dx2}, \eqref{e:du2},  and \eqref{e:da2} as follows:
$$
\begin{array}{ll}
\frac{x^\ast_{jj}}{h_k} &= \frac{p^{xk}_{j+1}-p^{xk}_j}{h_k} - \lm^kw^{xk}_j \;\mbox{ for }\;j=0,\ldots,k-1, \\
\frac{u^\ast_{jj}}{h_k} &= \frac{p^{uk}_{j+1}-p^{uk}_j}{h_k} - \lm^kw^{uk}_j \;\mbox{ for }\;j=0,\ldots,k-1,\\
\frac{a^\ast_{jj}}{h_k} &= -\lm^kw^{ak}_j - \frac{1}{h_k}\lm^k\th^{ak}_j \;\mbox{ for }\;j=0,\ldots,k-1,\\
\frac{X^\ast_{jj}}{h_k} &= -\lm^kv^{xk}_j-\frac{1}{h_k}\lm^k\th^{Xk}_j + p^{xk}_{j+1} \;\mbox{ for }\;j=0,\ldots,k-1.\\
\end{array}
$$
Combining this with the first inclusion in \eqref{e:5.18*} allows us to arrive at \eqref{e:dac6}. Moreover, the inclusion in \eqref{e:dac7} follows from that in \eqref{e:dT} and the second inclusion of \eqref{e:5.18*}.\vspace*{-0.05in}

To finish the proof of the theorem, it remains to justify the nontriviality condition \eqref{e:nontri-da}. Suppose on the contrary that $\lm^k=0, \;\al^k=0,\; p^{uk}_0=0, \; p^{xk}_j=0 $, and $\psi^{k}_j=0$ for $j=0, \ldots, k$.  This leads to $x^\ast_{0k} = p^{xk}_0=0, \; u^\ast_{0k}=p^{uk}_0=0$.  Consequently, it follows from \eqref{e:dac5} that $p^{uk}_j=0$ for $j=0, \ldots, k$. Moreover, applying \eqref{e:dx2}, \eqref{e:dx3}, \eqref{e:du2}, \eqref{e:da2}, \eqref{e:dX}, and \eqref{e:dT} tells us that $(x^\ast_{jj}, u^\ast_{jj}, a^\ast_{jj}, X^\ast_{jj}, T^\ast_k)=0$ for $j=0, \ldots, k$. Combining the latter with \eqref{e:5.18*}, observe that all the  components of $y^\ast_j$ are identically zero for $j=0, \ldots, k$,  which is a contradiction completing the proof.$\h$\vspace*{-0.03in}

Observe that the discrete-time Euler-Lagrange inclusion in \eqref{e:dac6} are expresses in terms of the normal cone to the graph (i.e., coderivative) of the velocity mapping \eqref{F} associated with the sweeping dynamics. The next theorem incorporates the coderivative evaluation for $F$ given in Proposition~\ref{Th:co-nc}.\vspace*{-0.12in}

\begin{theorem}\label{Th9}
Let $(\ox^k,\ou^k,\oa^k, \T^k)$ be an optimal solution to discrete-time problem $(P_k)$ with any fixed $k\in\N$ and with the sweeping velocity mapping $F$ defined in \eqref{F}. Suppose
that the functions $g_i$ in \eqref{e:MS} are of class ${\cal C}^2$ and the perturbation mapping $f(\cdot,a)$ is of class $C^1$ around the optimal points. Then there exist dual elements $\lm^k,\; \al^k, \;p^k,\;\psi^{uk},\; \psi^{ak}$ from Theorem~{\rm\ref{Th8}} together with vectors $\eta^k_j\in\R^m_+$ as $j=0,\ldots,k-1$ and $\gg^k_j\in\R^m$ as $j=0,\ldots,k-1$ such that the following conditions are satisfied:\\
{\sc Nontriviality condition}
\begin{equation}\label{e:dac14}
\lm^k+\|\al^k\|+\sum^{k}_{j=0}\|p^{xk}_j\|+\|p^{uk}_0\| + \sum^{k-1}_{j=0}\| \psi^k_j\|\not=0.
\end{equation}
{\sc Primal-dual dynamic relationship} for all $j=0,\ldots,k-1$:
\begin{equation}\label{e:dac15}
\begin{array}{ll}
\frac{\ox^k_{j+1}-\ox^k_j}{h_k}-f(\ox^k_j,\oa^k_j)=\disp\sum_{i\in I(\ox^k_j-\ou^k_j)}\eta^k_{ji}\nabla g_i(\ox^k_j-\ou^k_j),
\end{array}
\end{equation}
\begin{equation}\label{e:dac16}
\begin{array}{ll}
&\frac{p^{xk}_{j+1}-p^{xk}_j}{h_k}-\lm^kw^{xk}_j= -\nabla_xf(\ox^k_j,\oa^k_j)^*(\Lm^k_j) - \xi^k_j -\disp\sum^m_{i=1}\gg^k_{ji}\nabla g_i(\ox^k_j-\ou^k_j),
\end{array}
\end{equation}
\begin{equation}\label{e:dac17}
\begin{array}{ll}
\frac{p^{uk}_{j+1}-p^{uk}_j}{h_k}-\lm^kw^{uk}_j - \frac{1}{h_k}\disp\sum^s_{i=1}\psi^{uk}_{ji}\nabla\nu_i(\ou^k_j)&=\xi^k_j +\disp\sum^m_{i=1}\gg^k_{ji}\nabla g_i(\ox^k_j-\ou^k_j),
\end{array}
\end{equation}   
\begin{equation}
\label{e:dac18}
\begin{array}{ll}
&-\lm^kw^{ak}_j -\frac{1}{h_k}\lm^k\th^{ak}_j -\frac{\psi^{ak}_j}{h_k}=-\nabla_af(\ox^k_j,\oa^k_j)^*(\Lm^k_j),
\end{array}
\end{equation}
\begin{equation}
\label{e:dac18a}
\begin{array}{ll}
p^{uk}_{j+1}=\lm^k\Big(v^{uk}_j+\frac{\th^{Uk}_j}{h_k}\Big)
\end{array}
\end{equation}
with $(w^{xk}_j,w^{uk}_j,w^{ak}_j,v^{xk}_j,v^{uk}_j,v^{ak}_j)$ taken from \eqref{e:dac9a}, while the active constraint index set $I(\cdot)$ and the triples $(\th^{Xk}_j,\th^{Uk}_j, \th^{ak}_j)$  are specified in \eqref{e:a-index} and \eqref{e:dac9}, respectively. Additionally, $\Lm^k_j$ and $\xi^k_j$ are defined by
\begin{equation}
\label{e:dac18c}
\begin{cases}
\Lm^k_j:= \lm^k(h^{-1}_k\th^{Xk}_j+v^{xk}_j)-p^{xk}_{j+1},\\
\xi^k_j:= \disp\l(\sum^m_{i=1}\eta^k_{ji}\nabla^2g_i(\ox^k_j-\ou^k_j)\r)\Lm^k_j.
\end{cases}
\end{equation}
{\sc Adjoint Inclusions:} $\psi^{uk}_j$ and $\psi^{ak}_j$ satisfy \eqref{psi} and \eqref{rho}, respectively. 
\noindent {\sc Transversality conditions}:
\begin{equation}
\label{e:dac19}
\begin{cases}
&p^{uk}_k=-\sum\limits^m_{i=1}\al^k_{i}\nabla g_i(\ox^k_k-\ou^k_k),\\
&\(-p^{xk}_k+ \sum\limits^m_{i=1}\al^k_i\nabla g_i(\ox^k_{k}-\ou^k_{k}), \H^k +\lm^k(\T-\T^k)+ \lm^k\varrho^k \)\\
&\in  \lm^k \partial\varphi\(\ox^k_{k}, \T^k\) +  N_{\Xi^k_x\times \Xi^k_T }\Big(\big(\ox_{k},\T\big) \Big),
\end{cases}
\end{equation}
where $\H^k$ and $\varrho^k$ are taken from \eqref{e:Hk} and \eqref{e:varrho}, respectively.

\noindent{\sc Complementarity slackness conditions}
\begin{equation}
\label{e:dac20}
[g_i(\ox^k_j-\ou^k_j)>0]\Lto\eta^k_{ji}=0,
\end{equation}
\begin{equation}
\label{e:dac21}
\left\{
\begin{array}{ll}
[i\in I_1(\Lm^k_j)],\; \mbox{i.e.,}\; [g_i(\ox^k_j-\ou^k_j)>0\;\mbox{or} \\
\eta^k_{ji}=0, \la \nabla g_i(\ox^k_j-\ou^k_j),\Lm^k_j)\ra<0]\Lto[\gg^k_{ji}=0],
\end{array}\right.
\end{equation}
\begin{equation}
\label{e:dac22}
\left\{
\begin{array}{ll}
[i\in I_2(\Lm^k_j)],\; \mbox{i.e.,}\; [g_i(\ox^k_j-\ou^k_j)=0,\eta^k_{ji}=0,\;\mbox{and}\\[1ex]
\la \nabla g_i(\ox^k_j-\ou^k_j),\Lm^k_j\ra>0]\Lto[\gg^k_{ji}\ge0]
\end{array}\right.
\end{equation}
for $j=0,\ldots,k-1$ and $i=1,\ldots,m$, where $I_1(\cdot)$ and $I_2(\cdot)$ are defined in \eqref{e:indexes}, together with 
\begin{equation}
\label{e:dac23}
[g_i(\ox^k_j-\ou^k_j)>0]\Lto\gg^k_{ji}=0\;\mbox{for}\;\;j=0,\ldots,k-1\;\;\mbox{and}\;\;i=1,\ldots,m,
\end{equation}
\begin{equation}\label{e:dac24}
[g_i(\ox^k_k-\ou^k_k)>0]\Lto\al^k_{i}=0\;\;\mbox{for}\;\;i=1,\ldots,m,\;\mbox{ and}
\end{equation}
\begin{equation}\label{e:dac25}
\eta^k_{ji}>0\Lto[\la\nabla g_i(\ox^k_j-\ou^k_j),\Lm^k_j)\ra=0].
\end{equation}
Furthermore, the fulfillment of the {\sc Enhanced nontriviality condition}
\begin{equation}
\label{e:dac26}
\lm^k+\sum^{k-1}_{j=0}\n\psi^{uk}_j\en + \| p^{xk}_k\| +\|p^{uk}_0\|\not=0
\end{equation}
is guaranteed under the assumption that the Jacobian matrix $\{\nabla g(\ox^k_j-\ou^k_j)\}$  is surjective.
\end{theorem}\vspace*{-0.1in}
{\bf Proof}. It follows from \eqref{e:dac6} and the coderivative definition \eqref{e:cor} that
\begin{equation}\label{e:dac27}
\begin{array}{ll}
&\bigg(\frac{p^{xk}_{j+1}-p^{xk}_j}{h_k}-\lm^kw^{xk}_j,\frac{p^{uk}_{j+1}-p^{uk}_j}{h_k}-\lm^kw^{uk}_j -\frac{\psi^{uk}_j}{h_k}, -\lm^kw^{ak}_j -\frac{1}{h_k}\lm^k\th^{ak}_j -\frac{\psi^{ak}_j}{h_k}\bigg)\\
&\in D^*F\bigg(\ox^k_j,\ou^k_j,\oa^k_j,\frac{\ox^k_{j+1}-\ox^k_j}{-h_k}\bigg)(\Lm^k_j),\;\;j=0,\ldots,k-1.
\end{array}
\end{equation}
Employing the inclusion $\frac{\ox^k_{j+1}-\ox^k_j}{-h_k}+f(\ox^k_j,\oa^k_j)\in N_{C(u_j}(\ox^k_j-\ou^k_j)$ for $j=0,\ldots,k-1$ along with the representation of $F$ in \eqref{e:F-rep} guarantees the existence of vectors $\eta^k_j\in\R^m_+$, $j=0,\ldots,k-1$ that fulfill the conditions in \eqref{e:dac15} and \eqref{e:dac20}. Using further the second-order upper estimate from Theorem~\ref{Thm_Co} with $x:=\ox^k_j,u:=\ou^k_j,\;a:=\oa^k_j,\;w:=\frac{\ox^k_{j+1}-\ox^k_j}{-h_k}$, and $y:=\Lm^k_j$, and then incorporating \eqref{e:dom-co-est}, we find $\gg^k_j\in\R^m$ such that
 $$
\begin{array}{ll}
&\bigg(\frac{p^{xk}_{j+1}-p^{xk}_j}{h_k}-\lm^kw^{xk}_j,\frac{p^{uk}_{j+1}-p^{uk}_j}{h_k}-\lm^kw^{uk}_j - \frac{1}{h_k}\disp\sum^s_{i=1}\psi^{uk}_{ji}\nabla\nu_i(\ou^k_j),-\lm^kw^{ak}_j -\frac{1}{h_k}\lm^k\th^{ak}_j -\frac{\psi^{ak}_j}{h_k}\bigg)\\
=&\bigg(-\nabla_xf(\ox^k_j,\oa^k_j)^*(\Lm^k_j)- \xi^k_j 
- \disp\sum^m_{i=1}\gg^k_{ji}\nabla g_i(\ox^k_j-\ou^k_j), \xi^k_j+\sum^m_{i=1}\gg^k_{ji}\nabla g_i(\ox^k_j-\ou^k_j), -\nabla_af(\ox^k_j,\oa^k_j)^*(\Lm^k_j)\bigg),
\end{array}
$$
for $j=0,\ldots,k-1$. As a result, all the conditions stated in \eqref{e:dac16}, \eqref{e:dac17}, \eqref{e:dac18}, \eqref{e:dac21}, and \eqref{e:dac22} are satisfied with $\eta^k_j\in\R^m_+$ for $j=0,\ldots,k-1$ and $\al^k\in \R^m_+$. This allows us to derive the nontriviality condition \eqref{e:dac14} from \eqref{e:nontri-da} as well as the transversality conditions \eqref{e:dac19} from \eqref{e:dac4} and \eqref{e:dac7}. Moreover, \eqref{e:dac24} arises naturally from \eqref{e:dac1} and the definition of $\eta^k_k$,whereas \eqref{e:dac27} leads us to the conclusion that
$$
\lm^k(h^{-1}_k\th^{Xk}_j+v^{xk}_j)-p^{xk}_{j+1}\in\dom D^*N_{C(u_j)}\bigg(\ox^k_j-\ou^k_j,\dfrac{\ox^k_{j+1}-\ox^k_j}{-h_k}+f(\ox^k_j,\oa^k_j)\bigg),
$$
It follows from \eqref{e:dom-co-est} that \eqref{e:dac25} holds. The remaining task is to establish the enhanced nontriviality condition \eqref{e:dac26} under the assumption that the Jacobians $\{\nabla g(\ox^k_j-\ou^k_j)\}$ are surjective. Proceeding by contradiction, suppose that \eqref{e:dac26} fails,\ meaning that $\lm^k=0,\; p^{xk}_k=0,\; p^{uk}_0=0$, and $\psi^{uk}_j=0$ for $j=0, \ldots, k-1$. Then we deduce from \eqref{e:dac18a} that $p^{uk}_j=0$ for $j=0, \ldots, k$. Applying the first equation in the transversality condition \eqref{e:dac7} with $p^{uk}_k=0$ results in $\sum\limits^m_{i=1}\al^k_{i}\nabla g_i(\ox^k_k-\ou^k_k)=0$. Due to the surjectivity of $\{\nabla g(\ox^k_k-\ou^k_k)\}$, the latter implies that
$\al^k=(\eta^k_{1}, \ldots, \eta^k_{m})=0$. Furthermore, it follows from \eqref{e:dac17} that
$$
\xi^k_j +\sum^m_{i=1}\gg^k_{ji}\nabla g_i(\ox^k_j-\ou^k_j)=0,\quad j=0,\ldots,k-1,
$$
since $p^{uk}_j=0$ for $j=0, \ldots, k$, $\psi^{uk}_j=0$ for $j=0,\ldots, k-1$, and $\lm^k=0$. Combining the latter with \eqref{e:dac16} and the fact that $p^{xk}_k=0$ brings us to the equalities
$$
p^{xk}_{j+1}-p^{xk}_j=-h_k\nabla_xf(\ox^k_j,\oa^k_j)^*(-p^{xk}_{j+1})\
$$
for $j=0, \ldots, k-1$, and hence $p^{xk}_j=0$ for $j=0,\ldots, k$. It readily follows from the above that $\lm^k=0, \;p^{uk}_0, \;\psi^{uk}_j=0, \;\psi^{ak}_j=0$ for $j=0, \ldots, k-1$, and that $p^{xk}_j=0$ for $j=0,\ldots,k$. The latter contradicts the fulfillment of \eqref{e:dac7} and thus completes the proof of the theorem. $\h$\vspace*{-0.2in}
	
\section{Optimal Control of Motion Models}\label{sec:motion}
\setcounter{equation}{0}\vspace*{-0.1in}

This section is devoted to the application of the obtained well-posedness and optimality conditions to solve sweeping optimal control problems arising in motion modeling. An example of such models is the {\em crowd motion model in a corridor} whose dynamics was described by a sweeping process in \cite{mv0}. A considerably simplified optimal control model for corridor crowd motions was first studied in \cite{cm2}, with polyhedral sweeping sets on the fixed time interval. Based on necessary suboptimality conditions obtained in \cite{cm1}, some relations for calculating optimal controls in the corridor crowd motion model were established in \cite{cm2} in restrictive settings under involved assumptions for the case of two agents.\vspace*{-0.05in}

Having in hand the sweeping control theory developed above, we can do much better to solve {\em realistic moving motion models} described by {\em free-time} sweeping control processes with {\em uniformly prox-regular} moving sets. For two agents moving in a corridor, their interactions during the optimal time interval are carefully investigated by using the obtained (sub)optimality conditions before and after the contact time with calculating optimal control strategies needed to reach the target in the {\em shortest time} with {\em minimal control efforts}. The established relationships allow us to design a {\em numerical algorithm} aimed at calculating optimal parameters of the model, to provide computations for random data, and to reach practical conclusions about the agent's optimal behavior during the dynamic process. \vspace*{-0.05in}

Imagine two agents navigating a corridor, each modeled as a rigid disk in the plane $\R^2$ with radius $L_i$ for $i=1,2$. The position of the $i$-th disk is represented by its center $x_i\in \R^2$. To prevent the agents from overlapping, we introduce  the set of {\em admissible configuration}
\begin{equation}\label{ex1}
C: =\big\{x = (x_1, x_2)\in \R^4\;\big|\;g(x):=\| x_2-x_1\|- (L_1+L_2)\ge 0\big\}.
\end{equation}
Assume that the two agents $x_1$ and $x_2$ are oriented such that the destination is always directly to their right, with $x_2$ being closer to the destination. Under this orientation, their coordinates are given by $x_i=(x_i^1,0)$ for $i=1,2$. Consequently, the set $C$ in \eqref{ex1} can be reformulated as 
\begin{equation*}
C=\big\{(x_1^1, x_2^1)\in \R^2\;\big|\;x_2^1-x_1^1\ge L_1+L_2\big\}
\end{equation*}
Each agent possesses a {\em desired spontaneous velocity}, representing the speed and direction they would naturally adopt if unimpeded by the other. Assuming that both agents seek to reach their destination using the {\em shortest path}, their spontaneous velocity can be described as
$$
U(x_i): =-s_i\frac{x_i-(x^d_i,0)}{\| x_i-(x^d_i,0)\|} = \l( -s_i\frac{x^1_i-x^d_i}{|x^1_i-x^d_i|},0\r)=(s_i,0),
$$
where $0<s_i=\frac{x^d_i-x^1_i}{\T}$ represents the speed of agent $x_i$ for $i=1,2$, and  where $(x^d_i, 0)$ denotes the destination coordinates. When there is no interaction between the agents ($g(x)>0$), each agent moves at the desired velocity meaning that agent's actual velocity match the desired velocity, i.e., $\dot x_i(t)=U(x_i(t))$. However, if the agents come into {\em contact} in the sense that $g(x)=0$, they both must adjust the velocities to {\em avoid collision}. To reflect this scenario, the {\em actual velocities} of two agents should be selected from the following set of {\em admissible velocities} defined by
\begin{equation}\label{ex2}
V(x): =\big\{v = (v^2,0)\in \R^2\;\big|\;g(x)=0\implies \la \nabla g(x), v\ra\geq 0\big\}.
\end{equation}
Suppose that $g(x(t))=0$ and $\dot x(t)\in V(x(t))$. Then the time derivative of $g(x(t))$ satisfies the equation
$$
\frac{dg(x(t))}{dt}=\la \nabla g(x(t)),\dot x(t)\ra\geq 0,
$$
which indicates that selecting the actual velocities in this manner will always increase the distance between the two agents. To establish a connection between the actual and desired velocities, we project the desired velocity onto the set of admissible velocities $V(x)$. By the convexity of $V(x)$ in \eqref{ex2}, this yields 
\begin{equation}\label{ex3}
\dot x(t) = \Pi(U(x(t)); V(x(t)).
\end{equation}
Employing the orthogonal decomposition via the sum of mutually polar cones gives us
$$
U(x(t)) = \Pi(U(x(t)); V(x(t))) + \Pi(U(x(t)); V^0(x(t))) = \dot x(t)+ \Pi(U(x(t)); V^0(x(t))),
$$ 
where $V^0$ stands for the {\em polar} to $V(x)$ defined by $V^0(x):=\{w\;\mid \la w,v\ra\leq 0 \mbox{ for all } v\in V(x)\}$. Then it follows from \cite[Proposition~4.1]{ve} that $V^0 (x) = N_C(x)$. Thus the projected differential equation \eqref{ex3} can be written in the {\em sweeping process form} $\dot x(t)\in -N_C(x(t)) + U(x(t))$.\vspace*{-0.05in}

Now we are in a position to formulate the following optimal control problem $(M)$:
\begin{equation}\label{ex4}
\mbox{ minimize }\;J[x,a,T]:= \l\| x(T)-x^{des}\r\|_{1}+\tau T+\frac{1}{2}\int_0^T\l\|a(t)\r\|^2dt
 \end{equation}
 over the control functions $a(\cdot)=(a_1(\cdot),a_2(\cdot))$ and the corresponding trajectories $x(\cdot)=(x_1^1(\cdot),0,x_2^1(\cdot),0)$ of the nonconvex sweeping process written in the form
 \begin{equation}
 \label{ex5}
 \begin{cases}
 \dot x(t) \in f(x(t),a(t)) - N_{C_{u(t)}}(x(t)) \mbox{ for a.e. } t\in [0,T], \\
 u(t)=0\;\mbox{ for }\;t\in [0,T],\\
 a_i(t)\in A_i:= [\alpha_1, \alpha_2]\;\mbox{ for }\;i=1,2\;\mbox{ and	}\;t\in [0,T],\\
 f(x(t),a(t)): = a(t)U(x(t)) = (a_1(t)s_1, 0, a_2(t)s_2,0), \\
 x(0) = x_0\in C_{u(0)},\\
 (x(T),T) \in \Xi_x\times\Xi_T = \R^n\times [0,\infty),
 \end{cases}
 \end{equation}
 where $C$ is taken from \eqref{ex1}, $\tau$ is a given positive parameter, $\| x(T)-x^{des}\|_{1}:=\l|x_1^1(T)-x^d\r|+\l|x_2^1(T)-x^d\r|$, and where $x^{des}= (x^d,0,x^d,0)$ stands for the destination. The role of control $a$ is to adjust the speed of each agent. The meaning of the cost functional in \eqref{ex4} is driving the agents to the destination as soon as possible with the minimum effort. In this model, the terminal and running costs are given by 
 $$
 \begin{cases}
 \varphi(x(T),T): =\l\| x(T)-x^{des}\r\|_{1}+\tau T,\\
 \ell(t,x,a,\dot x):= \frac{\| a\|^2}{2}.
 \end{cases}
 $$
 We approximate an optimal solution to the original problem by finding optimal solutions to problems $(P_k)$ using a set of necessary optimality conditions obtained in Theorem~\ref{Th9}. For $k$ sufficiently large, the last four terms in the cost functional of $(P^k)$ in \eqref{DP} can be neglected. The gradient of the running cost is calculated by $\nabla_{x,a,\dot x}\ell(t,x,a,\dot x)=(0, a, 0)$.\vspace*{-0.1in}

 To proceed further, fix $k$ and consider an optimal solution $(\ox^k, \oa^k, \T^k)$ $=(\ox^k_0, \ldots, \ox^k_k, \oa^k_0, \ldots, \oa^k_{k-1}, \T^k )$ to problem $(P_k)$ with the uniform grid $\Delta_k(\T^k)$. Applying the necessary optimality conditions from Theorem~\ref{Th9} gives us the dual elements $\lambda^k \geq 0, \; \eta^k_j\in \R_+, \; p^k_j = (p^{1,k}_{1,j}, 0, p^{1,k}_{2,j},0) $, and $\gamma^k_j\in \R$ satisfying the following conditions while dropping the quantities
$(\th^{Xk}_j, \th^{ak}_j) \approx (0,0), \; T-T^k\approx 0,  $ and $\varrho^k\approx 0$:\vspace*{-0.15in}
\begin{enumerate}
\item[\bf(1)] $(w^{xk}_j, w^{ak}_j, v^{xk}_j) = (0,\oa^k_j, 0)$ and $j=0,\ldots, k-1$.\vspace*{-0.05in}
\item[\bf(2)] The discrete velocity of each agent is given by
$$
\begin{cases}
\dot\ox_1^k(t) = \frac{\ox^k_{1,j+1}-\ox^k_{1,j}}{h_k}=(\oa^k_{1,j}s_1,0) + \eta^k_j(-1,0), \\
\dot\ox_2^k(t) = \frac{\ox^k_{2,j+1}-\ox^k_{2,j}}{h_k}=(\oa^k_{2,j}s_2,0) + \eta^k_j(1,0) 
\end{cases}
$$
for $j=0,\ldots, k-1$ and $t\in [t^k_j, t^k_{j+1})$.\vspace*{-0.05in}
\item[\bf(3)] $\ox^{1,k}_{2,j}-\ox^{1,k}_{1,j}>L_1+L_2 \implies \eta^k_j=0 $\;\mbox{ for }\;$j=0,\ldots, k-1$.\vspace*{-0.05in}
\item[\bf(4)] $\eta^k_j>0 \implies 0=
\la (-1,0,1,0),-p^{xk}_{j+1}\ra 
= -\la (-1,0,1,0), (p^{1,xk}_{1,j+1},0,p^{1,xk}_{2,j+1},0)\ra $, 
i.e., \\
$\big[\eta^k_j>0 \implies p^{1,xk}_{1,j+1}= p^{1,xk}_{2,j+1}\big]$\;\mbox{ for }\;$j=0,\ldots, k-1$.\vspace*{-0.05in}
\item[\bf(5)] 
$
\begin{cases}
\Lm^k_j = -p^{xk}_{j+1}, \; \xi^k_j=0,
\frac{p^{xk}_{j+1}-p^{xk}_j}{h_k}= -\gamma^k_j(-1,0, 1,0) 
\end{cases}
$ for $j=0,\ldots, k-1$.\vspace*{-0.05in}
\item[\bf(6)] 
$
\begin{cases}
\lm^k \oa^k_j+\frac{\psi^{ak}_j}{h_k}=\(s_1p^{1,xk}_{1,j+1}, s_2p^{1,xk}_{2,j+1}\),\\
\psi^{ak}_j\in N_A(\oa^k_j)
\end{cases}
$ for $j=0,\ldots, k-1$.\vspace*{-0.05in}
\item[\bf(7)] 
$
\begin{cases}
\H^k = \frac{1}{k}\sum\limits^{k-1}_{j=0}\[\la p^{xk}_{j+1}, \frac{\ox^k_{j+1}-\ox^k_j}{h_k}\ra  -\frac{1}{2}\lm^k\|\oa^k_j\|^2 \],\\
\nabla\varphi(\ox^k_k,\T^k)= (-1,0,-1,0, \tau), \;
N_{\Xi^k_x\times \Xi^k_T}((\ox^k_k, \T^k))=\{(0,0)\}, \\
p^{xk}_k + \al^k \(-1,0,1,0\), \H^k) = \lm^k\nabla\varphi(\ox^k_k,\T^k) = \lm^k(-1,0,-1,0, \tau ).
\end{cases}
$

which readily implies that \\
\begin{equation}
\label{ex6}
\begin{cases}
p^{xk}_k =  \al^k \(-1,0,1,0\) + \lm^k(1,0,1,0), \\
\H^k =\lm^k\tau.
\end{cases}
\end{equation}
\item[\bf(9)] $\lm^k  + \| p^{xk}_k\| >0$.
\end{enumerate}\vspace*{-0.1in}
\begin{figure}
\centering
\includegraphics[width=0.7\linewidth]{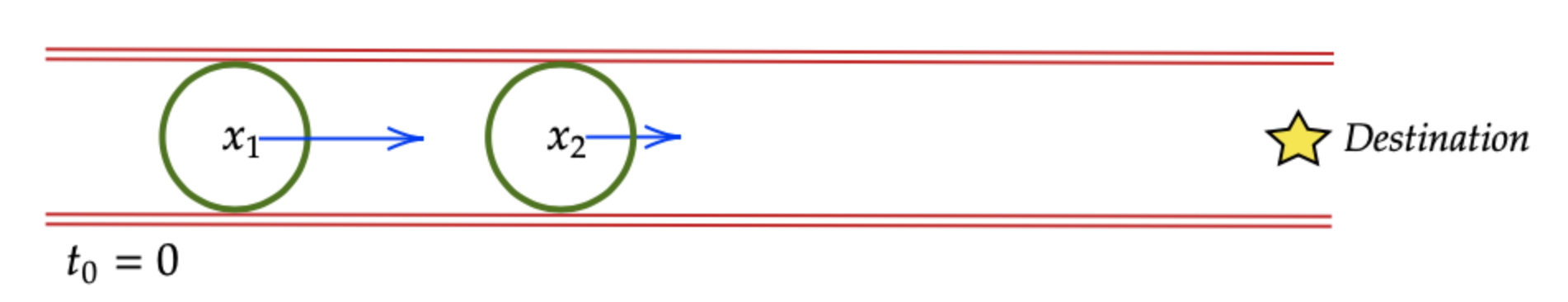}\vspace*{-0.15in}
\caption{Crowd motion model in a corridor before the contacting time.}
\label{f1}
\end{figure}
In what follows, we analyze the obtained optimality conditions ${\bf(1)}$--${\bf(9)}$ and deduce some conclusions from them about the agents' behavior before and after the contact time. Observe first that {\bf(6)} yields
\begin{equation}
\label{ex7}
\begin{array}{ll}
\(p^{1,xk}_{1,j+1}, p^{1,xk}_{2,j+1}\) =
\(\frac{\lm^k\oa^k_{1,j}}{s_1} + \frac{\psi^{ak}_{1,j}}{s_1h_k}, \frac{\lm^k\oa^k_{2,j}}{s_2}+\frac{\psi^{ak}_{2,j}}{s_2h_k}\) 
\end{array}
\end{equation}
for $j=0,\ldots, k-1$. Using the equations from {\bf(5)} tells us that
$$
\begin{cases}
\frac{p^{1,xk}_{1,j+1}-p^{1,xk}_{1,j}}{h_k} = \gg^k_j,\\
\frac{p^{1,xk}_{2,j+1}-p^{1,xk}_{2,j}}{h_k}  = -\gg^k_j,
\end{cases}
$$
which implies in turn that
\begin{equation*}
 \begin{array}{ll}
 p^{1,xk}_{1,j+1}+p^{1,xk}_{2,j+1}= p^{1,xk}_{1,j}+p^{1,xk}_{2,j}, 
 \end{array}
 \end{equation*}
 and thus ensures the representation
\begin{equation}
\label{ex8}
\begin{array}{ll}
\frac{\lm^k}{s_1}\oa^k_{1,j+1}  +\frac{\psi^{ak}_{1,j+1}}{s_1h_k}+\frac{\lm^k}{s_2}\oa^k_{2,j+1}  +\frac{\psi^{ak}_{2,j+1}}{s_2h_k}= \frac{\lm^k}{s_1}\oa^k_{1,j}  +\frac{\psi^{ak}_{1,j}}{s_1h_k}+\frac{\lm^k}{s_2}\oa^k_{2,j}  +\frac{\psi^{ak}_{2,j}}{s_2h_k}
\end{array}
\end{equation} 
for $j=0, \ldots, k-2$  due to \eqref{ex7}. It follows from {\bf (6)} and {\bf(7)} that 
\begin{equation}
\label{ex9}
\begin{cases}
\frac{\lm^k}{s_1}\oa^k_{1,k-1} + \frac{\psi^{ak}_{1,k-1}}{s_1h_k}  -\lm^k = -\al^k, \\
\frac{\lm^k}{s_2}\oa^k_{2,k-1} + \frac{\psi^{ak}_{2,k-1}}{s_2h_k} -\lm^k=\al^k.
\end{cases}
\end{equation} 
For simplicity, consider the case where $\oa^k_{i,j}\in A_i=(-\infty, \infty)$. Then all the quantities $\psi^{ak}_j$ vanish in \eqref{ex7} and \eqref{ex8}, and hence  we can rewrite \eqref{ex8} in the form
\begin{equation}
\label{ex10}
\begin{array}{ll}
\frac{\lm^k}{s_1}\oa^k_{1,j+1} +\frac{\lm^k}{s_2}\oa^k_{2,j+1} = \frac{\lm^k}{s_1}\oa^k_{1,j} +\frac{\lm^k}{s_2}\oa^k_{2,j}  
\end{array}
\end{equation}
for $j=0, \ldots, k-2$.
Let $t=t^\ast \in \[t^k_{j_0}, t^k_{j_0+1}\)$ be the first time that two agents are in contact, i.e., $\ox^{1,k}_2(t^\ast)-\ox^{1,k}_1(t^\ast) = L_1+L_2$ for some $j_0$. Again for simplicity, we assume that $t^\ast=t^k_{j_0}$ and investigate now the behavior of both agents before the contact time $t^\ast$. \vspace*{-0.05in}

{\bf Agents' behavior before the contact time} (see Figure~\ref{f1}): 
In the time interval $t\in [0,t^\ast)$, the velocities of the two agents are represented by 
\begin{equation*}
\begin{cases}
\dot\ox_1^k(t) = (\oa^k_1(t)s_1,0),\\
\dot\ox_2^k(t) = (\oa^k_2(t)s_2,0)
\end{cases}
\end{equation*}
due to $\eta^k_j=0$ for $j=0,\ldots, j_0-1$.
Using the equations from {\bf (5)} with $\gg^k_j=0$, we get 
\begin{equation}
\label{ex11}
\begin{cases}
\frac{\lm^k}{s_1}\oa^k_{1,j+1} =\frac{\lm^k}{s_1}\oa^k_{1,j}\\[1ex]
\frac{\lm^k}{s_2}\oa^k_{2,j+1} =\frac{\lm^k}{s_2}\oa^k_{2,j}\end{cases}
\mbox{ for }\;j=0, \ldots, j_0-2.
\end{equation}
The nontriviality condition {\bf (9)} ensures that $\lm^k\neq 0$. Then it follows from \eqref{ex11} that 
\begin{equation*}
\begin{cases}
\oa^k_{1,j+1} =\oa^k_{1,j},\\
\oa^k_{2,j+1} =\oa^k_{2,j}
\end{cases}
\mbox{ for }\;j=0, \ldots, j_0-2, 
\end{equation*}
which tells us therefore that 
\begin{equation}
\label{ex12}
\begin{cases}
\oa^k_{1,j}=\oa^k_{1,0}, \\
\oa^k_{2,j}= \oa^k_{2,0}
\end{cases}
\mbox{ for }\;j=0, \ldots, j_0-1. 
\end{equation}
Using the discrete velocity formulas in {\bf(2)} yields 
\begin{equation}
\label{ex13}
\begin{cases}
\ox^{1,k}_{1,j_0}= \ox^{1,k}_{1,0}+h_k\sum\limits^{j_0-1}_{j=0}\oa^k_{1,j}s_1 = x^d-s_1\T^k + h_kj_0\oa^k_{1,0}s_1,\\
\ox^{1,k}_{2,j_0}= \ox^{1,k}_{2,0}+h_k\sum\limits^{j_0-1}_{j=0}\oa^k_{2,j}s_2 = x^d-s_2\T^k + h_kj_0\oa^k_{2,0}s_2,
\end{cases}
\end{equation}
which being combined with the fact that $\ox^{1,k}_{2,j_0} - \ox^{1,k}_{1,j_0}=L_1+L_2 $ gives us
\begin{equation*}
\ox^{1,k}_{2,0}-\ox^{1,k}_{1,0} + h_kj_0\underbrace{\( \oa^k_{2,0}s_2-\oa^k_{1,0}s_1\)}_{\leq0} = \ox^{1,k}_{2,j_0} - \ox^{1,k}_{1,j_0}=L_1+L_2.
\end{equation*}
The latter can be equivalently rewrite as 
\begin{equation*}
(s_1-s_2)\T^k+ h_kj_0\underbrace{\( \oa^k_{2,0}s_2-\oa^k_{1,0}s_1\)}_{\leq0} = \ox^{1,k}_{2,j_0} - \ox^{1,k}_{1,j_0}=L_1+L_2,
\end{equation*}
which implies in turn that 
\begin{equation*}
\begin{array}{ll}
\Lm_1-\Lm_2&=(s_1-s_2)\T^k= h_kj_0\(\oa^k_{1,0}s_1-\oa^k_{2,0}s_2\)+L_1+L_2\\
&=\frac{1}{k}j_0\(\oa^k_{1,0}\T^ks_1-\oa^k_{2,0}\T^ks_2\)+L_1+L_2\\
&=\frac{1}{k}j_0\(\oa^k_{1,0}\Lm_1-\oa^k_{2,0}\Lm_2\)+L_1+L_2
\end{array}
\end{equation*}
with the notation  $\Lm_i:=x^d-\ox^{1,k}_{i,0}=s_i\T^k$ for $i=1, 2$.

{\bf Agents' behavior after the contact time} (see Figure~\ref{f2}): 
\begin{figure}
\centering
\includegraphics[width=0.7\linewidth]{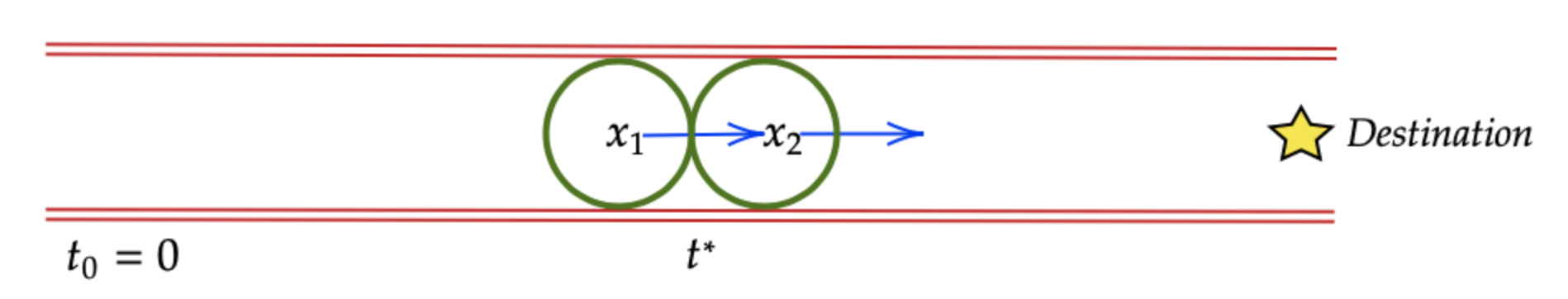}\vspace*{-0.15in}
\caption{Crowd motion model in a corridor after the contacting time.}
\label{f2}
\end{figure}
Next we investigate the behavior of both agents on the time interval $[t^k_{j_0},\T^k]$ after the contact time. The velocities of the two agents in this case are  
\begin{equation}
\label{ex14}
\begin{cases}
\dot\ox^k_1(t)=(\oa^k_{1,j}s_1,0) + \eta^k_{j_0},(-1,0), \\
\dot\ox^k_2(t)=(\oa^k_{2,j}s_2,0) + \eta^k_{j_0}(1,0)
\end{cases}
\end{equation}
for $t\in [t^k_{j_0}, \T^k]$. We have $\ox^{1,k}_{2}(t)-\ox^{1,k}_1(t) =L_2+L_1$, and thus $\dot\ox^{1,k}_2(t)-\dot\ox_1^{1,k}(t) =  0$ for $t\in [t^k_{j_0}, \T^k]$ to avoid collision of the agents. On the other hand, it follows from \eqref{ex9} that 
\begin{equation}
\label{ex15}
\begin{array}{ll}
\frac{\oa^k_{1,k-1}}{s_1}+\frac{\oa^k_{2,k-1}}{s_2}=2.
\end{array}
\end{equation}
Let us show that the numbers $\eta^k_j$ always take positive values after the contact time $[t^\ast, \T^k]$. Indeed, we can calculate the energy $\(\oa^k_{1,j}\)^2+\( \oa^k_{2,j}\)^2 $ explicitly in terms of $\oa^k_{1,j}$ by
\begin{equation*}
\(\oa^k_{1,j}\)^2+\( \oa^k_{2,j}\)^2 = 
\begin{cases}
\(\frac{s_1^2+s_2^2}{s_1^2}\)\(\oa^k_{1,j}\)^2 \mbox{ if } \eta^k_j>0, \\
\(\frac{s_1^2+s_2^2}{s_2^2}\)\(\oa^k_{1,j}\)^2 \mbox{ if } \eta^k_j=0. 
\end{cases}
\end{equation*}
It is obvious that the energy in the first case is less than the one in the second case since $s_1>s_2$. Therefore, we can rule out the case $\eta^k_j=0$ on the interval after the contact time, and hence 
\begin{equation}
\label{ex16}
\begin{array}{ll}
\frac{\oa^k_{1,j}}{s_1}=\frac{\oa^k_{2,j}}{s_2}\;\mbox{ or }\;\frac{\oa^k_{1,j}}{\Lm_1}=\frac{\oa^k_{2,j}}{\Lm_2}\;\mbox{ for }\;j=j_0, \ldots, k-1.
\end{array}
\end{equation}
This reflects that the optimal sweeping motion reaches the boundary of the state constraints and remains there until the process concludes; see Figure~\ref{f2}.
Furthermore, we deduce from \eqref{ex15} that 
\begin{equation}
\label{ex17}
\begin{array}{ll}
\frac{\oa^k_{1,k-1}}{s_1} = \frac{\oa^k_{2,k-1}}{s_2}=1.
\end{array}
\end{equation}
Using \eqref{ex10}, \eqref{ex16}, and \eqref{ex17} gives us the equality
\begin{equation}
\label{ex18}
\begin{array}{ll}
\frac{\oa^k_{1,j}}{s_1}=\frac{\oa^k_{2,j}}{s_2}=1\;\mbox{ for }\;j = j_0, \ldots, k-1.
\end{array}
\end{equation}
Then it follows from \eqref{ex18}, \eqref{ex10}, and \eqref{ex12} that 
\begin{equation}
\label{ex18a}
\begin{array}{ll}
\frac{\oa^k_{1,0}}{s_1} + \frac{\oa^k_{2,0}}{s_2}=2.
\end{array}
\end{equation}
Using again the discrete velocity equations, we get
\begin{equation*}
\begin{cases}
\ox^{1,k}_{1,k}&= \ox^{1,k}_{1,j_0}+h_k\frac{s^2_1+s^2_2}{2}\sum\limits^{k-1}_{j=j_0}\frac{\oa^k_{1,j}}{s_1} = \ox^{1,k}_{1,j_0}+ \frac{\T^k}{k}(k-j_0)\frac{s^2_1+s^2_2}{2}\\
&=\ox^{1,k}_{1,j_0}+ \frac{k-j_0}{k}\frac{\Lm^2_1+\Lm^2_2}{2\T^k},
\\
\ox^{1,k}_{2,k}&= \ox^{1,k}_{2,j_0}+h_k\frac{s^2_1+s^2_2}{2}\sum\limits^{k-1}_{j=j_0}\frac{\oa^k_{2,j}}{s_2} = \ox^{1,k}_{2,j_0}+(k-j_0)\frac{s^2_1+s^2_2}{2}\\
&=\ox^{1,k}_{2,j_0}+ \frac{k-j_0}{k}\frac{\Lm^2_1+\Lm^2_2}{2\T^k}.
\end{cases}
\end{equation*}
Combining this with \eqref{ex13} and \eqref{ex14} brings us to 
\begin{equation}
\label{ex19}
\begin{cases}
\ox^{1,k}_{1,k}-x^d&=-\Lm_1+h_kj_0\oa^k_{1,0}s_1 +\frac{k-j_0}{k}\frac{\Lm^2_1+\Lm^2_2}{2\T^k}\\
&= -\Lm_1+\frac{1}{k}j_0\oa^k_{1,0}\Lm_1 +\frac{k-j_0}{k}\frac{\Lm^2_1+\Lm^2_2}{2\T^k},\\
\ox^{1,k}_{2,k}-x^d&=-\Lm_2+h_kj_0\oa^k_{2,0}s_2 +\frac{k-j_0}{k}\frac{\Lm^2_1+\Lm^2_2}{2\T^k}\\
&= -\Lm_2+\frac{1}{k}j_0\oa^k_{2,0}\Lm_2 +\frac{k-j_0}{k}\frac{\Lm^2_1+\Lm^2_2}{2\T^k}.
\end{cases}
\end{equation}\vspace*{-0.1in}

Now we are able to calculate the {\em optimal time}, {\em optimal cost function value}, {\em optimal contact time}, and parameter of {\em optimal controls}. First compute $\H^k$ from {\bf(7)} by
\begin{equation*}
\begin{array}{ll}
\H^k &= \frac{1}{k}\sum\limits^{k-1}_{j=0}\[\la p^{xk}_{j+1}, \frac{\ox^k_{j+1}-\ox^k_j}{h_k}\ra  -\frac{\lm^k}{2}\sum\limits^{k-1}_{j=j_0}\|\oa^k_j\|^2\]\\
&= \frac{\lm^k}{2k}\[\sum\limits^{j_0-1}_{j=0}\l\{\(\oa^k_{1,j}\)^2+\(\oa^k_{2,j}\)^2\r\} +\frac{s^2_1+s^2_2}{s_1^2}\sum\limits^{k-1}_{j=j_0}\(\oa^k_{1,j}\)^2\]\\[1ex]
&= \frac{\lm^k}{2 k}\[\sum\limits^{j_0-1}_{j=0}\l\{\(\oa^k_{1,j}\)^2+\(\oa^k_{2,j}\)^2\r\} +\frac{\Lm_1^2+\Lm_2^2}{\Lm_1^2}\sum\limits^{k-1}_{j=j_0}\(\oa^k_{1,j}\)^2\]\\
&= \frac{\lm^k}{k}\[j_0\frac{\(\oa^k_{1,0}\)^2+\(\oa^k_{2,0}\)^2}{2}  +\frac{\Lm_1^2+\Lm_2^2}{2\(\T^k\)^2}(k-j_0)\].
\end{array}
\end{equation*}
Combining this with \eqref{ex6} gives us the formula  
\begin{equation}
\label{ex21}
\begin{array}{ll}
\tau = \frac{j_0}{k}\frac{\(\oa^k_{1,0}\)^2+\(\oa^k_{2,0}\)^2}{2}  +\frac{k-j_0}{k}\frac{\Lm_1^2+\Lm_2^2}{2\(\T^k\)^2}.
\end{array}
\end{equation}
It follows from the computation of $\H^k$ that $\lm^k\tau =\H^k=\frac{\lm^k}{k}\sum\limits^{k-1}_{j=0} \ell\(t^k_j, x^k_j,u^k_j,a^k_j, \frac{x^k_{j+1}-x^k_j}{h_k}, \frac{u^k_{j+1}-u^k_j}{h_k}\)$, i.e., 
\begin{equation*}
\begin{array}{ll}
\sum\limits^{k-1}_{j=0} h_k\ell\(t^k_j, x^k_j,u^k_j,a^k_j, \frac{x^k_{j+1}-x^k_j}{h_k}, \frac{u^k_{j+1}-u^k_j}{h_k}\) = kh_k\tau= \tau\T^k.
\end{array}
\end{equation*}
This allows to compute the cost function $ J_k[x^k,u^k,a^k, T^k]$ by
\begin{equation}
\label{ex22}
\begin{array}{ll}
J_k[\ox^k,\oa^k, \T^k] &= \ph\(\ox^k_{k},\T^k\) + \sum\limits^{k-1}_{j=0} h_k\ell\(t^k_j, \ox^k_j,\ou^k_j,\oa^k_j, \frac{\ox^k_{j+1}-\ox^k_j}{h_k}\) \\
&=\l|\ox^{1,k}_{1,k}-x^d \r|+ \l|\ox^{1,k}_{2,k}-x^d\r|\ + 2\tau\T^k \\
&= \Lm_1- \frac{1}{k}j_0\oa^k_{1,0}\Lm_1 -\frac{k-j_0}{k}\frac{\Lm^2_1+\Lm^2_2}{2\T^k} 
+ \Lm_2-\frac{1}{k}j_0\oa^k_{2,0}\Lm_2 -\frac{k-j_0}{k}\frac{\Lm^2_1+\Lm^2_2}{2\T^k} +2\tau\T^k \\
&=\Lm_1+\Lm_2 -\frac{1}{k}j_0\(\oa^k_{1,0}\Lm_1+\oa^k_{2,0}\Lm_2\)  -\frac{k-j_0}{k}\frac{\Lm^2_1+\Lm^2_2}{\T^k}+ 2\tau\T^k.
\end{array}
\end{equation}
Let $\zeta\in (-1,1)$ be a scalar such that 
\begin{equation*}
\begin{array}{ll}
J_k[\ox^k,\oa^k, \T^k] &= \Lm_1+\Lm_2 -\frac{1}{k}j_0\(\oa^k_{1,0}\Lm_1+\oa^k_{2,0}\Lm_2\)  -\frac{k-j_0}{k}\frac{\Lm^2_1+\Lm^2_2}{\T^k}+ \tau\T^k +\tau \T^k \\
&=\Lm_1+\Lm_2+\frac{j_0}{k}\[\frac{\T^k\(\oa^k_{1,0}\)^2}{2}-\Lm_1\oa^k_{1,0}\] + \frac{j_0}{k}\[\frac{\T^k\(\oa^k_{2,0}\)^2}{2} - \Lm_2\oa^k_{2,0}\] - \frac{k-j_0}{k}\frac{\Lm_1^2+\Lm_2^2}{2\T^k}+\tau\T^k\\
&=\Lm_1+\Lm_2+\frac{j_0}{k}\frac{\T^k}{2}\(\oa^k_{1,0}-\frac{\Lm_1}{\T^k}\)^2+ \frac{j_0}{k}\frac{\T^k}{2}\(\oa^k_{2,0}-\frac{\Lm_2}{\T^k}\)^2 - \frac{j_0}{k}\frac{\Lm_1^2+\Lm_2^2}{2\T^k } - \frac{k-j_0}{k}\frac{\Lm_1^2+\Lm_2^2}{2\T^k}+\tau\T^k\\
&= \Lm_1+\Lm_2  +\frac{j_0}{k}\frac{\T^k}{2}\zeta^2 s_1^2+\frac{j_0}{k}\frac{\T^k}{2}\zeta^2 s_2^2  - \frac{\Lm_1^2+\Lm_2^2}{2\T^k } + \tau\T^k\\
&=  \Lm_1+\Lm_2  +\frac{j_0}{k}\frac{\zeta^2}{2\T^k}\(\Lm_1^2+\Lm_2^2\) - \frac{\Lm_1^2+\Lm_2^2}{2\T^k } +\tau \T^k\\
&=\Lm_1+\Lm_2 + \frac{\Lm_1^2+\Lm_2^2}{2\T^k }\[\underbrace{\frac{j_0}{k}\zeta^2-1}_{<0}\] + \tau\T^k.
\end{array}
\end{equation*}
Employing \eqref{ex18a} brings us to 
\begin{equation*}
\begin{array}{ll}
\tau=\frac{j_0}{k}\frac{(1+\zeta)^2s_1^2+(1-\zeta)^2s_2^2}{2} + \frac{k-j_0}{k}\frac{\Lm_1^2+\Lm_2^2}{2\(\T^k\)^2}
= \frac{j_0}{k}\frac{(1+\zeta)^2\Lm_1^2+(1-\zeta)^2\Lm_2^2}{2\(\T^k\)^2} + \frac{k-j_0}{k}\frac{\Lm_1^2+\Lm_2^2}{2\(\T^k\)^2}.
\end{array}
\end{equation*}
This yields therefore the calculation formula
\begin{equation}
\label{ex23}
\begin{array}{ll}
\(\T^k\)^2 = \frac{1}{\tau}\[\frac{j_0}{k}\frac{(1+\zeta)^2\Lm_1^2+(1-\zeta)^2\Lm_2^2}{2} +\frac{k-j_0}{k}\frac{\Lm_1^2+\Lm_2^2}{2}\]: = h(\zeta).
\end{array}
\end{equation}
Consider the case where agent~1 is forced to move faster and agent!2 should slow down on the way to the destination, which implies that $\zeta\geq 0$. It allows us to deduce from \eqref{ex23} that $\(\T^k\)^2=h(\zeta)\geq h(0)=\frac{\Lm_1^2+\Lm_2^2}{2\tau}$. As a consequence, we get the inequality
$$
\begin{array}{ll}
J_k\[\ox^k,\oa^k, \T^k\]\geq \Lm_1+\Lm_2-\frac{\Lm_1^2+\Lm_2^2}{2\sqrt{h(0)}}+\sqrt{h(0)}. 
\end{array}
$$
where the equality holds if $\zeta=0$. Hence the values of the optimal time and cost functional are $\T^k=\sqrt{\frac{\Lm_1^2+\Lm_2^2}{2\tau}}$ and $J_k\[\ox^k,\oa^k,\T^k\]=\Lm_1+\Lm_2$, respectively. The optimal controls are $\oa^k_{i,j}=s_i=\frac{\Lm_i}{\T^k} =\frac{\sqrt{2\tau}\Lm_i}{\sqrt{\Lm_1^2+\Lm_2^2}}$ for $i=1,2$ and $j=0,\ldots, k-1$, while the contact time is computed by 
$$
\begin{array}{ll}
t^k_{j_0}=h_kj_0=\frac{\(\T^k\)^2\Lm}{\Lm_1^2-\Lm_2^2}=\frac{\(\Lm_1^2+\Lm_2^2\)\Lm}{2\tau\(\Lm_1^2-\Lm_2^2\)}
\end{array}
$$
with $\Lm:=\Lm_1-\Lm_2 - (L_1+L_2)$.\vspace*{-0.05in}

To further demonstrate the application of the necessary conditions from Theorem~\ref{Th9} to this model, we not only analyze, compute, and determine optimal solutions based on the dynamical system, but also develop a {\em numerical algorithm} in Python; see Algorithm~1 below. Specifically, consider the data $x^d = 0$, $x^1_0 = (-48, 0)$, $x^2_0 = (-24, 0)$, and $L_1 = L_2 = 3$ when executing the code.

\begin{algorithm}
\caption{Discrete dynamic optimization of two controlled agents}
\textbf{Input:} $\tau$, $x_d$, $x_{1\text{int}}$, $x_{2\text{int}}$, $L_1$, $L_2$ \\
\textbf{Output:} \\
\hspace*{10pt} $\triangleright$ Approximating optimal time value $\T^k$\\
\hspace*{10pt} $\triangleright$ Approximating optimal controls $\oa^k_1, \oa^k_2$\\ 
\hspace*{10pt} $\triangleright$ Approximating contact time value $t^\ast$
\begin{algorithmic}[1]
\State $\Lm_1 \gets x_d - x_{1\text{int}}$
\State $\Lm_2 \gets x_d - x_{2\text{int}}$
\State $\Lm\gets \Lm_1-\Lm_2 - L_1-L_2$
\State $\T^k \gets \sqrt{\frac{\Lm_1^2 + \Lm_2^2}{2\tau}}$
\State $\oa^k_1 \gets \frac{\Lm_1}{\T^k}$
\State $\oa^k_2 \gets \frac{\Lm_2}{\T^k}$
\State $t^\ast \gets \frac{\(\T^k\)^2 \Lm}{\Lm_1^2 - \Lm_2^2}$ 
\State \textbf{return} $\T^k, \oa^k_1, \oa^k_2, t^\ast$
\end{algorithmic}
\end{algorithm}\vspace*{-0.07in}

In this scenario, agent~1 and agent~2 are initially positioned 48 and 24 meters away from the destination, respectively, while maintaining a minimum safe distance of $L_1 + L_2 = 6$ meters between them. Consequently, agent~1 is inclined to move faster than agent~2. Setting $\tau = 1$, we calculate the optimal time and controls as $T^k = 37.94$ seconds and $(\oa^k_1, \oa^k_2) = (1.26, 0.63)$; see Table~\ref{Ta1}. The uncontrolled speeds of the two agents are $s_1 = \oa^k_1 = 1.26$ and $s_2 = \oa^k_2 = 0.63$. To enhance their performance, agent~1 accelerates while agent~2 decelerates. These changes are reflected in their controlled speeds before and after contact being calculated as $s^b_1 = \oa^k_1 s_1 = 1.6$, $s^b_2 = \oa^k_2 s_2 = 0.4$, and $s = \frac{s^b_1 + s^b_2}{2} = \frac{\oa^k_1 s_1 + \oa^k_2 s_2}{2} = \tau = 1$. Here $s^b_1$, $s^b_2$, and $s$ represent the speeds of the two agents before and after contact. To reach the destination faster, both agents exert additional control efforts, balancing the trade-off between time efficiency and energy expenditure. The optimal times, controls, and related data for various $\tau$ values are presented in Table~\ref{Ta1}. \vspace*{-0.05in}

\begin{table}[h!]
\centering
\caption{Optimal time, controls, and associated parameters for different values of $\tau$ in the controlled motion problem.}

\label{Ta1}
\begin{tabular}{|c|c|c|c|c|c|c|c|}
\hline
\rowcolor{blue!20} $\tau $ & $\oa^k_1$ & $\oa^k_2$ & $s^b_1$ & $s^b_2$ & $s$ & $\T^k$ & $t^\ast$\\
\hline
\hline
\rowcolor{blue!10} 1  & 1.26 & 0.63 & 1.6 & 0.4 & 1 & 37.94 & 15 \\
\hline
2 & 1.78 & 0.89 & 3.2 & 0.8 & 2 & 26.83 & 7.5 \\
\hline
\rowcolor{blue!10} 3 & 2.19 & 1.09 & 4.8 & 1.2 & 3 & 21.9 & 5 \\
\hline
4 & 2.52 & 1.26 & 6.4 & 1.6 & 4 & 18.97 & 3.75 \\
\hline
\rowcolor{blue!10} 5 & 2.82 & 1.41 & 8 & 2 & 5 & 16.97 & 3\\
\hline
6 & 3.09 & 1.54 & 9.6 & 2.4 & 6 & 15.49 & 2.5\\
\hline
\rowcolor{blue!10} 7 & 3.34 & 1.67 & 11.2 & 2.8 & 7 & 14.34 & 2.14 \\
\hline
8 &  3.57 & 1.78 & 12.8 & 3.2 & 8 & 13.41 & 1.875 \\
\hline		
\rowcolor{blue!10} 9 & 3.79 & 1.89 & 14.4 & 3.6 & 9 & 12.64 & 1.66 \\
\hline
10 & 4 & 2 & 16 & 4 & 10 & 12 & 1.5 \\
\hline
\end{tabular} 
\end{table} 
 In this model, the two agents start from different initial positions, with agent 2 closer to the target. Consequently, agent 1, being farther from the target, requires more energy than agent 2 to reach the destination. Although the rigid disks representing the agents have identical radii, agent 1 moves at a speed four times greater than that of agent 2 before contact occurs. However, once contact is made, both agents synchronize their speeds and proceed toward the target together.\vspace*{-0.05in}

Observe that the optimal time to reach the target increases if the agents have less energy or lower initial speed. This finding underscores the importance of energy and speed in determining the efficiency of their motions toward the target. Any reduction in these parameters can significantly prolong the time required to complete the task. The synchronization of speeds after contact indicates a {\em cooperative dynamic}, where both agents adjust their velocities to achieve the goal. The  considerations above emphasize the interplay between energy, speed, and position in optimizing the agents' motions, highlighting the crucial role which these factors play in controlling the overall process.\vspace*{-0.25in}

\section{Conclusions and Future Research}\label{sec:conc}\vspace*{-0.1in}

This paper introduces a novel approach to tackling optimal control problems involving discontinuous dynamic constraints, with a particular emphasis on the timing of dynamic processes. The developed new version of discrete approximation method has been demonstrated to be both reliable and effective, showing strong convergence toward optimal solutions. Through the application of advanced variational analysis techniques, we derived new necessary optimality conditions in discrete-time problems with the nonconvex dynamics and free-time constraints. These theoretical findings are further illustrated in applications to motion models, where we are able to calculate optimal parameters of the control model.\vspace*{-0.05in}

Our future research aims at deriving necessary optimality conditions for the continuous-time sweeping control problem $(P)$ with free time. It can be done by passing to the limit from the necessary optimality conditions for discrete-time problems $(P_k)$ developed in Section~\ref{nec}. Then we plan to apply the obtained conditions to various practical models governed by the constrained sweeping dynamics with free time. Besides motion models (in both corridor and planar versions). this includes more realistic models of robotics, traffic equilibria, unmanned surface vehicles, and nanotechnology in comparison with their simplified polyhedral versions studied recently in \cite{CMNNO2024,mnn}. 
\end{document}